\documentclass[11pt,reqno,final]{amsart}
\usepackage{subeqnarray}
\usepackage{cases}
\usepackage{stmaryrd}
\usepackage{color}
\usepackage{amsfonts}
\usepackage{mathrsfs}
\usepackage{amssymb,amsmath,amsthm}
\usepackage{epsfig}
\usepackage{graphicx}
\usepackage{appendix}
\usepackage{lineno}
\usepackage{enumerate}
\usepackage{comment}
\usepackage{caption}
\usepackage[notcite,notref]{showkeys}
\usepackage[numbers,sort&compress]{natbib}
\usepackage{empheq}
\newcommand{\dps}{\displaystyle}
\newcommand{\tps}{\textstyle}
\newtheorem{theorem}{\indent Theorem}[section]
\newtheorem{lemma}{\indent Lemma}[section]

\newtheorem{remark}{\indent Remark}[section]

\newcommand{\ba}{\begin{array}}\newcommand{\ea}{\end{array}}
\newcommand{\be}{\begin{eqnarray}}\newcommand{\ee}{\end{eqnarray}}
\newcommand{\beq}{\begin{equation*}}\newcommand{\eeq}{\end{equation*}}
\newcommand{\bex}{\begin{eqnarray*}}
\newcommand{\eex}{\end{eqnarray*}}

\def\bq{\begin{equation}}
\def\eq{\end{equation}}
\def\beq{\begin{equation*}}
\def\eeq{\end{equation*}}
\def\br{\begin{eqnarray}}
\def\er{\end{eqnarray}}
\def\brr{\bq\begin{array}{r@{}l}}
\def\err{\end{array}\eq}
\def\bry{\beq\begin{array}{r@{}l}}
\def\ery{\end{array}\eeq}
\def\brl{\bq\begin{array}{l}}
\def\erl{\end{array}\eq}
\def\bryl{\beq\begin{array}{l}}
\def\eryl{\end{array}\eeq}

\font\tenbi=cmmib10   at 11 pt
\font\sevenbi=cmmib10 at 9pt
\font\fivebi=cmmib7 at 6pt
\newfam\bifam
\textfont\bifam=\tenbi \scriptfont\bifam=\sevenbi  \scriptscriptfont\bifam=\fivebi

\def\bi{\fam\bifam\tenbi}
\font\sixtdb=msbm10 at 16 pt \font\tendb=msbm10 at 12 pt  \font\sevendb=msbm7
\newfam\dbfam
\textfont\dbfam=\sixtdb

\textfont\dbfam=\tendb \scriptfont\dbfam=\sevendb

\def\Dt {\triangle t}

\def\n{{\bi n}}
\def\x{{\bi x}}
\def\vPhi{{\vec{\Phi}}}
\def\vPsi{{\vec{\Psi}}}

\def\Dt {\tau}

\title[A linear MBP-preserving  BDF2 scheme for the Allen-Cahn equation]
{A linear second-order maximum bound principle-preserving BDF scheme for the Allen-Cahn equation with a general mobility$^*$}
\author[Dianming Hou, Lili Ju and Zhonghua Qiao]
{Dianming Hou$^{1}$
\quad
Lili Ju$^{2}$
\quad
Zhonghua Qiao$^{3}$
}
\thanks{\hskip -12pt
$^{1}$School of Mathematics and Statistics, Jiangsu Normal University, Xuzhou, Jiangsu 221116, China. Email: {\tt dmhou@stu.xmu.edu.cn}. Current address: Department of Applied Mathematics, The Hong Kong Polytechnic University, Hung Hom, Kowloon, Hong Kong. D. Hou's work is partially supported by Natural Science Foundation of China grant 12001248,  Jiangsu Province Higher Education Institutions grant
BK20201020,  Jiangsu Province Universities  Science Foundation grant 20KJB110013 and  Hong Kong Polytechnic University grant 1-W00D. \\
$^{2}$Department of Mathematics, University of South Carolina, Columbia, SC 29208, USA. Email:  {\tt ju@math.sc.edu}. L. Ju's work is partially supported by US National Science Foundation grant DMS-2109633.\\
$^{3}$Department of Applied Mathematics, The Hong Kong Polytechnic University, Hung Hom, Kowloon, Hong Kong. Email:  {\tt zqiao@polyu.edu.hk}. Z. Qiao's work   is partially supported by the Hong Kong Research Grants Council RFS grant RFS2021-5S03 and GRF grant 15302122, the Hong Kong Polytechnic University grant 4-ZZLS, and CAS AMSS-PolyU Joint Laboratory of Applied Mathematics.
}
\begin{document}
\graphicspath{{figures/},}
\date {\today}
\maketitle

\begin{abstract}
In this paper, we propose and analyze a linear second-order numerical method for solving the  Allen-Cahn equation with a general mobility. The proposed fully-discrete scheme is carefully constructed based on the combination of first and second-order backward differentiation formulas with nonuniform time steps for temporal approximation and  the central finite difference for spatial discretization.  The discrete maximum bound principle is proved of the proposed scheme by using the kernel recombination technique under certain mild constraints on the time steps and the ratios of adjacent time step sizes. Furthermore, we  rigorously derive the discrete $H^{1}$ error estimate  and energy stability for the classic constant mobility  case and the $L^{\infty}$ error estimate for the general mobility case. Various numerical experiments are also presented to validate the theoretical results and demonstrate the performance of the proposed method with a time adaptive strategy.
\end{abstract}

\section{Introduction}
\setcounter{equation}{0}

In this paper, we consider the following Allen-Cahn equation with a general mobility:
\brr\label{prob}
\begin{cases}
\dps\frac{\partial \phi}{\partial t}=\dps-M(\phi)\mu,&\quad (\x,t)\in\Omega\times (0,T],\\
\mu=\dps-\varepsilon^{2}\Delta\phi+F'(\phi),&\quad (\x,t)\in\Omega\times (0,T],
\end{cases}
\err
with the initial condition
$\phi(\x,0)=\phi_{0}(\x)$ for any $\x\in\Omega$
and subject to the homogeneous Neumann  or the periodic boundary condition,  where $\Omega$ is a  bounded Lipschitz domain in $\mathbb{R}^{d}$ $(d=1,2,3)$,
 $T>0$ is the terminal time,
 $\phi(\x,t)$ is the unknown function, $\varepsilon>0$ represents the interfacial width parameter, $M(\phi)\geq 0$ is  a general mobility function,
 and $F(\phi)=\frac14(1-\phi^{2})^{2}$ is  the double-well potential function.
This problem has a structure of $L^2$ gradient flow corresponding to the following free energy functional $E(\phi)$, defined by
\bq\label{energy}
E(\phi)=\int_{\Omega}\Big(\frac{\varepsilon^{2}}{2}|\nabla\phi|^{2}+F(\phi)\Big)d\x.
\eq
This structure implies that solution of  \eqref{prob} will approach to a steady state as $t\rightarrow\infty$, provided all steady states are isolated.
It is a physically attractive and thermodynamically-consistent model often used to describe the transitions of the phases in the binary alloys.
More specifically, the Allen-Cahn equation \eqref{prob}  satisfies the following energy dissipation law
 \bq\label{EDlaw}
 \frac{d}{dt}E(\phi)
 =-\int_{\Omega}M(\phi)\mu^{2}d\x\leq0,
 \eq
which indicates that the free energy $E(\phi)$ monotonically decreases in time.
Furthermore, the Allen-Cahn equation \eqref{prob} satisfies the maximum bound principle (MBP), i.e., $|\phi(\x,t)|\leq1$ if $|\phi(\x,0)|\leq1$ for any $\x\in\Omega$ and $t\geq 0$, and we refer to \cite{STY16} for more discussions.
The MBP and energy dissipation law are two important features of the equation \eqref{prob},  and thus it is highly desired for the numerical schemes to preserve these physical properties in the discrete level.

During the past decades, there have been extensive works devoted to the development of numerical methods for the Allen-Cahn equation \eqref{prob} with preservation of discrete MBP and energy stability, especially for the  constant mobility case.
First-order (in time) linear stabilized schemes with central finite difference method for spatial discretization were obtained  for the Allen-Cahn equation \eqref{prob} with a constant mobility in \cite{TY16} and the generalized Allen-Cahn equation with an advection term in \cite{STY16}, which are unconditionally energy stable and preserve the MBP simultaneously.  A second-order convex splitting scheme based on Crank-Nicolson approach was investigated for fractional-in-space Allen-Cahn equation  in \cite{HTY17}, in which the discrete MBP and energy dissipation were rigorously established. However, it results in a nonlinear system to be solved at each time step. Hou et al. \cite{HL20} developed a stabilized second-order Crank-Nicolson/Adams-Bashforth scheme for the Allen-Cahn equation, which preserves the discrete MBP and energy stability conditionally, and leads to solutions of only linear Poisson-type equations with constant coefficients at each time step. Recently, Cheng et al. \cite{CS22_1,CS22_2} proposed a Lagrange multiplier approach to construct positivity and bound preserving schemes for a class of semi--linear and quasi--linear parabolic equations. They have provided a new interpretation for the cut-off approach.
Based on cut-off approach and the scalar auxiliary variable (SAV) method \cite{Shen17_1,ALL19}, Yang et al. developed a class of arbitrarily high-order energy-stable and maximum bound preserving schemes for Allen-Cahn equation with a constant mobility in \cite{YYZ22}.

Du et al. developed first-order exponential time differencing (ETD) and second-order ETD Runge-Kutta (ETDRK2) schemes for the nonlocal Allen-Cahn equation, which preserves the discrete MBP unconditionally in \cite{DJLQ19}, and later they also established an abstract framework on the MBP  for a class of semilinear parabolic equations in \cite{DJLQ21}.
These ETD approaches were also successfully applied to the conservative Allen-Cahn equations in \cite{JJLL21,LJCF21} of preserving the MBP and mass conservation in the discrete level, and the molecular beam epitaxial model \cite{CQW19,CLLWW20} of maintaining the discrete energy stability.
Combining SAV technique with the stabilized first-order ETD and ETDRK2 methods, Ju et al. \cite{JLQ22_1,JLQ22_2} successfully constructed both the energy dissipation law and the MBP preserving schemes for a class of
Allen-Cahn type gradient flows.
The unconditional energy stability of the stabilized ETDRK2 scheme for the gradient flows are also established in \cite{FY22}.
Based on integrating factor Runge-Kutta (IFRK) method, high-order MBP preserving schemes  in time were recently developed for the semilinear parabolic equations in \cite{JLQY21}.
Subsequently,  a family of  stabilized IFRK schemes (up to the third-order and fourth-order) were proposed in \cite{LLJF21,ZYQGS21,ZYQS21} to preserve the discrete MBP unconditionally.
Recently, an arbitrarily high-order multistep exponential integrator method
was presented in \cite{LYZ20} by enforcing the maximum bound via a cut-off operation.
However, these high-order MBP-preserving ETD and IFRK methods seem difficultly to be extended to the problems with variable mobilities, since they are derived from either the variation-of-constant formula or an exponential transformation of the solution.
We also would like to remark  that all above MBP-preserving and energy stable scheme are based on the single time-stepping approach.
There also exist few research and results on the MBP preservation of multiple time-stepping method, such as the popular high-order BDF schemes.
Liao et al. studied the two-step second-order backward differentiation formula (BDF2) scheme  for the time discretization of the Allen-Cahn equation with a constant mobility in \cite{LTZ20}, in which
the MBP preservation and energy stability are established under certain mild constraints on the time steps and the ratios of adjacent time step sizes. However, it uses fully implicit treatment for the nonlinear term and thus leads to solving a nonlinear system at each time step. There also have been a lot of research work \cite{LQW21,HHW20,MQWZ20,LCWYH18,CCWWW12,CWWW22,CFWW19,YCWW18} on high-order BDF schemes for gradient flows, which maintain certain discrete energy stability.

 Another common feature of the Allen-Cahn equation \eqref{prob} is that its evolution process often takes quite long time before it settles at a steady state.
 Moreover, it usually undergoes both fast and slow changing stages during the whole evolution process.
Therefore, it is also highly useful to develop high-order structure-preserving numerical schemes with variable time steps for the Allen-Cahn equation, so that
some existing time adaptive strategies can be easily applied.
In this paper, we will propose and analyze an efficient  linear second-order numerical method with nonuniform time steps for solving the Allen-Cahn equation with a general (constant or variable) mobility,
which is based on the nonuniform BDF2 approach \cite{BJ98,CWYZ19,LTZ20,HQ21}  and preserves the discrete MBP under some mild constraints
like \cite{LTZ20}.

The rest of the paper is organized as follows: In Section 2, we first review some preliminaries on the temporal and spatial discretization, and then
propose the linear second-order BDF scheme for Allen-Cahn equation \eqref{prob}. Next we  establish  the discrete  MBP  of the proposed scheme using the kernel recombination technique in Section 3. In Section 4,  some results on error estimates in the $L^{\infty}$ and $H^{1}$ norms and energy stability are rigorously derived.
Several  examples are tested in Section 5 to numerically validate the theoretical prediction and demonstrate the performance of the proposed scheme.
Finally,  some concluding remarks are drawn in section 6.

\section{The linear second-order BDF scheme with nonuniform time steps}
\setcounter{equation}{0}

We first briefly review the BDF2 formula for approximating time derivative and the central finite difference for discretizing the Laplacian, and then propose a linear second-order BDF scheme for the Allen-Cahn equation with a general mobility \eqref{prob}. Without loss of generality, we focus on the two-dimensional problem ($d=2$) with the homogenous Neumann boundary condition, i.e., $\frac{\partial\phi}{\partial \n}\big|_{\partial\Omega}=0 $ in what follows. It is easy to extend the corresponding results to the cases of higher dimensional spaces and/or  the periodic boundary condition.

\subsection{The BDF2 formula with nonuniform time steps  and its reformulation through kernel recombination}\label{lerrec}

 Let $\{\Dt_{n}=t_{n}-t_{n-1}>0\}_{n=1}^N$ denote the time step sizes of a general partition of the time interval $[0,T]$ such that $t_0=0$ and $\sum_{n=1}^{N} \Dt_n=T$, and $\{\gamma_{n+1}=\frac{\Dt_{n+1}}{\Dt_{n}}>0\}_{n=1}^{N-1}$ denote the ratios of the corresponding two adjacent time step sizes. Define $\tau= \max\limits_{1\leq n\leq N}\Dt_{n}$ as the maximum time step size of such time partition and $\gamma_{max}=\max\limits_{1\leq n\leq N}\gamma_{n}$ as the maximum adjacent time-step ratio.

 For any function $\phi(t)$ defined on $[0,T]$,   denote $\Pi_{2,n}\phi(t)$ as its quadratic interpolation operator using the three points
$(t_{n-1},\phi(t_{n-1})),$ $(t_{n},\phi(t_{n}))$ and $(t_{n+1},\phi(t_{n+1}))$, and we then have
\beq
\Pi_{2,n}\phi(t)=\phi(t_{n-1})\frac{(t-t_{n})(t-t_{n+1})}{\Dt_{n}(\Dt_{n}+\Dt_{n+1})}-\phi(t_{n})\frac{(t-t_{n-1})(t-t_{n+1})}{\Dt_{n}\Dt_{n+1}}+\phi(t_{n+1})\frac{(t-t_{n-1})(t-t_{n})}{(\Dt_{n}+\Dt_{n+1})\Dt_{n+1}}
\eeq
for any $ t\in [t_{n-1},t_{n+1}]$  and consequently
\beq
\dps\frac{\partial \Pi_{2,n}\phi}{\partial t}(t_{n+1})=\frac{1}{\Dt_{n+1}}\Big(\frac{1+2\gamma_{n+1}}{1+\gamma_{n+1}}\phi(t_{n+1})-(1+\gamma_{n+1})\phi(t_{n})+\frac{\gamma^{2}_{n+1}}{1+\gamma_{n+1}}\phi(t_{n-1})\Big).
\eeq
Thus  the correspondingly derived second-order  BDF  approximation  to $\phi'(t)$ at $t=t_{n+1}$ reads:
\brr\label{BDF_s}
\phi'(t_{n+1})\approx \dps F^{n+1}_{2}\phi=\;&\dps\frac{1}{\Dt_{n+1}}\Big(\frac{1+2\gamma_{n+1}}{1+\gamma_{n+1}}\phi^{n+1}-(1+\gamma_{n+1})\phi^{n}+\frac{\gamma^{2}_{n+1}}{1+\gamma_{n+1}}\phi^{n-1}\Big)\\[10pt]
=\;&\dps b^{n}_{0}\delta_{\Dt}\phi^{n+1}+b^{n}_{1}\delta_{\Dt}\phi^{n},\qquad  \,n=1,2,\cdots,N-1,
\err
where $\delta_{\Dt}\phi^{n+1}=\phi^{n+1}-\phi^{n}$, $\phi^{n}$ is a certain approximation to $\phi(t_{n})$, and the discrete convolution kernels
$$b^{n}_{0}=\frac{1+2\gamma_{n+1}}{\Dt_{n+1}(1+\gamma_{n+1})}>0,\quad b^{n}_{1}=-\frac{\gamma^{2}_{n+1}}{\Dt_{n+1}(1+\gamma_{n+1})}<0.$$
For $n=0$, if we set $b^{0}_{0}=1/\Dt_{1}$ and $b^{0}_{1}=0$, then $F^{1}_{2}\phi = \frac{\delta_{\Dt}\phi^{1}}{\Dt_{1}}$  degrades to the  first-order BDF approximation
to $\phi'(t)$ at $t_1$, i.e., the well-know backward Euler approximation
\beq
\phi'(t_{n+1})\approx F^{n+1}_{1}\phi =\frac{ \delta_{\Dt}\phi^{n+1}}{\Dt_{n+1}},\quad  n=0,1,\cdots,N-1.
\eeq

A novel technique through variable-weights recombination of a new specially-created variable was first proposed in \cite{LX16} to achieve $3-\alpha$ order accuracy for the discrete form of $\alpha$-th order fractional Caputo derivative under the uniform time partition,  in which the reformed convolution kernels are positive and monotone and play an important role in stability and convergence analysis. Also see \cite{LMZ19,LTZ20} for some recent developments in this direction.
Following this  kernel recombination technique, we define a new variable $\psi$ as
 \bq\label{reform1}
  \psi^{n+1}=\phi^{n+1}-\eta\phi^{n},\quad n=0,1,\cdots,N-1,
 \eq
  with $\psi^{0}=\phi^{0}$, where  $\eta$ is a constant parameter to be determined such that the reformed discrete convolution kernels are positive and monotone.
  Then we have for $n=0,1,\cdots,N-1,$
  \beq
  \phi^{n+1}=\dps\sum_{k=0}^{n+1}\eta^{n+1-k}\psi^{k}, \quad \delta_{\Dt}\phi^{n+1}=\dps\sum_{k=0}^{n}\eta^{n-k}\delta_{\Dt}\psi^{k+1}+\eta^{n+1}\phi^{0}.
  \eeq
  Combing \eqref{BDF_s} and the above identities, we can equivalently reform the BDF2 formula \eqref{BDF_s}  as follows
 \bq\label{reform-BDF2}
 \dps F^{n+1}_{2}\phi=\sum_{k=0}^{n}d^{n}_{n-k}\delta_{\Dt}\psi^{k+1}+d^{n}_{n+1}\psi^{0},\quad 1\leq n\leq N-1,
 \eq
 where the reformed discrete convolution kernels are defined by
 \bq\label{eqn1_1}
 d^{n}_{0}=b^{n}_{0}, \quad d^{n}_{k}=\eta^{k-1}\big(b^{n}_{0}\eta+b^{n}_{1}\big),\quad 1\leq k\leq n+1.
 \eq
 Thus we have
 \bq\label{equn1}
 d^{n}_{k+1}=\eta d^{n}_{k},\quad 1\leq k\leq n.
 \eq
 In order to make $\{d_k^n\}_{k=0}^{n+1}$ positive and decreasing, i.e., $ d^{n}_{0}\geq d^{n}_{1}\geq\cdots\geq d^{n}_{n+1}\geq0,$ we need to require  $\eta$ to satisfy that
 \beq
 0<\frac{\gamma^{2}_{n+1}}{1+2\gamma_{n+1}}=-\frac{b^{n}_{1}}{b^{n}_{0}}\leq\eta<1
 \eeq
 for all $n=1, 2,\cdots,N-1$.
 Since $0<\gamma_{n+1}\leq\gamma_{max}$ and $\frac{x^{2}}{1+2x}$ is increasing in $(0,+
 \infty)$,  we then have
 \bq\label{eta}
 \frac{\gamma^{2}_{max}}{1+2\gamma_{max}}\leq\eta<1,
 \eq
  which also implies $0<\gamma_{max}<1+\sqrt{2}$.

\subsection{The central finite difference for the Laplacian}

We firstly recall some notations and results of the discrete function spaces and operators from \cite{WISE10,BLWW13,BHLWWZ13,SWWW12,HWWL09,WWL09,LSR19,WW88}.
Let $\Omega=(0,L_{x})\times(0,L_{y}),$ and we also assume $L_{x}=L_{y}=L$ and the spatial grid spacing $h=L/M$ for simplicity.
We first define the  following two finite grid sets:
$$
\mathbf{E}=\{x_{i+\frac{1}{2}}=ih\;\big|\; i=0,1, \cdots,M\},\qquad \mathbf{C}=\{x_{i}=\big(i-\textstyle\frac{1}{2}\big)h\;\big|\; i=1, \cdots,M\},
$$
and then we introduce the following discrete function spaces:
\bry
\mathcal{C}_{h}=&\dps\{U: \mathbf{C}\times\mathbf{C}\rightarrow\mathbb{R}\;\big|\;U_{i,j}, \;1\leq i,j\leq M\},\\[4pt]
 e^{x}_{h}=&\dps\{U: \mathbf{E}\times\mathbf{C}\rightarrow\mathbb{R}\;\big|\;U_{i+\frac{1}{2},j}, \;0\leq i\leq M,\;1\leq j\leq M\},\\[4pt]
 e^{y}_{h}=&\dps\{U: \mathbf{C}\times\mathbf{E}\rightarrow\mathbb{R}\;\big|\;U_{i,j+\frac{1}{2}}, \;1\leq i\leq M,\;0\leq j\leq M\},\\[4pt]
 e^{x}_{0,h}=&\dps\{U\in e^{x}_{h}\;\big|\;U_{\frac{1}{2},j}=U_{M+\frac{1}{2},j}=0,\; 1\leq j\leq M\},\\[4pt]
 e^{y}_{0,h}=&\dps\{U\in  e^{y}_{h}\;\big|\;U_{i,\frac{1}{2}}=U_{i,M+\frac{1}{2}}=0,\; 1\leq i\leq M\}.
\ery
Under the homogeneous Neumann boundary condition, the discrete gradient operator $\nabla_{h} =(\nabla^x_{h}, \nabla^y_{h}): \mathcal{C}_h\rightarrow( e_{0,h}^{x},  e^{y}_{0,h})$
is defined by
\bq
(\nabla^x_{h}U)_{i+\frac{1}{2},j}=\dps\frac{U_{i+1,j}-U_{i, j}}{h}, \quad 1\leq i\leq M-1,\;1\leq j\leq M,
\eq
\bq
(\nabla^y_{h}U)_{i,j+\frac{1}{2}}=\dps\frac{U_{i,j+1}-U_{i, j}}{h}, \quad 1\leq i\leq M,\;1\leq j\leq M-1
\eq
for any $U\in \mathcal{C}_h$, and the discrete divergence operator $\nabla_{h}\cdot:( e_{h}^{x},  e^{y}_{h})\rightarrow\mathcal{C}_h$ is represented by
\bq
(\nabla_{h}\cdot (U^{x},U^{y})^T)_{i,j}=\tps\frac{U^x_{i+1/2,j}-U^x_{i-1/2,j}}{h}+\tps\frac{U^y_{i,j+1/2}-U^y_{i,j-1/2}}{h},\quad 1\leq i,j\leq M
\eq
for any $(U^{x},U^{y})^{T}\in ( e_{h}^{x},  e^{y}_{h}).$ Then the  discrete LapLacian $\Delta_{h}:\mathcal{C}_{h}\rightarrow\mathcal{C}_{h}$ by the central finite difference is  defined by
\beq
(\Delta_{h}U)_{i,j}=(\nabla_{h}\cdot(\nabla_{h} U))_{i,j},\quad 1\leq i,j\leq M.
\eeq
The two average operators $a_{x}: e^{x}_{h}\rightarrow\mathcal{C}_h$ and $a_{y}: e^{y}_{h}\rightarrow\mathcal{C}_h$ are defined by
$$(a_{x}U)_{i,j}=\tps\frac12({U_{i+\frac{1}{2},j}+U_{i-\frac{1}{2},j}}),\quad 1\leq i,j\leq M,$$
 and
$$(a_{y}U)_{i,j}=\tps\frac12({U_{i,j+\frac{1}{2}}+U_{i,j-\frac{1}{2}}}),\quad 1\leq i,j\leq M$$
for any $U\in  e^{y}_{h}$.
We define some related discrete inner-products as follows:
\bry
&\big<U,V\big>_{\Omega}=\dps h^{2}\sum_{i,j=1}^{M}U_{i,j}V_{i,j},\quad \forall\,U,V\in\mathcal{C}_{h},\\[7pt]
&[U^x,V^x]_{x}=\big<a_{x}(U^xV^x),1\big>_{\Omega},\quad \forall\,U^x,V^x\in e^{x}_{h},\\[7pt]
&[U^y,V^y]_{y}=\dps\big<a_{y}(U^yV^y),1\big>_{\Omega},\quad \forall\,U^y,V^y\in e^{y}_{h},\\[7pt]
& [(U^{x},U^{y})^{T},(V^{x},V^{y})^{T}]_{\Omega}=[U^{x},V^{x}]_{x}+[U^{y},V^{y}]_{y}.\\
\ery
Then we  have the following result for the discrete analogue of integration by parts.
\begin{lemma}[\cite{LSR19,WW88}]\label{intpat}
For any $U,V\in\mathcal{C}_{h}$, it holds
\beq
-\big<\Delta_{h}U,V\big>_{\Omega}=[\nabla_{h}U,\nabla_{h}V]_{\Omega}.
\eeq
\end{lemma}
For any $U\in\mathcal{C}_{h}$, we define the following discrete $L^2$, $H^1$ and $L^{\infty}$ norms/semi-norms:
\bry
&\|U\|^{2}_{h}=\big<U,U\big>_{\Omega},\\[7pt]
& \|\nabla_{h}U\|^{2}_{h}=[\nabla_{h}U,\nabla_{h}U]_{\Omega}=[d_{x}U,d_{x}U]_{x}+[d_{y}U,d_{y}U]_{y},\\[7pt]
&\dps\|U\|^{2}_{H^{1}_{h}}=\|U\|^{2}_{h}+\|\nabla_{h}U\|^{2}_{h},\quad \dps\|U\|_{\infty}=\max_{1\leq i\leq M} \sum_{j=1}^{M}|U_{i,j}|
\ery
For convenience of description, we also define $\vec{U}\in{\mathbb R}^{M^2}$ as the vector representation of $U\in\mathcal{C}_{h}$,  in which the elements are arranged first along the $x$-direction then along the $y$-direction. Note that we do not differ them  in places  there is no ambiguity.

\subsection{The linear second-order BDF scheme for the Allen-Cahn equation}

Denote by $\Pi_{\mathcal{C}_h}$ the operator pointwisely limiting a function onto $\mathcal{C}_h$.
Let us first  recall the fully-discrete linear   first-order  BDF scheme (called ``BDF1'') proposed in \cite{STY16,TY16} for solving the Allen-Cahn equation with a general mobility \eqref{prob}: given $\Phi^0=\Pi_{\mathcal{C}_h}\phi_0$,  for $n=0,1,\cdots,N-1$, find $\Phi^{n+1}\in\mathcal{C}_{h}$ such that
 \bq\label{BDF_1}
 F^{n+1}_{1}\Phi-\varepsilon^{2}M(\Phi^{n})\Delta_{h}\Phi^{n+1}+f(\Phi^{n})+S(\Phi^{n+1}-\Phi^{n})=0,
  \eq
  where $F^{n+1}_{1}\Phi= \frac{\Phi^{n+1}-\Phi^{n}}{\Dt_{n+1}}$ and $f(\phi) =M(\phi)F'(\phi)$ and $S\ge 0$ is a constant stabilizing parameter. We will denote the scheme \eqref{BDF_1} as $\Phi^{n+1} = {\rm BDF1}(\Phi^n,\tau_{n+1})$.
  The above linear BDF1 scheme \eqref{BDF_1} also be rewritten in the following vector form:
  \bq\label{BDF1_tsor}
 F^{n+1}_{1}\vPhi-\varepsilon^{2}\Lambda^{n}D_{h}\vPhi^{n+1}+f(\vPhi^{n})+S(\vPhi^{n+1}-\vPhi^{n})=0,
  \eq
   where
$D_{h}=I\otimes G+G\otimes I\in {\mathbb R}^{M^2\times M^2}$ with $I$ denoting the identity matrix (with the matched dimensions) and
\[ G=\frac{1}{h^{2}}
\begin{pmatrix}
-1&1&&&&\\
1&-2&1&&&\\
&\ddots &\ddots &\ddots&\\
&&1&-2&1&\\
&&&1&-1&\\
\end{pmatrix}_{M\times M},\]
 and   $f(\vPhi^{n})=\Lambda^{n}\big(\big(\vPhi^{n}\big)^{.3}+\vPhi^{n}\big)$  is defined elementwise with  the diagonal matrix $\Lambda^{n}=\mbox{diag}(M(\vPhi^{n}))$.
Clearly, $D_h$ is  the corresponding matrix representation of $\Delta_h$.

In analogous to the energy $E(\phi)$ defined in \eqref{energy}, we define the discrete energy $E_{h}(\vPhi^{n})$ as
\bq\label{dis_eg}
E_{h}(\vPhi^{n})=-h^{2}\frac{\varepsilon^{2}}{2}(\vPhi^{n})^{T}D_{h}\vPhi^{n}+h^{2}\sum_{i=1}^{M^{2}}F(\vPhi^{n}_{i})=\frac{\varepsilon^{2}}{2}[\nabla_{h}\Phi^n,\nabla_{h}\Phi^n]_{\Omega}+\big<F(\Phi^{n}),1\big>_{\Omega}.
\eq
Then the unconditional energy stability and the discrete maximum bound principle of the fully-discrete BDF1 scheme \eqref{BDF_1} hold as stated in the following lemma, and we refer to Theorem 3.2 in \cite{STY16} and Theorem 3 in \cite{TY16} for details.
\begin{lemma}[\cite{STY16,TY16}]\label{lem1}
Asssume that $\|\vPhi^{0}\|_{\infty}\leq 1$ and the stabilizing parameter
\bq\label{eqn1_3}
S\geq \max_{\rho\in[-1,1]}\big( M'(\rho)F'(\rho)+M(\rho)F''(\rho)\big),
\eq
then  it unconditionally holds for the BDF1 scheme \eqref{BDF_1} that $\|\vPhi^{n+1}\|_{\infty}\leq1$ for $n=0,1,\cdots,N-1$.
 Particularly, if the mobility function $M(\phi)\equiv1$, then
\bq
E_{h}(\vPhi^{n+1})\leq E_{h}(\vPhi^{n})
\eq
for all $n=0,1,\cdots,N-1,$ provided that $S\geq2$.
\end{lemma}

Now we are ready to construct a fully-discrete linear  second-order BDF scheme with nonuniform time steps (called ``BDF2'' hereafter)
for the Allen-Cahn equation with a general mobility \eqref{prob} under the homogenous Neumann boundary condition: given $\Phi^0=\Pi_{\mathcal{C}_h}\phi_0$, compute $\Phi^{1} = {\rm BDF1}(\phi^0,\tau_{1})$ and for $n=1,2\cdots,N-1$, find $\Phi^{n+1}\in\mathcal{C}_{h}$ such that
	\begin{subequations}\label{BDF_2}
	      \begin{empheq}[left=\empheqlbrace]{align}
		&\Phi^{*,n+1} = {\rm BDF1}(\Phi^n,\tau_{n+1}),\label{eqn41}\\
		&F^{n+1}_{2}\Phi-\varepsilon^{2}M(\Phi^{*,n+1})\Delta_{h}\Phi^{n+1}+f(\Phi^{*,n+1})+S(\Phi^{n+1}-\Phi^{*,n+1})=0,\label{eqn4}
		\end{empheq}
       \end{subequations}
      where $F^{n+1}_{2}\Phi = b^{n}_{0}(\Phi^{n+1}-\Phi^{n})+b^{n}_{1}(\Phi^{n}-\Phi^{n-1})$.
 We denote the scheme \eqref{BDF_2} as $\phi^{n+1} = {\rm BDF2}(\phi^n,\phi^{n-1},\tau_{n+1},\tau_n)$.
 The above linear BDF2 scheme \eqref{BDF_2} can be rewritten in the following vector form:
 	\begin{subequations}\label{BDF2_tsor}
	      \begin{empheq}[left=\empheqlbrace]{align}
		&\vPhi^{*,n+1} = {\rm BDF1}(\vPhi^n,\tau_{n+1}),\label{BDF2_tsor1}\\
		&F^{n+1}_{2}\vPhi-\varepsilon^{2}\Lambda^{*,n+1}D_{h}\vPhi^{n+1}+f(\vPhi^{*,n+1})+S(\vPhi^{n+1}-\vPhi^{*,n+1})=0,\label{BDF2_tsor2}
		\end{empheq}
       \end{subequations}
      where $\Lambda^{*,n+1}=\mbox{diag}(M(\vPhi^{*,n+1}))$.

  \section{The discrete maximum bound principle}
  \setcounter{equation}{0}

  In this section, we will prove the discrete maximum bound principle of the proposed BDF2 scheme \eqref{BDF_2} using the kernel recombination technique described in Section \ref{lerrec}. Define $\vPsi^{n}=\vPhi^{n}-\eta\vPhi^{n-1}$, and then we can  combine \eqref{reform1} and \eqref{reform-BDF2} to obtain
  the following kernel recombination form for \eqref{BDF2_tsor2}: for $n=1,2,\cdots,N-1$,
 \bq\label{eqn1}
 \big((d_{0}^{n}+S)I-\varepsilon^{2}\Lambda^{*,n+1}D_{h}\big)\vPhi^{n+1}=\eta d^{n}_{0}\vPhi^{n}+\sum_{k=0}^{n}(d^{n}_{n-k}-d^{n}_{n-k+1})\vPsi^{k}+S\vPhi^{*,n+1}-f(\vPhi^{*,n+1}).
 \eq
 Substituting $\vPhi^{n+1}=\dps\sum_{k=0}^{n+1}\eta^{n+1-k}\vPsi^{k}$ into \eqref{eqn1} yields

  \brr\label{eqn_1}
 \big((d_{0}^{n}+S)I-\varepsilon^{2}\Lambda^{*,n+1}D_{h}\big)\vPsi^{n+1}
 =\dps\sum_{k=0}^{n}Q^{n}_{n-k}\vPsi^{k}+S\vPhi^{*,n+1}-f(\vPhi^{*,n+1}),
 \err
 where
 \bq\label{eqn1_2}
 Q^{n}_{k}=(d^{n}_{k}-d^{n}_{k+1}-S\eta^{k+1})I+\eta^{k+1}\varepsilon^{2}\Lambda^{*,n+1}D_{h}, \quad 0\leq k\leq n.
 \eq
 The following result for the estimation of $Q^{n}_{k}$ holds (we also refer to Lemma 4.1 in \cite{LTZ20} which is only for the specific case $M(\phi)\equiv 1$).
 \begin{lemma}\label{lemm1}
 Let $n$ be any fixed integer such that $1\leq n\leq N-1$ and suppose $\|\vPhi^{*,n+1}\|_{\infty}\leq 1$.
 Assume that $0<\gamma_{n+1}<1+\sqrt{2},$ the parameter $\eta$ satisfies \eqref{eta}, and
 \bq\label{res_tau}
 \Dt_{n+1}\leq \frac{g(\gamma_{n+1},\eta)}{S+4L\varepsilon^{2}/h^{2}},
 \eq
 where $L=\dps\max_{\rho\in [-1,1]}M(\rho)$ and
  \beq
   g(s,z)=\frac{(1-z)\big((1+2s)z-s^{2}\big)}{z^{2}(1+s)},\quad s\in(0,\gamma_{max}],~z\in\Big[\frac{\gamma_{max}^{2}}{1
   +2\gamma_{max}},1\Big).
 \eeq
Then it holds
\bq\label{est_q}
 \|Q^{n}_{k}\|_{\infty}\leq d^{n}_{k}-d^{n}_{k+1}-S\eta^{k+1}, \quad  \forall\,0\leq k\leq n.
 \eq
 \end{lemma}
 \begin{proof}
 From the definition of $Q^{n}_{k}$ in \eqref{eqn1_2} and \eqref{eqn1_1}, it follows
 \bry
 Q^{n}_{k}=\;&\big(\eta^{k-1}(b^{n}_{0}\eta+b^{n}_{1})-\eta^{k}(b^{n}_{0}\eta+b^{n}_{1})-S\eta^{k+1}\big)I+\eta^{k+1}\varepsilon^{2}\Lambda^{*,n+1}D_{h},\\[5pt]
=\;& \eta^{k+1}\Big(\big(\eta^{-2}(1-\eta)(b^{n}_{0}\eta+b^{n}_{1})-S\big)I+\varepsilon^{2}\Lambda^{*,n+1}D_{h}\Big)\\[5pt]
=\;& \eta^{k+1}\Big(\Big(\dfrac{g(\gamma_{n+1},\eta)}{\tau_{n+1}}-S\Big)I+\varepsilon^{2}\Lambda^{*,n+1}D_{h}\Big),\quad 1\leq k\leq n,
 \ery
 which means that all the entries of $Q^{n}_{k}$ are nonnegative based on the definition of $D_h$, the fact of $\|\Lambda^{*,n+1}\|_{\infty}\leq L$,  and \eqref{res_tau}. Thus we deduce that
 \beq
 \|Q^{n}_{k}\|_{\infty}= \max_{1\leq i\leq M^{2}}\tps\sum_{j=1}^{ M^{2}}\big(Q^{n}_{k}\big)_{i,j}\leq d^{n}_{k}-d^{n}_{k+1}-S\eta^{k+1}, \quad  \forall\,1\leq k\leq n,
 \eeq
 by using the fact
$ \sum_{j=1}^{ M^{2}}\big(\Lambda^{*,n+1}D_{h}\big)_{i,j}= 0$ for any $1\leq i\leq M^{2}.$

 For the case of $k=0$, using \eqref{eqn1_1} and \eqref{res_tau},  we get
  \bry
 Q^{n}_{0}=\;&\dps \big(d^{n}_{0}-d^{n}_{1}-S\eta\big)I+\eta\varepsilon^{2}\Lambda^{*,n+1}D_{h}\\[5pt]
 =\;&\dps \big(b^{n}_{0}-b^{n}_{0}\eta-b^{n}_1-S\eta\big)I+\eta\varepsilon^{2}\Lambda^{*,n+1}D_{h}\\[5pt]
  =\;&\dps\eta\big(\big(\eta^{-2}(b_{0}^{n}\eta+b^{n}_1-b^{n}_0\eta^2-b_{1}^{n}\eta)-\eta^{-2}b^{n}_1-S\big)I+\varepsilon^{2}\Lambda^{*,n+1}D_{h}\big)\\[5pt]
 \geq\;&\eta\big(\big(\eta^{-2}(1-\eta)(b_{0}^{n}\eta+b_{1}^{n})-S\big)I+\varepsilon^{2}\Lambda^{*,n+1}D_{h}\big),
 \ery
 which means that all the entries of $Q^{n}_{0}$ are also nonnegative and consequently we obtain $\|Q^{n}_{0}\|_{\infty}\leq d^{n}_{0}-d^{n}_{1}-S\eta$ by
 similar arguments as above.
 \end{proof}

 This lemma plays an important role in deriving the MBP property of the BDF2 scheme \eqref{BDF_2}.
We also remark that the inequality \eqref{res_tau} doesn't explicitly give a principle for determining the range of feasible time step size $\Dt_{n+1}$ in practice, since $\gamma_{n+1}$ in the righthand side of  \eqref{res_tau} depends on $\Dt_{n+1}$.
 Next we drive a uniform upper bound for the time step size $\Dt_{n+1}$ independent on $\gamma_{n+1}$ such that the estimate \eqref{est_q} for the matrix $Q^{n}_{k}$ holds.
In the numerical simulations, one can always set a pre-determined maximum adjacent time-step ratio  $\gamma_{*}$ such that $\gamma_{n+1}\leq\gamma_{*}$ for all $n\geq1$. Since it is required that $0<\gamma_{n+1}<1+\sqrt{2}$ (see Section \ref{lerrec}), we  choose $\gamma_{*}$ from $[1,1+\sqrt{2})$.
Noting that
 \beq
 \frac{\partial g}{\partial s}(s,z)=\frac{(1-z)(-s^{2}-2s+z)}{z^{2}(1+s)^{2}},\quad s\in(0,\gamma_{*}],~z\in\Big[\frac{\gamma_{*}^{2}}{1+2\gamma_{*}},1\Big),
 \eeq
 and combining with $\sqrt{1+z}-1<\sqrt{2}-1<1\leq \gamma_{*}$, it can be verified that for any fixed $z$, $g(s,z)$ is increasing in $(0,\sqrt{1+z}-1)$ and decreasing in $(\sqrt{1+z}-1,\gamma_{*})$ with respect to $s$. Furthermore, since $g(0,z)=\frac{1-z}{z}>\frac{(1-z)(3z-1)}{2z^{2}}=g(1,z)\geq g(\gamma_{*},z)$ for $z\in\Big[\frac{\gamma_{*}^{2}}{1+2\gamma_{*}},1\Big)$, we have
 \beq
 g(\gamma_{*},z)\leq g(\gamma_{n+1},z)
 \eeq
 for all $\gamma_{n+1}\in(0,\gamma_{*})$ and $z\in\Big[\frac{\gamma_{*}^{2}}{1+2\gamma_{*}},1\Big).$
 Thus, it follows from Lemma \ref{lemm1} that the estimate \eqref{est_q} for the matrix $Q^{n}_{k}$ holds for $0<\gamma_{n+1}\leq\gamma_{*}<1+\sqrt{2}$,  and
 \bq\label{est_tau2}
 \Dt_{n+1}\leq \frac{g(\gamma_{*},\eta)}{S+4L\varepsilon^{2}/h^{2}}, \quad \forall\,\eta\in\Big[\frac{\gamma_{*}^{2}}{1
   +2\gamma_{*}},1\Big).
 \eq
 Taking the fact
 \beq
 \frac{\partial g}{\partial \eta}(\gamma_{*},\eta)=\frac{2\gamma_{*}^{2}-(1+\gamma_{*})^{2}\eta}{(1+\gamma_{*})\eta^{3}},
 \eeq
 together with $\frac{\gamma_{*}^{2}}{1+2\gamma_{*}}<\frac{2\gamma_{*}^{2}}{(1+\gamma_{*})^{2}}<1$, we see that $g(\gamma_{*},\eta)$ is increasing in $\big(\frac{\gamma_{*}^{2}}{1+2\gamma_{*}}, \frac{2\gamma_{*}^{2}}{(1+\gamma_{*})^{2}}\big)$ and decreasing in $\big(\frac{2\gamma_{*}^{2}}{(1+\gamma_{*})^{2}},1\big)$ with respect to  $\eta$.
 Thus, the optimal value of $\eta$ for \eqref{est_tau2} is
 \begin{equation}\label{etasta}
  \eta_{*}=\frac{2\gamma_{*}^{2}}{(1+\gamma_{*})^{2}}.
 \end{equation}
 Summarizing the above discussions, we obtain the following result.
 \begin{lemma}\label{lem_1}
 Let $n$ be any fixed integer such that $1\leq n\leq N-1$ and suppose $\|\vPhi^{*,n+1}\|_{\infty}\leq 1$.
 Assume that  $0<\gamma_{n+1}\leq \gamma_{*}<1+\sqrt{2}, \eta=\eta_{*},$ and the time step size $\Dt_{n+1}$ satisfies
\bq\label{est_tau3}
 \Dt_{n+1}\leq \frac{\mathcal{G}(\gamma_{*})}{S+4L\varepsilon^{2}/h^{2}}
 \eq
 with
 \beq
 \mathcal{G}(\gamma_{*})=g\Big(\gamma_{*}, \frac{2\gamma_{*}^{2}}{(1+\gamma_{*})^{2}}\Big)=\frac{(1+2\gamma_{*}-\gamma_{*}^{2})^{2}}{4\gamma^{2}_{*}(1+\gamma_{*})}.
 \eeq
 Then it holds
 \bq\label{est_q1}
 \|Q^{n}_{k}\|_{\infty}\leq d^{n}_{k}-d^{n}_{k+1}-S\eta_{*}^{k+1}, \quad \forall\,0\leq k\leq n. \eq
 \end{lemma}
 \begin{remark}
Note that $\mathcal{G}(\gamma_{*})$ is decreasing with respect to $\gamma_{*}\in(1,1+\sqrt{2}).$
 Especially, we have $\mathcal{G}(1)=\frac{1}{2}$ for the case of uniform time steps ($\gamma_{*}=1$), and $\mathcal{G}(2)=\frac{1}{48}$ for the case of  $\gamma_{*}=2.$
 \end{remark}
 In what follows, by default we always set $\eta=\eta_{*}$ which is defined in \eqref{etasta}.
We next state the following useful lemmas.
\begin{lemma}[\cite{TY16,LTZ20,HTY17}]\label{lemm2}
Suppose $B= (b_{i,j})$ is a real $P\times P$ matrix satisfying
\beq
b_{i,i}<0,\quad |b_{i,i}|\geq\max_{1\leq i \leq P}\tps\sum_{j\neq i}^{P}|b_{i,j}|.
\eeq
Let  $A=aI-B$ where $a>0$ is a constant, then
\beq
\|A\overrightarrow{U}\|_{\infty}\geq a \|\overrightarrow{U}\|_{\infty},\quad \forall\, \overrightarrow{U}\in\mathbb{R}^{P}.
\eeq
\end{lemma}
\begin{lemma}[\cite{TY16}]\
\label{para}
If the  stabilizing parameter $S$ satisfies \eqref{eqn1_3}, then
\bq\label{eqn1_4}
\big|S\rho-f(\rho)\big|\leq S,\quad\forall\,\rho\in[-1,1].
\eq
\end{lemma}
\begin{proof}
Let $h(\rho)=S\rho-f(\rho)$. From \eqref{eqn1_3}, we have
\bry
h'(\rho)=S-[M'(\rho)F'(\rho)+M(\rho)F''(\rho)]\geq0,\quad\forall\, \rho\in[-1,1].
\ery
Together with $h(-1)=-S$ and $h(1)=S$, we obtain \eqref{eqn1_4}.
\end{proof}
Now, we are ready to show the MBP of the BDF2 scheme \eqref{BDF_2}.
\begin{theorem}\label{thmMBP}
Assume that the stabilizing parameter $S$ satisfies \eqref{eqn1_3} and $ 0<\gamma_{n+1}\leq \gamma_{*}<1+\sqrt{2}$ for all $1\leq n \leq N-1$. In addition,  assume that
\bq\label{constr_tau1}
\Dt_{1}\leq\frac{1-\eta_{*}}{\eta_{*}(S+4L\varepsilon^{2}/h^{2})},
\eq
 and $\Dt_{n+1}$ satisfies \eqref{est_tau3} for  $n=1,2,\cdots,N-1.$
 If $\|\vPhi^{0}\|_{\infty}\leq 1$, then it holds for the BDF2 scheme \eqref{BDF_2} that
   $\|\vPhi^{n+1}\|_{\infty}\leq1$ for $n=0,1,\cdots,N-1$.
\end{theorem}
\begin{proof}
For the first step, i.e., $\vPhi^{1} = {\rm BDF1}(\vPhi^0,\tau_{1})$ when $n=0$, it follows directly from Lemma \ref{lem1} that $\|\vPhi^{1}\|_{\infty}\leq 1.$ Substituting $\vPhi^{1}=\vPsi^{1}+\eta_{*}\vPhi^{0}$ into \eqref{BDF_1} gives
\bq\label{eqn2}
\Big(\Big(\frac{1}{\Dt_{1}}+S\Big)I-\varepsilon^{2}\Lambda^{0}D_{h}\Big)\vPsi^{1}
=\Big(\Big(\frac{1-\eta_{*}}{\Dt_{1}}-\eta_{*}S\Big)I+\eta_{*}\varepsilon^{2}\Lambda^{0}D_{h}\Big)\vPsi^{0}+S\vPhi^{0}-f(\vPhi^{0}).
\eq
Noting the constraint \eqref{constr_tau1} together with the definition of $D_{h}$ and a similar analysis used in Lemma \ref{lemm1}, we derive that
\beq
\Big(\Big(\frac{1-\eta_{*}}{\Dt_{1}}-\eta_{*}S\Big)I+\eta_{*}\varepsilon^{2}\Lambda^{0}D_{h}\Big)_{i,j}\geq0,\quad 1\leq i,j\leq M^2,
\eeq
and consequently
\bq\label{eqn3}
\Big\|\Big(\frac{1-\eta_{*}}{\Dt_{1}}-\eta_{*}S\Big)I+\eta_{*}\varepsilon^{2}\Lambda^{0}D_{h}\Big\|_{\infty}\leq\frac{1-\eta_{*}}{\Dt_{1}}-\eta_{*}S.
\eq
From \eqref{eqn2}, \eqref{eqn3} and Lemma \ref{lemm2}, it follows that
\bry
\tps\Big(\frac{1}{\Dt_{1}}+S\Big)\|\vPsi^{1}\|_{\infty}\leq\;&\tps\Big\|\Big(\Big(\frac{1}{\Dt_{1}}+S\Big)I-\varepsilon^{2}\Lambda^{0}D_{h}\Big)\vPsi^{1}\Big\|_{\infty}\\[8pt]
\leq\;&\tps\Big\|\Big(\Big(\frac{1-\eta_{*}}{\Dt_{1}}-\eta_{*}S\Big)I+\eta_{*}\varepsilon^{2}\Lambda^{0}D_{h}\Big)\vPsi^{0}\Big\|_{\infty}+\big\|S\vPhi^{0}-f(\vPhi^{0})\big\|_{\infty}\\[8pt]
\leq\;&\tps\Big(\frac{1-\eta_{*}}{\Dt_{1}}-\eta_{*}S\Big)+S\\[5pt]
=\;&\dps\Big(\frac{1}{\Dt_{1}}+S\Big)(1-\eta_{*}),
\ery
where we have used Lemma \ref{para}. Thus we have $\|\vPsi^{1}\|_{\infty}\leq 1-\eta_{*}$.

Next, for any $1\leq n\leq N-1$, we assume $\|\vPhi^{k}\|_{\infty}\leq1$ and $\|\vPsi^{k}\|_{\infty}\leq 1-\eta_{*}$ for  $1\leq k\leq n.$
Using $\vPhi^{*,n+1} = {\rm BDF1}(\vPhi^n,\tau_{n+1})$, $\|\vPhi^{n}\|_{\infty}\leq 1,$ and Lemma \ref{lem1}, we obtain $\|\vPhi^{*,n+1}\|_{\infty}\leq1.$
Thus, together with \eqref{eqn1}, \eqref {equn1} and  Lemmas \ref{lemm2} and \ref{para}, we have
\bry
(d_{0}^{n}+S)\|\vPhi^{n+1}\|_{\infty} \leq\;&\dps\|\big((d_{0}^{n}+S)I-\varepsilon^{2}\Lambda^{*,n+1}D_{h}\big)\vPhi^{n+1}\|_{\infty}\\
\leq\;&\dps \eta_{*} d^{n}_{0}\|\vPhi^{n}\|_{\infty}+\sum_{k=0}^{n}(d^{n}_{n-k}-d^{n}_{n-k+1})\|\vPsi^{k}\|_{\infty}+\|S\vPhi^{*,n+1}-f(\vPhi^{*,n+1})\|_{\infty}\\
\leq\;&\dps  \eta_{*} d^{n}_{0}+\sum_{k=1}^{n}(d^{n}_{n-k}-d^{n}_{n-k+1})(1-\eta_{*})+(d^{n}_{n}-d^{n}_{n+1})+S\\
=\;&\dps  \eta_{*} d^{n}_{0}+(d^{n}_{0}-d^{n}_{n})(1-\eta_{*})+(1-\eta_{*})d^{n}_{n}+S\\[5pt]
=\;&d^{n}_{0}+S,
 \ery
which gives $\|\vPhi^{n+1}\|_{\infty}\leq1.$  Using \eqref{eqn_1} together with \eqref {equn1}, Lemmas \ref{lem_1}, \ref{lemm2} and \ref{para}, we get
 \bry
\dps(d_{0}^{n}+S)\|\vPsi^{n+1}\|_{\infty}\leq\;& \|\big((d_{0}^{n}+S)I-\varepsilon^{2}\Lambda^{*,n+1}D_{h}\big)\vPsi^{n+1}\|_
{\infty}\\
\leq\;&\dps\sum_{k=0}^{n}\|Q^{n}_{n-k}\|_{\infty}\|\vPsi^{k}\|_{\infty}+\|S\vPhi^{*,n+1}-f(\vPhi^{*,n+1})\|_{\infty}\\
\leq\;&\dps(1-\eta_{*})\sum_{k=1}^{n}(d_{n-k}^{n}-d_{n+1-k}^{n}-S\eta^{n+1-k}_{*})+(d_{n}^{n}-d_{n+1}^{n}-S\eta^{n+1}_{*})+S\\[9pt]
=\;&(d^{n}_{0}+S)(1-\eta_{*}).
 \ery
 which gives $\|\vPsi^{n+1}\|_{\infty}\leq1-\eta_{*}.$ The proof is completed.
\end{proof}

\section{Error analysis and energy stability}
\setcounter{equation}{0}

In this section, we investigate the error estimate and energy stability of the proposed BDF2 scheme \eqref{BDF_2}.
Let $\Phi(t)= \Pi_{\mathcal{C}_h}\phi(t)$ where $ \phi$
denotes the exact solution of \eqref{prob}.
We also use  $C$ and $C_i$'s to denote some needed generic positive constants independent of $h$ and $\Dt$.

\subsection{ Discrete $H^{1}$ error estimate and energy stability for the constant mobility case}

In this subsection, we study the discrete $H^{1}$ error estimate and energy stability of the BDF2 scheme \eqref{BDF_2}  for the Allen-Cahn equation with  constant mobility, i.e., $M(\phi)\equiv C>0$. Without loss of generality, we assume $M(\phi)\equiv1$ and thus 
\eqref{eqn1_3} becomes $S\geq 2$.
Firstly, we recall a useful  inequality (see \cite{HQ21,HQ22,HX22_1,LZ20,LJWZ21}) presented  below, which will play an important role in our error analysis and energy stability: for any $\{\gamma_n>0\}_{n=1}^{N+1}$,
\brr\label{ideq1}
\dps \big<F^{n+1}_{2}\vPhi, \vPhi^{n+1}-\vPhi^{n}\big>_{\Omega}\geq\;&\dps\Big(\frac{\gamma_{n+2}^{3/2}}{1+\gamma_{n+2}}\frac{\|\vPhi^{n+1}-\vPhi^{n}\|^{2}_{h}}{2\tau_{n+1}}
-\frac{\gamma_{n+1}^{3/2}}{1+\gamma_{n+1}}\frac{\|\vPhi^{n}-\vPhi^{n-1}\|^{2}_{h}}{2\tau_{n}}\Big)\\[10pt]
\;&\dps+\, G(\gamma_{n+1},\gamma_{n+2})\frac{\|\vPhi^{n+1}-\vPhi^{n}\|^{2}_{h}}{2\Dt_{n+1}}, \qquad n=1,2,\cdots, N-1,
\err
where $G(s,z)=\frac{2+4s-s^{3/2}}{1+s}-\frac{z^{3/2}}{1+z}.$  Note that $\gamma_{N+1}$ is not used in the BDF2 scheme.
It is easy to verify that for any fixed $z\in(0,+\infty)$, $G(s,z)$ is increasing in $(0,1)$ and decreasing in $(1,+\infty)$ with respect to $s$.
Then it follows from $G(0,z)=G(4,z)$ that for any $0<s,z\leq\gamma_{*}<1+\sqrt{2}$,
\beq
G(s,z)\geq \min\{G(0,\gamma_{*}),G(\gamma_{*},\gamma_{*})\}\geq G(0,\gamma_{*})>G(0,1+\sqrt{2})>0.
\eeq

Define the errors  $e^{n}=\vPhi^{n}-\vPhi(t_{n})$ and $e^{*,n}=\vPhi^{*,n}-\vPhi(t_{n})$.
With a reasonable requirement on the exact solution $\phi$ of the problem \eqref{prob}, we are able to establish a discrete $H^{1}$ error estimate for the BDF2 scheme \eqref{BDF_2}.

\begin{theorem}\label{th1}
 Assume that $0<\gamma_{n+1}\leq\gamma_{*}<1+\sqrt{2}$ for all $1\leq n \leq N-1$, $S\geq 2$, and the time step sizes satisfy \eqref{est_tau3} and \eqref{constr_tau1}. Let $\gamma_{N+1}$ be any number in $(0,\gamma_*)$.
 In addition, assume that $\Dt_{1}\leq C_1 \Dt^{\frac{4}{3}}$ and
$\phi\in W^{3,\infty}(0,T;L^{\infty}(\Omega))\cap L^{\infty}(0,T;W^{4,\infty}(\Omega))$.
Then it holds for the BDF2 scheme \eqref{BDF_2} in the constant mobility case that
\brr\label{equa4_1}
&\dps\frac{\gamma^{3/2}_{n+2}}{1+\gamma_{n+2}}\frac{\|e^{n+1}-e^{n}\|^{2}_{h}}{\Dt_{n+1}}
+\varepsilon^{2}\|\nabla_{h} e^{n+1}\|^{2}_{h}+S\|e^{n+1}\|_h^{2}\\[6pt]
&\qquad\qquad\leq \dps C\exp(T)\Big(\Dt^{4}\|\phi\|^{2}_{W^{3,\infty}(0,T;L^{\infty}(\Omega))}+h^{4}\|\phi\|^{2}_{L^{\infty}(0,T;W^{4,\infty}(\Omega))}\big)\Big)
\err
for all $0\leq n\leq N-1$.
\end{theorem}

\begin{proof}
It follows from $\|\vPhi\|_{\infty}\leq1, \|\vPhi^{n}\|_{\infty}\leq1$ (by the discrete MBP stated in Theorem \ref{thmMBP}), and $f(\cdot)\in C^{1}(\mathbb{R})$ that
\bq\label{equ1}
\max\{\|f(\vPhi)\|_{\infty},\|f^{'}(\vPhi)\|_{\infty},\|f(\vPhi^{n})\|_{\infty},\|f^{'}(\vPhi^{n})\|_{\infty}\}\leq C_2
\eq
 for all $n=0,1,\cdots, N.$
 Comparing \eqref{prob} and \eqref{BDF_2} gives the error equations of $e^{*,n+1}$ and $e^{n+1}$:
 \begin{subequations}\label{err}
  \begin{empheq}[left=\empheqlbrace]{align}
 \frac{ e^{*,n+1}-e^{n}}{\Dt_{n+1}}-\varepsilon^{2}\Delta_{h} e^{*,n+1}+S e^{*,n+1}=\;&S e^{n}-S(\vPhi(t_{n+1})-\vPhi(t_{n}))\nonumber\\
 &+f(\vPhi(t_{n+1}))-f(\vPhi^{n})+T_{1}^{n}+T_{2}^{n},\label{err_2}\\
 \dps F^{n+1}_{2}e-\varepsilon^{2}\Delta_{h} e^{n+1}+S e^{n+1}=\;&S e^{*,n+1}+f(\vPhi(t_{n+1}))-f(\vPhi^{*,n+1})+T_{2}^{n}+T_{3}^{n}\label{err_1}
 \end{empheq}
 \end{subequations}
 for $n=1,2,\cdots,N-1$,
where the truncation errors $T^{n}_{i}, i=1,2,3$ are given by
 \begin{equation*}\label{truerrs}
 \begin{aligned}
 T^{n}_{1}&\tps=\vPhi_{t}(t_{n+1})-\frac{\vPhi(t_{n+1})-\vPhi(t_{n})}{\Dt_{n+1}},\quad
  T^{n}_{2}&\tps=\varepsilon^2\Delta \vPhi(t_{n+1})-\varepsilon^2\Delta_{h}\vPhi(t_{n+1}),\\
T^{n}_{3}&\tps=\vPhi_{t}(t_{n+1})-\partial_{t}(\Pi_{2,n}\vPhi)(t_{n+1}).
\end{aligned}
\end{equation*}
Taking the discrete $L^{2}$ inner products of \eqref{err_2} and \eqref{err_1} with $2\Dt_{n+1}e^{*,n+1}$ and $2(e^{n+1}-e^{n})$, respectively, we obtain by
Lemma \ref{intpat} that
 \begin{subequations}\label{err1}
 \begin{align}
 &\dps \|e^{*,n+1}\|^{2}_{h}-\|e^{n}\|^{2}_{h}+2\Dt_{n+1}\varepsilon^{2}\|\nabla_{h} e^{*,n+1}\|^{2}_{h}+2\Dt_{n+1}S \|e^{*,n+1}\|^{2}_{h},\nonumber\\
&\quad=\dps 2\Dt_{n+1}\big<S e^{n}-S(\vPhi(t_{n+1})-\vPhi(t_{n}))+f(\vPhi(t_{n+1}))-f(\vPhi^{n})+T_{1}^{n}+T_{2}^{n},e^{*,n+1}\big>_{\Omega},\label{err1_2}\\[5pt]
 &\dps 2\big<F^{n+1}_{2}e,e^{n+1}-e^{n}\big>_{\Omega}+\varepsilon^{2}\big(\|\nabla_{h} e^{n+1}\|^{2}_{h}-\|\nabla_{h} e^{n}\|^{2}_{h}\big)+S\big( \|e^{n+1}\|^{2}_{h}-\|e^{n+1}\|^{2}_{h}\big)\nonumber\\
 &\quad\leq2\big<S e^{*,n+1}+f(\vPhi(t_{n+1}))-f(\vPhi^{*,n+1})+T_{2}^{n}+T_{3}^{n},e^{n+1}-e^{n}\big>_{\Omega}.\label{err1_1}
 \end{align}
 \end{subequations}
For \eqref{err1_1}, using the inequality \eqref{ideq1}, Cauchy-Schwarz inequality and Young's inequality, we have
\brr\label{rv1}
&\dps\frac{\gamma^{3/2}_{n+2}}{1+\gamma_{n+2}}\frac{\|e^{n+1}-e^{n}\|^{2}_{h}}{\Dt_{n+1}}-\frac{\gamma^{3/2}_{n+1}}{1+\gamma_{n+1}}\frac{\|e^{n}-e^{n-1}\|^{2}_{h}}{\Dt_{n}}+G(\gamma_{*},\gamma_{*})\frac{\|e^{n+1}-e^{n}\|^{2}_{h}}{\Dt_{n+1}}\\[11pt]
&\dps+\varepsilon^{2}\big(\|\nabla_{h} e^{n+1}\|^{2}_{h}-\|\nabla_{h} e^{n}\|^{2}_{h}\big)+S\big(\|e^{n+1}\|^{2}_{h}-\|e^{n+1}\|^{2}_{h}\big)\\[5pt]
&\quad\leq\dps C_3\Dt_{n+1}\big(\|e^{*,n+1}\|^{2}_{h}+\|f(\vPhi(t_{n+1}))-f(\vPhi^{*,n+1})\|^{2}_{h}+\|T^{n}_{2}\|^{2}_{h}+\|T^{n}_{3}\|^{2}_{h}\big)\\[5pt]
&\quad\quad\dps+G(\gamma_{*},\gamma_{*})\frac{\|e^{n+1}-e^{n}\|^{2}_{h}}{\Dt_{n+1}}\\[5pt]
&\quad\leq\dps C_3\max\{1,(C_2)^2\}\Dt_{n+1}\big(\|e^{*,n+1}\|^{2}_{h}+\|T^{n}_{2}\|^{2}_{h}+\|T^{n}_{3}\|^{2}_{h}\big)+G(\gamma_{*},\gamma_{*})\frac{\|e^{n+1}-e^{n}\|^{2}_{h}}{\Dt_{n+1}},
\err
where we have used the fact
$$\|f(\vPhi(t_{n+1}))-f(\vPhi^{*,n+1})\|^{2}_{h}\leq (C_2)^2\|e^{*,n+1}\|^{2}_{h}$$ derived from \eqref {equ1}.
Thus we deduce that
\brr\label{eqt1}
&\dps\frac{\gamma^{3/2}_{n+2}}{1+\gamma_{n+2}}\frac{\|e^{n+1}-e^{n}\|^{2}_{h}}{\Dt_{n+1}}-\frac{\gamma^{3/2}_{n+1}}{1+\gamma_{n+1}}\frac{\|e^{n}-e^{n-1}\|^{2}_{h}}{\Dt_{n}}+\varepsilon^{2}\big(\|\nabla_{h} e^{n+1}\|^{2}_{h}-\|\nabla_{h} e^{n}\|^{2}_{h}\big)\\[11pt]
&\dps+S\big(\|e^{n+1}\|^{2}_{h}-\|e^{n}\|^{2}_{h}\big)
\leq\;\dps C_3\max\{1,(C_2)^2\}\Dt_{n+1}\big(\|e^{*n+1}\|^{2}_{h}+\|T^{n}_{2}\|^{2}_{h}+\|T^{n}_{3}\|^{2}_{h}\big).
\err
In a similar way, we can obtain the following estimate from \eqref{err1_2}
\brr\label{eqnn1}
&\dps \|e^{*,n+1}\|^{2}_{h}-\|e^{n}\|^{2}_{h}+2\Dt_{n+1}\varepsilon^{2}\|\nabla_{h} e^{*,n+1}\|^{2}_{h}+2\Dt_{n+1}S \|e^{*,n+1}\|^{2}_{h}\\[5pt]
&\quad\leq\dps C_4\Dt_{n+1}^{2}\big(\|e^{n}\|^{2}_{h}+\|\vPhi(t_{n+1})-\vPhi(t_{n})\|^{2}_{h}+\|f(\vPhi(t_{n+1}))-f(\vPhi^{n})\|^{2}_{h}\\[5pt]
&\qquad\dps+\|T_{1}^{n}\|^{2}_{h}+\|T_{2}^{n}\|^{2}_{h}\big)+\tps\frac12{\|e^{*,n+1}\|^{2}_{h}}\\[5pt]
&\quad\leq\dps C_5\Dt_{n+1}^{2}\big(\|e^{n}\|^{2}_{h}+\Dt_{n+1}^{2}\|\phi_{t}\|^{2}_{L^{\infty}(0,T;L^{\infty}(\Omega))}+\|T_{1}^{n}\|^{2}_{h}+\|T_{2}^{n}\|^{2}_{h}\big)+\tps\frac12{\|e^{*,n+1}\|^{2}_{h}},\\
\err
where we have used the fact
\bry
\|f(\vPhi(t_{n+1}))-f(\vPhi^{n})\|^{2}_{h}\leq&\dps\|f(\vPhi(t_{n+1}))-f(\vPhi(t_{n}))\|^{2}_{h}+\|f(\vPhi(t_{n}))-f(\vPhi^{n})\|^{2}_{h}\\[5pt]
\leq&\dps (C_2)^2\big(\|\vPhi(t_{n+1})-\vPhi(t_{n})\|^{2}_{h}+\|e^{n}\|^{2}_{h}\big)\\[5pt]
\leq&\dps (C_2)^2\big(\Dt_{n+1}^{2}\|\phi_{t}\|^{2}_{L^{\infty}(0,T;L^{\infty}(\Omega))}+\|e^{n}\|^{2}_{h}\big).
\ery
Then it follows from \eqref{eqnn1} that
\brr\label{eqt2}
&\dps \|e^{*,n+1}\|^{2}_{h}
\leq\dps 2\|e^{n}\|^{2}_{h}+2C_5\Dt^{2}_{n+1}\big(\|e^{n}\|^{2}_{h}+\Dt^{2}_{n+1}\|\phi_{t}\|^{2}_{L^{\infty}(0,T;L^{\infty}(\Omega))}+\|T_{1}^{n}\|^{2}_{h}+\|T_{2}^{n}\|^{2}_{h}\big).
\err
Combining with \eqref{eqt1} and \eqref{eqt2}, gives
\brr\label{eqt3}
&\dps\frac{\gamma^{3/2}_{n+2}}{1+\gamma_{n+2}}\frac{\|e^{n+1}-e^{n}\|^{2}_{h}}{\Dt_{n+1}}-\frac{\gamma^{3/2}_{n+1}}{1+\gamma_{n+1}}\frac{\|e^{n}-e^{n-1}\|^{2}_{h}}{\Dt_{n}}+\varepsilon^{2}\big(\|\nabla_{h} e^{n+1}\|^{2}_{h}-\|\nabla_{h} e^{n}\|^{2}_{h}\big)\\[11pt]
&\dps+S\big( \|e^{n+1}\|^{2}_{h}-\|e^{n}\|^{2}_{h}\big)\\[9pt]
&\quad\leq\dps C_6\Dt_{n+1}\big(\|e^{n}\|^{2}_{h}+\Dt^{4}_{n+1}\|\phi_{t}\|^{2}_{L^{\infty}(0,T;L^{\infty}(\Omega))}+\Dt^{2}_{n+1}\|T_{1}^{n}\|^{2}_{h}+\|T^{n}_{2}\|^{2}_{h}+\|T^{n}_{3}\|^{2}_{h}\big).
\err
For the truncation errors $T^{n}_{i}, i=1,2,3$, we have the following estimates (see \cite{LSR19,LTZ20}):
\brr\label{equ8}
\|T^{n}_{1}\|^{2}_{h}\leq &\;\dps C_7\Dt_{n+1}^{2}\|\phi\|^{2}_{W^{2,\infty}(0,T;L^{\infty}(\Omega))},\quad
\|T^{n}_{2}\|^{2}_{h}\leq \;\dps C_8h^{4}\|\phi\|^{2}_{L^{\infty}(0,T;W^{4,\infty}(\Omega))},\\[5pt]
\|T^{n}_{3}\|^{2}_{h}\leq &\;\dps C_9(\Dt_{n}+\Dt_{n+1})^{4}\|\phi\|^{2}_{W^{3,\infty}(0,T;L^{\infty}(\Omega))}.
\err
Thus, summing up the inequality \eqref{eqt3} from 1 to $n$ gives
\brr\label{equ9}
&\dps\frac{\gamma^{3/2}_{n+2}}{1+\gamma_{n+2}}\frac{\|e^{n+1}-e^{n}\|^{2}_{h}}{\Dt_{n+1}}
+\varepsilon^{2}\|\nabla_{h}e^{n+1}\|^{2}_{h}+S\|e^{n+1}\|^{2}_{h}\\[11pt]
&\quad\leq\dps\frac{\gamma^{\frac{3}{2}}_{2}}{1+\gamma_{2}}\frac{\|e^{1}\|^{2}_{h}}{\Dt_{1}}+\varepsilon^{2}\|\nabla_{h} e^{1}\|^{2}_{h} +S\|e^{1}\|^{2}_{h}+C_{10}\sum_{k=1}^{n}\Dt_{k+1}\|e^{k}\|^{2}_{h}\\
&\qquad\dps+C_{11}\big(\Dt^{4}\|\phi\|^{2}_{W^{3,\infty}(0,T;L^{\infty}(\Omega))}+h^{4}\|\phi\|^{2}_{L^{\infty}(0,T;W^{4,\infty}(\Omega))}\big).
\err
For the case of  $n=0$, the corresponding error equation (by BDF1) reads as
 \beq
 \frac{e^{1}}{\Dt_{1}}-\varepsilon^{2}\Delta_{h} e^{1}+S e^{1}=f(\vPhi(t_{1}))-f(\vPhi^{0})+T^{0}_{2}+T^{0}_{3}.
  \eeq
  Similar to the arguments for the case $n\geq1$, the following estimate can be derived under the assumption $\Dt_{1}\leq C_1 \Dt^{4/3}$:
\bq\label{equa2}
 \dps\frac{\|e^{1}\|^{2}_{h}}{\Dt_{1}}+\varepsilon^{2}\|\nabla e^{1}\|^{2}_{h}+S\|e^{1}\|^{2}_{h}\leq C_{12} \Dt_{1}\big(\Dt_{1}^{2}+h^{4}\big)\leq C_{12}\max\{1,(C_1)^3\}\big(\Dt^{4}+h^{4}\big).
\eq
Combining \eqref{equ9} and \eqref{equa2}  and using the discrete Gronwall's lemma,  we then obtain the desired estimate \eqref{equa4_1}.
 \end{proof}
\begin{remark}
It often imposes a further  restriction on the time step size when using the Gronwall's inequality for the error analysis.
However, in the above proof of Theorem \ref{th1}, we note that
the term $G(\gamma_{*},\gamma_{*})\frac{\|e^{n+1}-e^{n}\|^{2}_{h}}{\Dt_{n+1}}$ on the right-hand side of the error inequality \eqref{rv1} can be eliminated by a term from  the left-hand side of the equation.
Consquently, we are able to  obtain the error inequality \eqref{equ9}, which only contains  the norm terms of $e^{n+1}$ with positive coefficients on the left-hand side.
Thus, there is no further time step restriction from the use of the Gronwall's inequality in our error analysis.
\end{remark}

With the help of the MBP property (Theorem \ref{thmMBP}) and the discrete $H^1$ error estimate (Theorem \ref{th1}), we are able to achieve
the energy stability property of the BDF2 scheme \eqref{BDF_2}.

\begin{theorem}\label{st_the}
Under the assumption of Theorem \ref{th1}, the BDF2 scheme \eqref{BDF_2}  in the constant mobility case is energy stable in the sense that
\bq\label{eq2_1}
\dps E^{n+1}_{h}-E^{n}_{h}\leq C(h^{4}+\Dt^{2})
\eq
for all $0\leq n\leq N-1$, where  the modified discrete  energy $E^{n}_{h}$ is defined by
\beq
\dps E^{n}_{h}=E_{h}(\vPhi^{n})+
\frac{\gamma_{n+1}^{3/2}}{1+\gamma_{n+1}}\frac{\|\vPhi^{n}-\vPhi^{n-1}\|^{2}_{h}}{2\tau_{n}}.
\eeq
\end{theorem}
\begin{proof}
Taking the discrete $L^{2}$-inner product of \eqref{BDF_2} with $\vPhi^{n+1}-\vPhi^{n}$, we get  that
\brr\label{eqt4}
&\dps \big<F^{n+1}_{2}\vPhi,\vPhi^{n+1}-\vPhi^{n}\big>_{\Omega}+\varepsilon^{2}\big<\nabla_{h}\vPhi^{n+1},\nabla_{h}(\vPhi^{n+1}-\vPhi^{n})\big>_{\Omega}+\big<f(\vPhi^{n+1}),\vPhi^{n+1}-\vPhi^{n}\big>_{\Omega}\\[5pt]
&\quad=\dps\big<f(\vPhi^{n+1})-f(\vPhi^{*,n+1}),\vPhi^{n+1}-\vPhi^{n}\big>_{\Omega}-S\big<\vPhi^{n+1}-\vPhi^{*,n+1},\vPhi^{n+1}-\vPhi^{n}\big>_{\Omega}\\[5pt]
&\quad\leq\dps \frac{(C_{2})^{2}+S^{2}}{2}\|\vPhi^{n+1}-\vPhi^{*,n+1}\|^{2}_{h}+\|\vPhi^{n+1}-\vPhi^{n}\|^{2}_{h}\\[5pt]
&\quad\leq\tps\max\Big\{\frac12((C_{2})^{2}+S^{2}),1,\|\phi_{t}\|_{L^{\infty}(0,T;L^{\infty}(\Omega))}\Big\}\big(\|e^{n}\|^{2}_{h}+\|e^{n+1}\|^{2}_{h}+\|e^{*,n+1}\|^{2}_{h}+\Dt^{2}_{n+1}\big),
 \err
 where we have used the following inequalities
 \bry
 \|\vPhi^{n+1}-\vPhi^{*,n+1}\|^{2}_{h}=\;&\dps\|\vPhi^{n+1}-\vPhi(t_{n+1})+\vPhi(t_{n+1})-\vPhi^{*,n+1}\|^{2}_{h}\\[5pt]
 \leq \;&\|e^{n+1}\|^{2}_{h}+\|e^{*,n+1}\|^{2}_{h},\\[5pt]
 \|\vPhi^{n+1}-\vPhi^{n}\|^{2}_{h}=\;&\dps \|\vPhi^{n+1}-\vPhi(t_{n+1})+\vPhi(t_{n+1})-\vPhi(t_{n})+\vPhi(t_{n})-\vPhi^{n}\|^{2}_{h}\\[5pt]
 \leq\;&\dps\|e^{n+1}\|^{2}_{h}+\Dt^{2}_{n+1}\|\phi_{t}\|_{L^{\infty}(0,T;L^{\infty}(\Omega))}+\|e^{n}\|^{2}_{h}.
 \ery
 Noting that
 \bry
& a(a-b)=\dps\frac{1}{2}\big(a^{2}-b^{2}+(a-b)^{2}\big), \quad a,b\in \mathbb{R},\\[8pt]
& \big<F(\vPhi^{n+1})-F(\vPhi^{n}),1\big>_{\Omega}\leq\dps \big<f(\vPhi^{n+1}), \vPhi^{n+1}-\vPhi^{n}\big>_{\Omega}+\frac{1}{2}\|\vPhi^{n+1}-\vPhi^{n}\|^{2}_{h},
 \ery
and using  \eqref{ideq1} and \eqref{eqt4}, we can derive
 \brr\label{ttt}
E^{n+1}_{h}-E^{n}_{h}\leq\tps\max\big\{\frac12((C_{2})^{2}+S^{2}),\frac{3}{2},\frac{3}{2}\|\phi_{t}\|_{L^{\infty}(0,T;L^{\infty}(\Omega))}\big\}\\[5pt]
\cdot\big(\|e^{n}\|^{2}_{h}+\|e^{n+1}\|^{2}_{h}+\|e^{*,n+1}\|^{2}_{h}+\Dt^{2}_{n+1}\big).
 \err
Combining  with \eqref{equa4_1}, \eqref{eqt2}, and \eqref{ttt}, we then obtain \eqref{eq2_1}.
\end{proof}

\begin{remark}
For the quasi-uniform temporal mesh, there exits a finite constant $\beta$ such that $\dps\max\limits_{1\leq n\leq N}\Dt_{n}/\min\limits_{1\leq n\leq N}\Dt_{n}\leq\beta$ and thus
$\Dt\leq\frac{\beta T}{N}$.
When $\Dt$ is sufficient small and $h = O(\sqrt{\Dt})$, we can obtain
\beq
\dps E_{h}(\vPhi^{n})\leq E^{n}_{h}\leq  E^{1}_{h}+C\leq E_{h}(\vPhi^{0})+C,\quad\forall\,1\leq n\leq N
\eeq
for the BDF2 scheme \eqref{BDF_2} in the constant mobility case.
\end{remark}

\begin{remark}\label{rek1}
The inequality \eqref{ideq1} plays an important role in the above error and  energy stability analysis. 
Unfortunately, we have not been able to prove  a similar result as \eqref{ideq1}
 for the estimate of $\big<(\Lambda^{*,n+1})^{-1} F^{n+1}_{2}\vPhi,\vPhi^{n+1}-\vPhi^{n}\big>_{\Omega}$  in the case of variable mobility.
Thus the results in Theorems \ref{th1} and \ref{st_the} could not be applied to the variable mobility case,  and deeper analysis for this issue certainly needs more efforts.
\end{remark}

\subsection{Error estimate in the $L^{\infty}$ norm for the general mobility case}

In this subsection, we study the discrete $L^{\infty}$ error estimate of the BDF2 scheme \eqref{BDF_2}  for the Allen-Cahn equation with  a general mobility $M(\phi)$.
Let us define
\beq
F_{2}\vPhi(t_{n+1})=b^{n}_{0}\big(\vPhi(t_{n+1})-\vPhi(t_{n})\big)+b^{n}_{1}\big(\vPhi(t_{n})-\vPhi(t_{n-1})\big)
\eeq
for $1\leq n\leq N-1$ and $ \Lambda(\vPhi(t_{n}))=\mbox{diag}(M(\vPhi(t_{n})))$ for $0\leq n\leq N-1$.

\begin{lemma}\label{lem1_2}
Assume that $\{g^{k}\}_{k=0}^{N-1}$ and $\{\omega^{k}\}_{k=0}^{N}$ are two non-negative sequences and there exist some  constants $\zeta>0$ and $\lambda\in(0,1)$ such that
\bq\label{rev1}
\sum_{k=1}^{n+1}d^{n}_{n-k+1}\delta_{\Dt}\omega^{k}\leq \zeta\sum_{k=0}^{n}\lambda^{n-k}\omega^{k}+g^{n}, \quad \forall\,0\leq n\leq N-1,
\eq
where the discrete kernels $\{d^{n}_{k}\}_{k=0}^{n}$ are defined in \eqref{eqn1_1}. Then it holds
\bq\label{rv2}
\omega^{n+1}\leq\exp\Big(\frac{\zeta t_{n+1}}{1-\lambda}\Big)\Big(\omega^{0}+\sum_{k=0}^{n}\frac{g^{k}}{b^{k}_{0}}\Big).
\eq
\end{lemma}
The proof  of this lemma is similar to that of  Lemma 5.1 in \cite{LTZ20} and Theorem 3.1 in \cite{LMZ19} by using the technique of the discrete complementary convolution kernels of $\{d^{n}_{k}\}_{k=0}^{n}$. We omit it here and  leave it for the interested readers. Comparing with Lemma 5.1 in \cite{LTZ20} and Theorem 3.1 in \cite{LMZ19}, there is no term $\omega^{n+1}$ on the right-hand side of the condition \eqref{rev1}. Then the time step restriction required in \cite{LTZ20} and \cite{LMZ19} for the result \eqref{rv2} can be removed.

\begin{theorem}\label{th3}
Assume  that $0<\gamma_{n}\leq\gamma_{*}<1+\sqrt{2}$ for all $1\leq n \leq N-1$,  $M(\cdot)\in C^{1}(\mathbb{R})$, the stabilizing parameter satisfies \eqref{eqn1_3},
and the time step sizes satisfy \eqref{est_tau3} and \eqref{constr_tau1}. Let $\gamma_{N+1}$ be any number in $(0,\gamma_*$). In addition, assume  $\phi\in W^{3,\infty}(0,T;L^{\infty}(\Omega))\cap L^{\infty}(0,T;W^{4,\infty}(\Omega))$. Then it holds  for the BDF2 scheme \eqref{BDF_2} in the general mobility case that
\bq
\|e^{n+1}\|_{\infty}
\leq\frac{C_1 t_{n+1}}{1-\eta_{*}}\exp\Big(\frac{{C}_{2} t_{n+1}}{1-\eta_{*}}\Big)\big(\Dt^{2}\|\phi\|_{W^{3,\infty}(0,T;L^{\infty}(\Omega))}+h^{2}\|\phi\|_{L^{\infty}(0,T;W^{4,\infty}(\Omega))}\big)
\eq
for all $0\leq n\leq N-1$.
\end{theorem}
\begin{proof}
From \eqref{prob}, we deduce that the exact solution $\vPhi$ satisfies the following equation: for any $1\leq n\leq N-1$,
\begin{equation}\label{FforExact}
F_{2}\vPhi(t_{n+1})+\Lambda(\vPhi(t_{n+1}))\big(-\varepsilon^{2}D_{h}\vPhi(t_{n+1})+F'(\vPhi(t_{n+1}))\big)+\mathcal{T}^{n}_{2}+\mathcal{T}^{n}_{3}=0,
\end{equation}
where
$$  \mathcal{T}^{n}_{2}\dps=\Lambda(\vPhi(t_{n+1}))(-\varepsilon^2\Delta \vPhi(t_{n+1})+\varepsilon^2D_{h}\vPhi(t_{n+1})),\quad
 \mathcal{T}^{n}_{3}\dps=\vPhi_{t}(t_{n+1})-\partial_{t}(\Pi_{2,n}\vPhi)(t_{n+1}).$$
It is easy to verify that
\begin{eqnarray*}
&&\|\mathcal{T}^{n}_{2}\|_{\infty}\leq C_3 h^{2}\|\phi\|_{L^{\infty}(0,T;W^{4,\infty}(\Omega))},\quad
\|\mathcal{T}^{n}_{3}\|_{\infty}\leq C_4(\Dt_{n}+\Dt_{n+1})^{2}\|\phi\|_{W^{3,\infty}(0,T;L^{\infty}(\Omega))}.
\end{eqnarray*}
Subtracting \eqref{BDF2_tsor2} from \eqref{FforExact}, we derive the error equation of $e^{n+1}$ as
\brr\label{eqn1_7}
&\dps F^{n+1}_{2}e+Se^{n+1}-\varepsilon^{2}\Lambda^{*,n+1}D_{h}e^{n+1}\\[2pt]
&\quad=\dps Se^{*,n+1}-\Lambda^{*,n+1}\big(F'(\vPhi^{*,n+1})-F'(\vPhi(t_{n+1}))\big)
-(\Lambda^{*,n+1}-\Lambda(\vPhi(t_{n+1})))\\[2pt]
&\qquad\dps\big(-\varepsilon^{2}D_{h}\vPhi(t_{n+1})+F'(\vPhi(t_{n+1}))\big)+\mathcal{T}^{n}_{2}+\mathcal{T}^{n}_{3}=:I^{n}.
\err
Since $F(\rho)=(1-\rho^{2})^{2}/4$, we have
$\max_{\rho\in[-1,1]}F'(\rho)=\frac{2}{3\sqrt{3}}$ and  $\max_{\rho\in[-1,1]}F''(\rho)=2.$
Therefore, we can get
$$
\|\Lambda^{*,n+1}-\Lambda(\vPhi(t_{n+1}))\|_{\infty}\leq\dps\max_{\rho\in[-1,1]}\big|M'(\rho)\big|\| e^{*,n+1}\|_{\infty},
$$
and thus
\brr\label{eqn1_12}
\|I^{n}\|_{\infty}\leq&\;\dps S\| e^{*,n+1}\|_{\infty}+2L\| e^{*,n+1}\|_{\infty}+\Big(\varepsilon^{2}\|\phi\|_{L^{\infty}(0,T,W^{2,\infty}(\Omega))}+\tps\frac{2}{3\sqrt{3}}\Big)\\
&\dps\|\Lambda^{*,n+1}-\Lambda(\vPhi(t_{n+1}))\|_{\infty}+\|\mathcal{T}^{n}_{2}\|_{\infty}+\|\mathcal{T}^{n}_{3}\|_{\infty}\\[10pt]
\leq&\;\dps {C}_{4}\| e^{*,n+1}\|_{\infty}+C_5\big[\Dt^{2}\|\phi\|_{W^{3,\infty}(0,T;L^{\infty}(\Omega))}+h^{2}\|\phi\|_{L^{\infty}(0,T;W^{4,\infty}(\Omega))}\big],
\err
where $L$ is defined in Lemma \ref{lemm1},
\bq\label{eqn1_9}
{C}_{4}=S+2L+\max_{\rho\in[-1,1]}\big|M'(\rho)\big|\Big(\varepsilon^{2}\|\phi\|_{L^{\infty}(0,T,W^{2,\infty}(\Omega))}+\tps\frac{2}{3\sqrt{3}}\Big)=S+2L+{C}_{6},
\eq
and we have used the fact
\bq\label{eqn1_10}
\|\varepsilon^{2}D_{h}\vPhi(t_{n+1})\|_{\infty}=\|\varepsilon^{2}\Delta\vPhi(t_{n+1})+\mathcal{T}^{n}_{2}\|_{\infty}\leq\varepsilon^{2}\|\phi\|_{L^{\infty}(0,T,W^{2,\infty}(\Omega))}+\|\mathcal{T}^{n}_{2}\|_{\infty}.
\eq
Following the similar process of deriving \eqref{eqn1_7}, we can easily obtain the error equation of $ e^{*,n+1}$ from \eqref{prob} and \eqref{BDF2_tsor1} as:
\brr\label{eqn1_8}
&\dps\frac{ e^{*,n+1}- e^{n}}{\Dt_{n+1}}+S e^{*,n+1}-\varepsilon^{2}\Lambda^{n}D_{h} e^{*,n+1}\\
&\quad=\dps S e^{n}-S(\vPhi(t_{n+1})-\vPhi(t_{n}))-\Lambda^{n}\big[F'(\vPhi^{n})-F'(\vPhi(t_{n+1}))\big]\\[5pt]
&\qquad\dps-\big[\Lambda^{n}-\Lambda(\vPhi(t_{n+1}))\big]\big[-\varepsilon^{2}D_{h}\vPhi(t_{n+1})+F'(\vPhi(t_{n+1}))\big]+\mathcal{T}^{n}_{1}+\mathcal{T}^{n}_{2},
\err
where
$\mathcal{T}^{n}_{1}=\vPhi_{t}(t_{n+1})-\frac{\vPhi(t_{n+1})-\vPhi(t_{n})}{\Dt_{n+1}}$
satisfies
\begin{equation*}\label{ts_trunc2}
\|\mathcal{T}^{n}_{1}\|_{\infty}\leq C_7\Dt_{n+1}\|\phi\|_{W^{2,\infty}(0,T;L^{\infty}(\Omega))}.
\end{equation*}
Noting that
\bry
\|F'(\vPhi^{n})-F'(\vPhi(t_{n+1}))\|_{\infty}=\;&\dps\|F'(\vPhi^{n})+F'(\vPhi(t_{n}))\|_{\infty}+\|F'(\vPhi(t_{n}))-F'(\vPhi(t_{n+1}))\|_{\infty}\\[5pt]
\leq\;&\dps 2\big[\| e^{n}\|_{\infty}+\|\phi\|_{W^{1,\infty}(0,T;L^{\infty}(\Omega))}\Dt_{n+1}\big],\\[5pt]
\big\|\Lambda^{n}-\Lambda(\vPhi(t_{n+1}))\big\|_{\infty}\leq\;&\dps\big\|\Lambda^{n}-\Lambda(\vPhi(t_{n}))\big\|_{\infty}+\big\|\Lambda(\vPhi(t_{n}))-\Lambda(\vPhi(t_{n+1}))\big\|_{\infty}\\[5pt]
\leq\;&\dps\max_{\rho\in[-1,1]}\big|M'(\rho)\big|\big[\| e^{n}\|_{\infty}+\Dt_{n+1}\|\phi\|_{W^{1,\infty}(0,T;L^{\infty}(\Omega))}\big].
\ery
Multiplying \eqref{eqn1_8} with $\Dt_{n+1}$, and combining it with \eqref{eqn1_9} and \eqref{eqn1_10}, we derive that
\bry
(1+S\Dt_{n+1})\| e^{*,n+1}\|_{\infty}\leq\;&\dps\big\| e^{*,n+1}+S\Dt_{n+1} e^{*,n+1}-\varepsilon^{2}\Dt_{n+1}\Lambda^{n}D_{h} e^{*,n+1}\big\|_{\infty}\\[5pt]
\leq\;&\dps \big(1+\Dt_{n+1}(S+2L+{C}_{6})\big)\| e^{n}\|_{\infty}\\[5pt]
\;&\dps+\,C_8\big(\Dt^{2}_{n+1}\|\phi\|_{W^{2,\infty}(0,T;L^{\infty}(\Omega))}+\Dt_{n+1} h^{2}\|\phi\|_{L^{\infty}(0,T;W^{4,\infty}(\Omega))}\big).
\ery
Therefore, we obtain
\bq\label{eqn1_13}
\| e^{*,n+1}\|_{\infty}\leq {C}_{9}\| e^{n}\|_{\infty}+C_{10}\big(\Dt^{2}_{n+1}\|\phi\|_{W^{2,\infty}(0,T;L^{\infty}(\Omega))}+\Dt_{n+1} h^{2}\|\phi\|_{L^{\infty}(0,T;W^{4,\infty}(\Omega))}\big),
\eq
where $
{C}_{9}=1+(2L+{C}_{6})\Dt_{n+1}\geq\tps\frac{1+\Dt_{n+1}(S+2L+{C}_{6})}{1+S\Dt_{n+1}}.$

Define $\overline{ e}^{n+1}= e^{n+1}-\eta_{*} e^{n}$ for $0\leq n\leq N-1$ with $\overline{ e}^{0}= e^{0}=0$. Then we have
\bq\label{eqn1_15}
\overline{ e}^{1}= e^{1}, \quad  e^{n+1}=\dps\sum_{k=0}^{n+1}\eta_*^{n+1-k}\overline{ e}^{k}, \quad 1\leq n\leq N-1,
\eq
and it is easy to check that
\beq
\|\overline{ e}^{1}\|_{\infty}=\| e^{1}\|_{\infty}\leq C_{11}\big(\Dt^{2}_{1}\|\phi\|_{W^{2,\infty}(0,T;L^{\infty}(\Omega))}+\Dt_{1} h^{2}\|\phi\|_{L^{\infty}(0,T;W^{4,\infty}(\Omega))}\big)
\eeq
with the BDF1 scheme as the starting step.
Using the similar process to derive \eqref{eqn_1} from \eqref{BDF2_tsor2}, we can obtain the following equation for $\overline{ e}^{n+1}$ from \eqref{eqn1_7}:
 \brr\label{eqn_11}
 \big((d_{0}^{n}+S)I-\varepsilon^{2}\Lambda^{*,n+1}D_{h}\big)\overline{ e}^{n+1}=&
\dps\sum_{k=0}^{n}Q^{n}_{n-k}\overline{ e}^{k}+I^{n},
 \err
where $Q^{n}_{k}$ is defined in \eqref{eqn1_2} with $\eta=\eta_{*}$.
Consequently, we can use Lemma \ref{lem_1} and \eqref{eqn_11} to get
\bry
 d_{0}^{n}\|\overline{ e}^{n+1}\|_{\infty}\leq\;&\dps\big\|\big((d_{0}^{n}+S)I-\varepsilon^{2}\Lambda^{*,n+1}D_{h}\big)\overline{ e}^{n+1}\big\|_{\infty}\\[5pt]
 \leq\;&\tps\sum_{k=0}^{n}(d^{n}_{n-k}-d^{n}_{n-k+1}-S\eta^{n+1-k}_{*})\|\overline{ e}^{k}\|_{\infty}+\|I^{n}\|_{\infty}\\[5pt]
\leq\;&\tps\sum_{k=0}^{n}(d^{n}_{n-k}-d^{n}_{n-k+1})\|\overline{ e}^{k}\|_{\infty}+\|I^{n}\|_{\infty}.
 \ery
Rewriting the above inequality gives
$\dps\sum_{k=1}^{n+1}d^{n}_{n-k+1}\delta_{\Dt}\|\overline{ e}^{k}\|_{\infty}\leq \dps\|I^{n}\|_{\infty}.$
Combining it with \eqref{eqn1_12}, \eqref{eqn1_13} and \eqref{eqn1_15}, we obtain
\be\label{eqt1_1}
\begin{aligned}
\sum_{k=1}^{n+1}d^{n}_{n-k+1}\delta_{\Dt}\|\overline{ e}^{k}\|_{\infty}\leq&\;\dps{C}_{2}\| e^{k}\|_{\infty}+C_{12}\big(\Dt^{2}\|\phi\|_{W^{3,\infty}(0,T;L^{\infty}(\Omega))}+h^{2}\|\phi\|_{L^{\infty}(0,T;W^{4,\infty}(\Omega))}\big)\\
\leq&\; \tps{C}_{2}\sum_{k=1}^{n}\eta_{*}^{n-k}\|\overline{ e}^{k}\|_{\infty}+C_{12}\big(\Dt^{2}\|\phi\|_{W^{3,\infty}(0,T;L^{\infty}(\Omega))}\\
&\dps+h^{2}\|\phi\|_{L^{\infty}(0,T;W^{4,\infty}(\Omega))}\big)
\end{aligned}
\ee
with ${C}_{2}={C}_{4}{C}_{8}$.
 Next, it follows from Lemma \ref{lem1_2} that
 \beq
\begin{aligned}
\|\overline{ e}^{n+1}\|_{\infty}\leq &\;\dps C_{1}\exp\Big(\frac{{C}_{2} t_{n+1}}{1-\eta_{*}}\Big)\big(\Dt^{2}\|\phi\|_{W^{3,\infty}(0,T;L^{\infty}(\Omega))}+h^{2}\|\phi\|_{L^{\infty}(0,T;W^{4,\infty}(\Omega))}\big)\tps\sum_{k=0}^{n}\frac{1}{b^{k}_{0}}\\
\leq&\;\dps C_{1} t_{n+1}\exp\Big(\frac{{C}_{2} t_{n+1}}{1-\eta_{*}}\Big)\big(\Dt^{2}\|\phi\|_{W^{3,\infty}(0,T;L^{\infty}(\Omega))}+h^{2}\|\phi\|_{L^{\infty}(0,T;W^{4,\infty}(\Omega))}\big),
\end{aligned}
\eeq
where we have used  the fact that
$ \tps\frac{1}{b^{0}_{0}}=\Dt_{1}$ and $\tps\frac{1}{b^{k}_{0}}=\frac{1+\gamma_{k+1}}{1+2\gamma_{k+1}}\Dt_{k+1}\leq \Dt_{k+1}$ for $1\leq k\leq n.$
Finally, we obtain
\beq
\begin{aligned}
\| e^{n+1}\|_{\infty}\leq\;&\dps \sum_{k=1}^{n+1}\eta^{n+1-k}_{*}\|\overline{ e}^{k}\|_{\infty}\\
\leq\;&\dps \frac{C_1 t_{n+1}}{1-\eta_{*}}\exp\Big(\frac{{C}_{2} t_{n+1}}{1-\eta_{*}}\Big)\big(\Dt^{2}\|\phi\|_{W^{3,\infty}(0,T;L^{\infty}(\Omega))}+h^{2}\|\phi\|_{L^{\infty}(0,T;W^{4,\infty}(\Omega))}\big),
\end{aligned}
\eeq
 which completes the proof.
 \end{proof}
\begin{remark}\label{rek2}
Similar to the derivation of linear BDF2 scheme \eqref{BDF_2}, one can also construct a linear second order in time scheme with variable time step sizes based on the Crank-Nicolson formulation as follows:
given $\Phi^0=\Pi_{\mathcal{C}_h}\phi_0$, and for $n=1,2\cdots,N-1$, find $\Phi^{n+1}\in\mathcal{C}_{h}$ such that
	\begin{subequations}\label{CN_2}
	      \begin{empheq}[left=\empheqlbrace]{align}
		&\Phi^{n+\frac{1}{2}} = {\rm BDF1}(\Phi^n,\tau_{n+1}/2),\label{eqn41}\\	&\frac{\Phi^{n+1}-\Phi^{n}}{\Dt_{n+1}}-\varepsilon^{2}M(\Phi^{n+\frac{1}{2}})\Delta_{h}\frac{\Phi^{n+1}+\Phi^{n}}{2}+f(\Phi^{n+\frac{1}{2}})+S\Big(\frac{\Phi^{n+1}+\Phi^{n}}{2}-\Phi^{n+\frac{1}{2}}\Big)=0,\label{eqn4}
		\end{empheq}
       \end{subequations}
The above theoretical analysis for the BDF2 scheme \eqref{BDF_2} can be applied to the Crank-Nicolson scheme \eqref{CN_2} to derive similar results obtained for the BDF2 scheme \eqref{BDF_2}, including the conditional MBP preserving and corresponding error estimates.
\end{remark}

\section{Numerical results}
\setcounter{equation}{0}

In this section we perform various experiments on the  Allen-Cahn equation \eqref{prob} to numerically validate the theoretical results of the proposed BDF2 scheme \eqref{BDF_2} in terms of accuracy and preservation of the MBP.  The homogenous Neumann boundary condition is always imposed.

\subsection{Test of temporal convergence}

We consider   two types of mobility functions: one is the constant mobility  $M(\phi)\equiv1$ and  the other is the nonlinear degenerate mobility $M(\phi)=1-\phi^{2}$.
We choose $\Omega = (0,1)^2$, $\varepsilon=0.1$, the initial value
$$\phi_0(x,y)=0.1(\cos3x\cos2y+\cos5x\cos5y),$$ and the terminal time $T=1$. The stabilizing parameter is set to be $S=2$ to satisfy the requirement \eqref{eqn1_3}
for both mobility functions.

\begin{figure*}[!t]
\centerline{\includegraphics[scale=0.4]{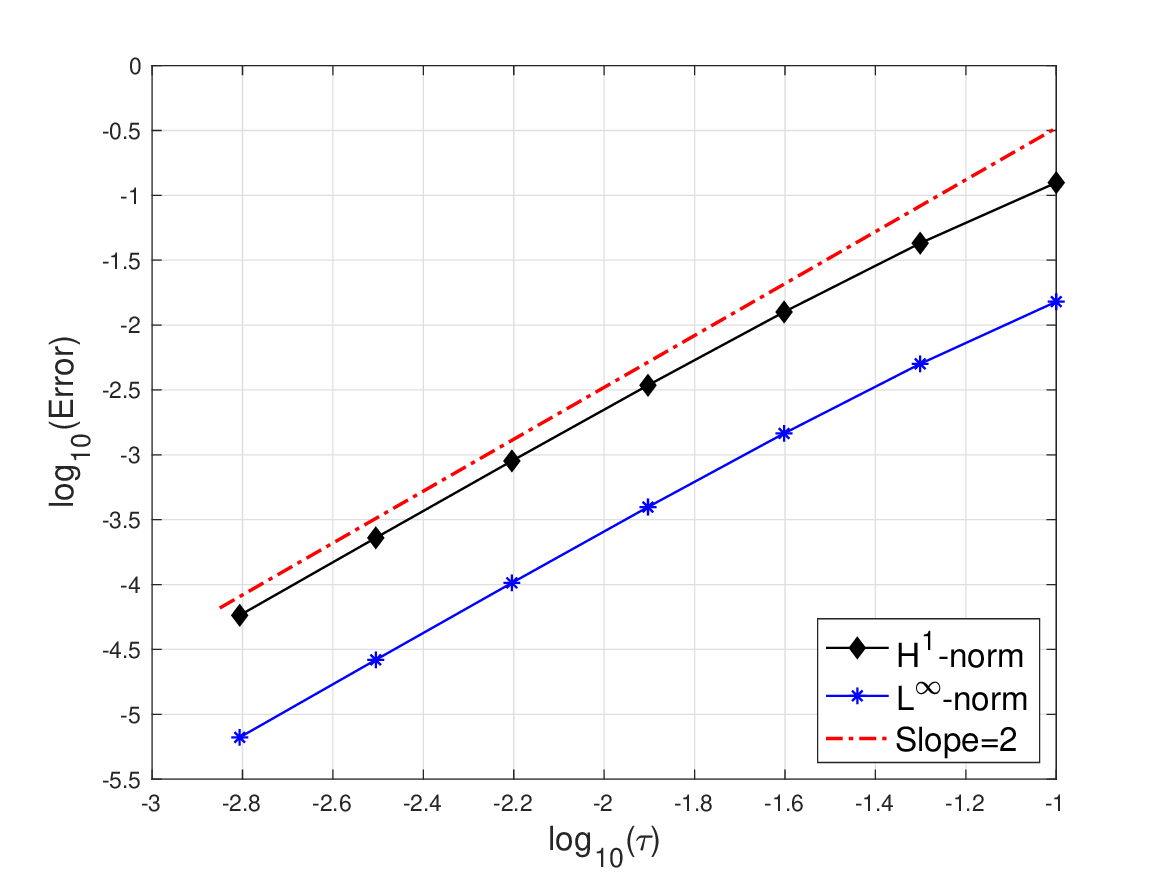}\includegraphics[scale=0.4]{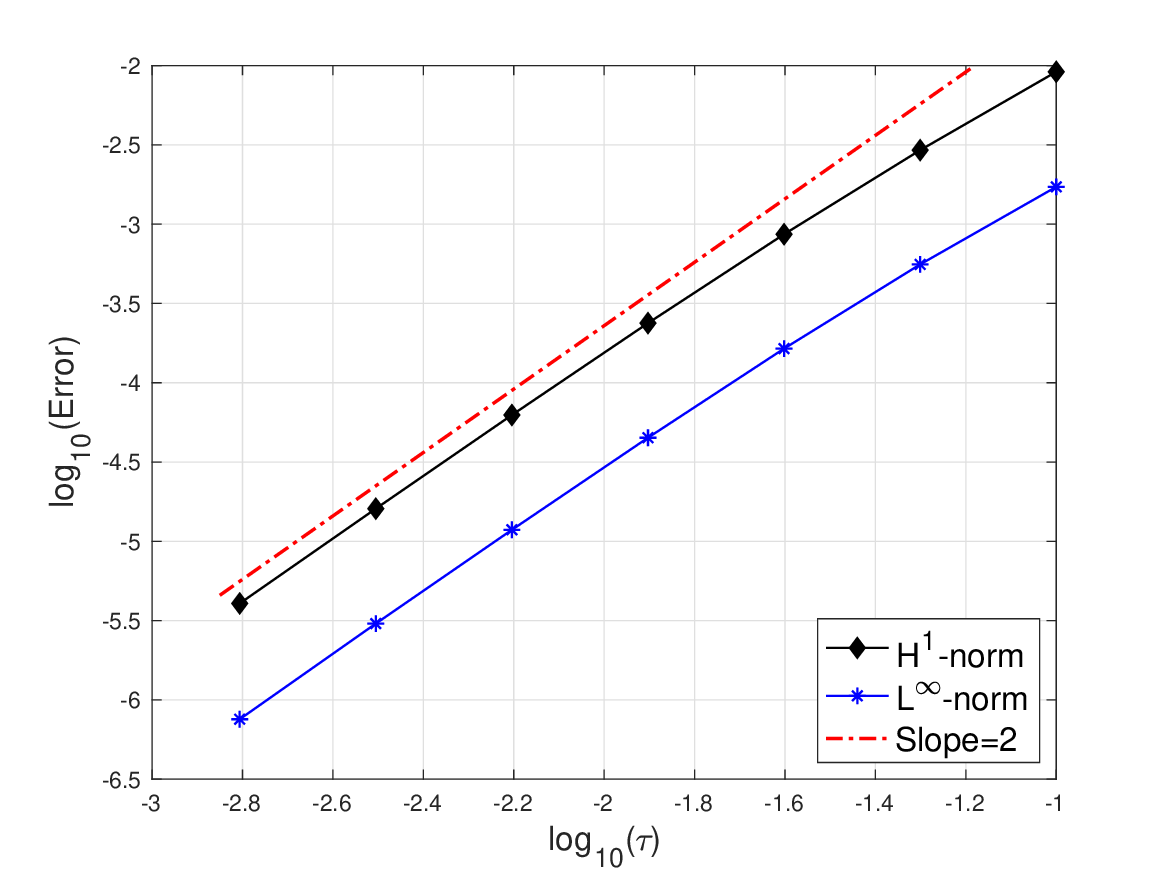}}
\caption{Plots of the numerical solution errors vs. the time step sizes in the log-log scale for the BDF2 scheme \eqref{BDF_2} with uniform time steps.
Left: $M(\phi)\equiv1$; right: $M(\phi)=1-\phi^{2}$.
}\label{fig1}
\end{figure*}

The central finite difference method is used for the spatial discretization with the fixed small mesh size $h=1/1024$.
Since there is no  analytical solution available for this example  to exactly evaluate  the numerical solution errors, we instead compute their approximations in the   discrete $L^{\infty}$ and $H^1$ norms, respectively:
\beq
e^T_{\infty}=\|\vPhi^{N}-\vPhi^{2N}\|_{\infty}, \quad e^T_{H^1}=\|\vPhi^{N}-\vPhi^{2N}\|_{H^1_h},
\eeq
where  $\vPhi^{N}$ and  $\vPhi^{2N}$ denote the numerical solution at the terminal time $T=1$ with $N$ and $2N$ subintervals for the time domain $[0,1]$, respectively.
To validate the theoretical temporal accuracy,
we firstly investigate the error behaviors of the BDF2 scheme \eqref{BDF_2} with the uniform time steps by repeatedly refining the time step size $\tau$
from $1/10$ to $1/640$ (i.e., $N$ changes from $10$ to $640$).
The solution errors vs. the time step sizes  are plotted in Fig. \ref{fig1} in the log-log scale for both mobility functions.
It is observed that the BDF2 scheme \eqref{BDF_2} achieves the expected second-order temporal accuracy for all test cases.
Next, we numerically study the error behaviors of the BDF2 scheme \eqref{BDF_2} with  nonuniform time steps.
The  nonuniform time step sizes $\{\widehat{t}_{n}\}_{n=0}^{N}$ used here is produced  by $25\%$ perturbation of the uniform ones $\{t_{n}=n/N\}_{n=0}^{N}$. As reported in Table \ref{table1}, the second-order temporal accuracy is still achieved by the BDF2 scheme for all cases.

\begin{table}[!t]\small
  \begin{center}
      \caption{Numerical solution errors and convergence rates of the BDF2 scheme \eqref{BDF_2} with nonuniform time steps.}\label{table1}
    \begin{tabular}{|ccc|cc|cc|cc|cc|} \hline
    \multicolumn{3}{|c|}{Time steps}&\multicolumn{4}{|c|}{$M(\phi)\equiv 1$ }&\multicolumn{4}{|c|}{$M(\phi)=1-\phi^{2}$ } \\ \hline
     $N$&$\tau$&$\max\{\gamma_{n}\}$&$e^T_{\infty}$&Order&$e^T_{H^1_h}$&Order&$e^T_{\infty}$&Order&$e^T_{H^1_h}$&Order\\ \hline
     $10$ &1.393e-1&2.282&1.669e-2&--   &1.375e-1&--    & 1.886e-3& --  &1.005e-2& --   \\
     $20$ &7.033e-2&2.358&5.834e-3&1.54&4.941e-2&1.50 & 6.399e-4&1.58 &3.371e-3&1.60 \\
     $40$ &3.408e-2&2.218&1.597e-3&1.79&1.376e-2&1.77 & 1.782e-4&1.76 &9.391e-4&1.76 \\
    $80$  &1.785e-2&2.656&3.728e-4&2.25&3.262e-3&2.23 & 4.255e-5&2.21 &2.239e-4&2.22 \\
    $160$ &9.061e-3&2.712&9.872e-5&1.96&8.648e-4&1.96 & 1.136e-5&1.95 &6.002e-5&1.94 \\
    $320$ &4.638e-3&2.832&2.511e-5&2.05&2.201e-4&2.04 & 2.904e-6&2.04 &1.537e-5&2.03  \\
   $640$ &2.289e-3 &2.717&6.293e-6&1.96&5.527e-5&1.96 & 7.160e-7&1.98 &3.849e-6&1.96 \\
      \hline
    \end{tabular}
  \end{center}
\end{table}

\subsection{Test of MBP preservation}

We demonstrate the MBP preservation of the proposed BDF2 scheme \eqref{BDF_2} through two well-known benchmark examples governed by the
Allen-Cahn equations. One is the shrinking bubble problem \cite{Chen98} and the other is the grain coarsening problem.

\paragraph*{{\bf The shrinking bubble problem}}
We consider the Allen-Cahn equation \eqref{prob} with $M(\phi)\equiv1$ and $\varepsilon=0.01$ in a rectangular domain $(-0.5,0.5)^2$. The initial bubble is given by
\beq
\phi_{0}(\x)=\begin{cases}
\begin{array}{r@{}l}
1,&\quad |\x|^2<0.2^{2},\\[1pt]
-1,&\quad |\x|^2\geq0.2^{2}.
\end{array}
\end{cases}
\eeq
As discussed in \cite{Chen98,JZZD15,HAX19,CGJMP19}, this model describes the evolution in time of a shrinking bubble  with the initial radius $R_{0}=0.2$, and the velocity of this circular moving interface approximately satisfies the following  relation
\bq\label{test_1}
R(t)=\sqrt{R^{2}_{0}-2\varepsilon^{2}t},
\eq
if $\varepsilon$ is sufficiently small.
Here, $R(t)$ is the radius of the circle at time t.

\begin{figure}[!t]
\centerline{\includegraphics[scale=0.29]{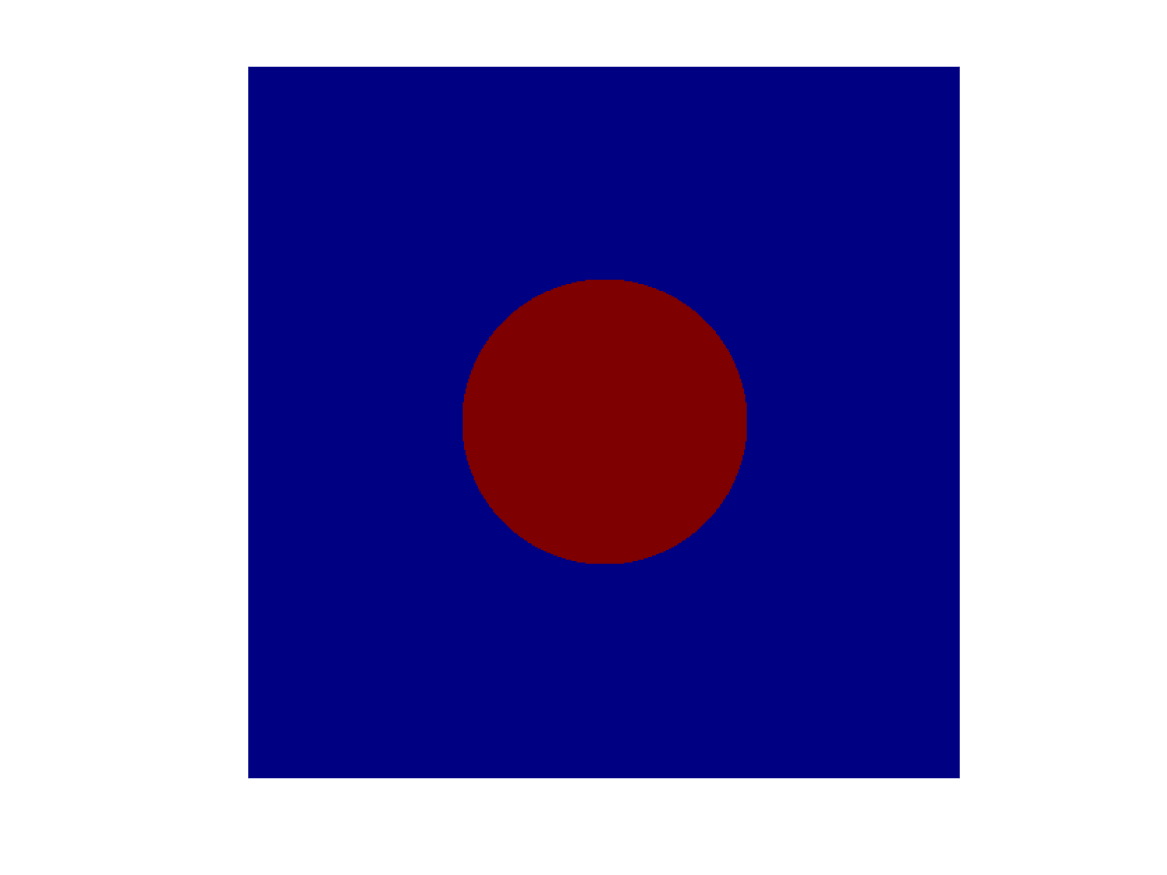}\hspace{-0.95cm}\includegraphics[scale=0.29]{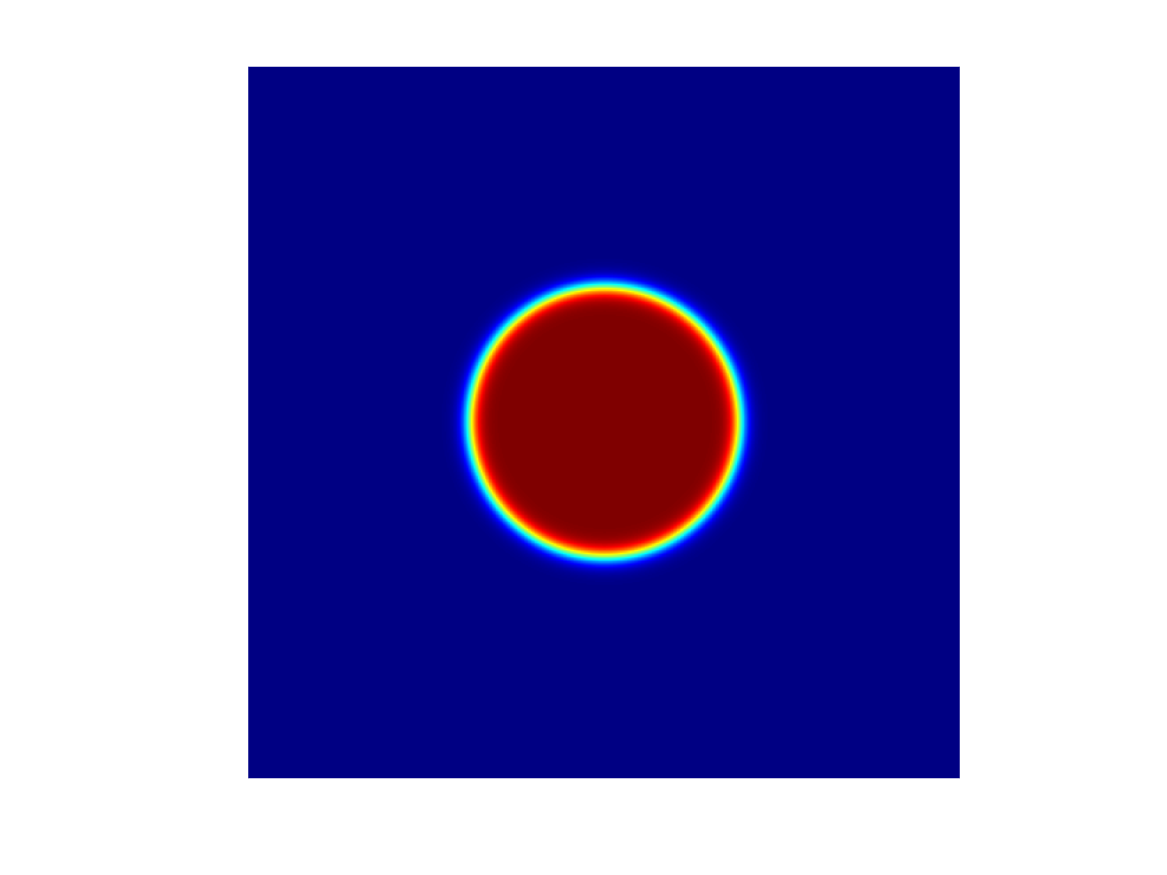}\hspace{-0.95cm}
\includegraphics[scale=0.29]{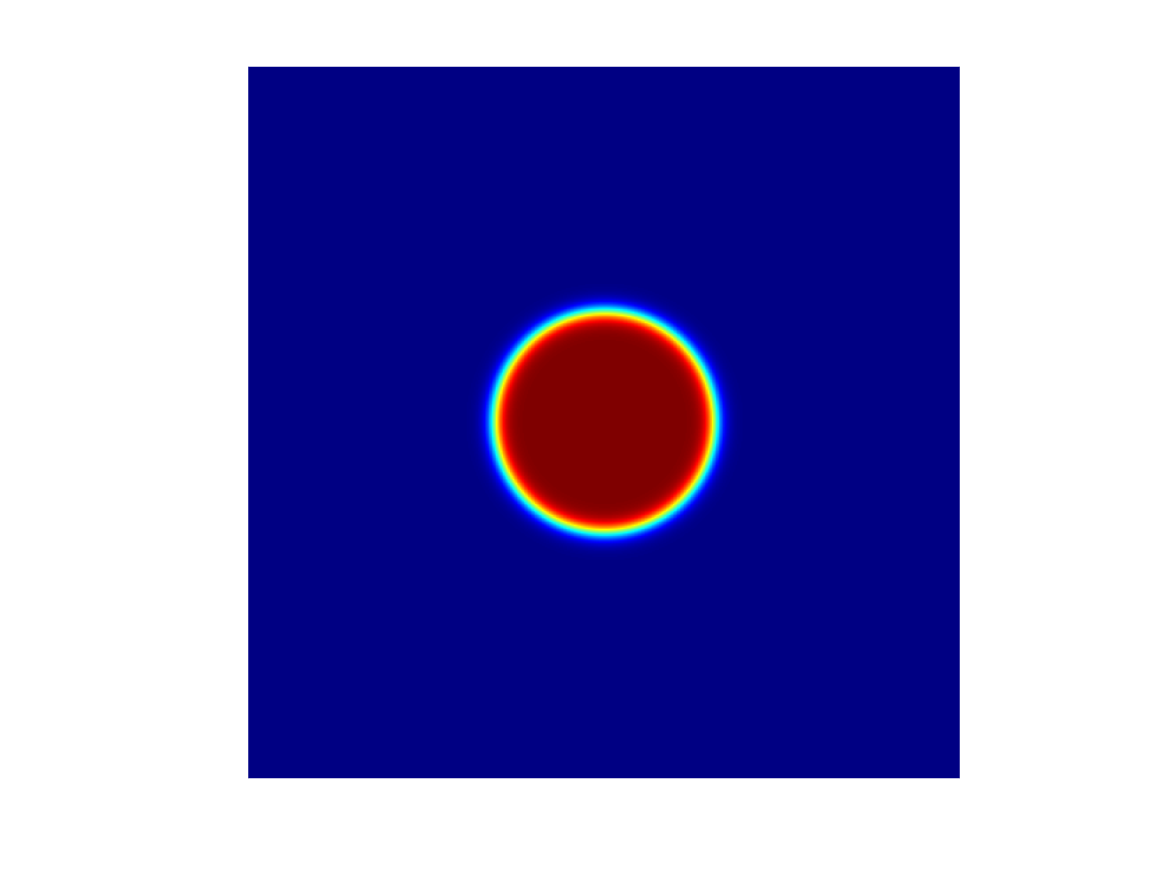}}
\centerline{\includegraphics[scale=0.29]{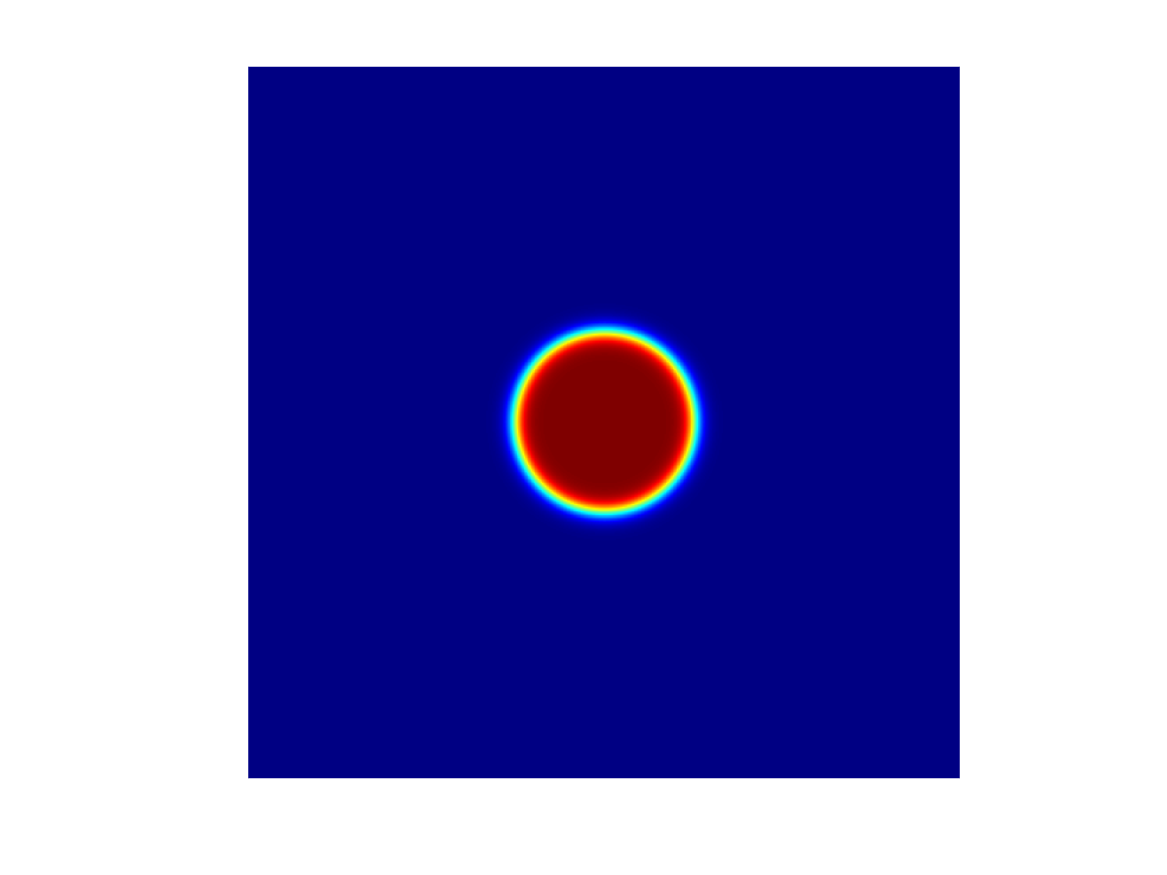}\hspace{-0.95cm}\includegraphics[scale=0.29]{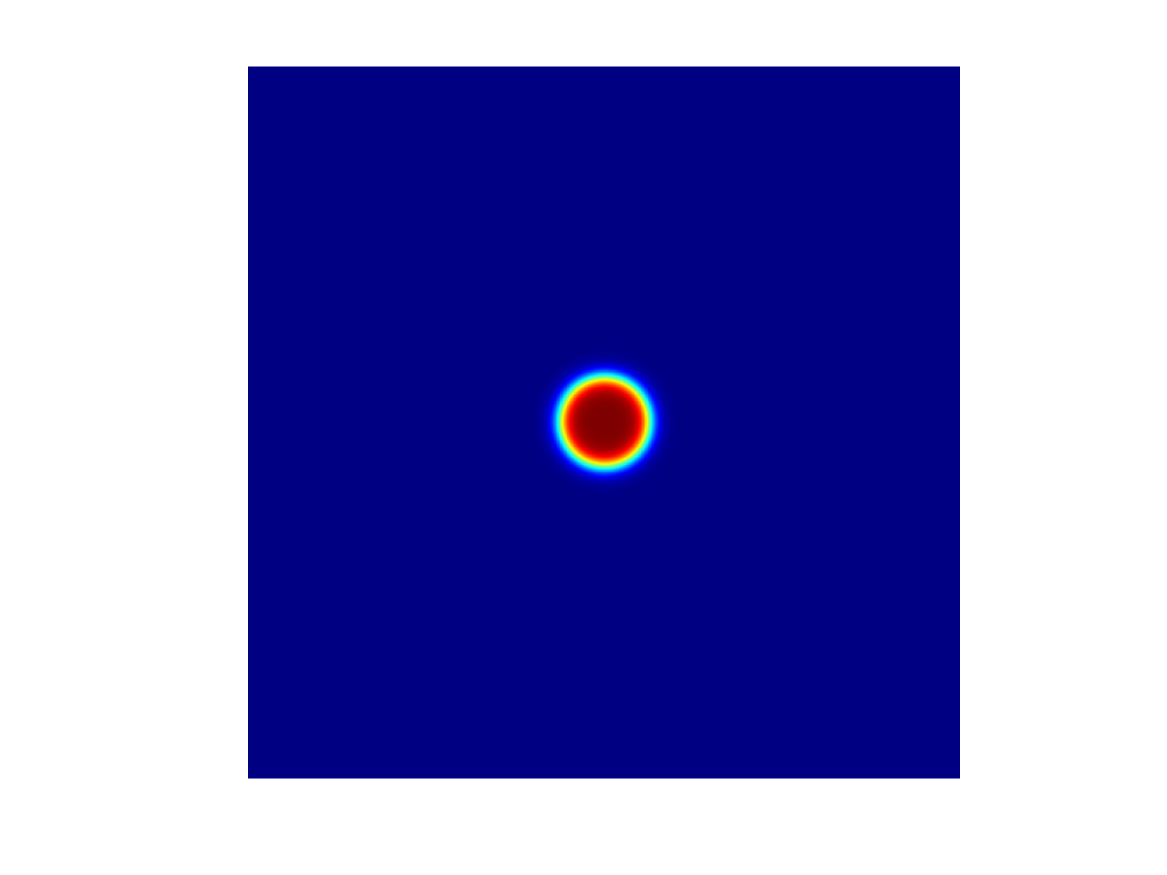}\hspace{-0.95cm}
\includegraphics[scale=0.29]{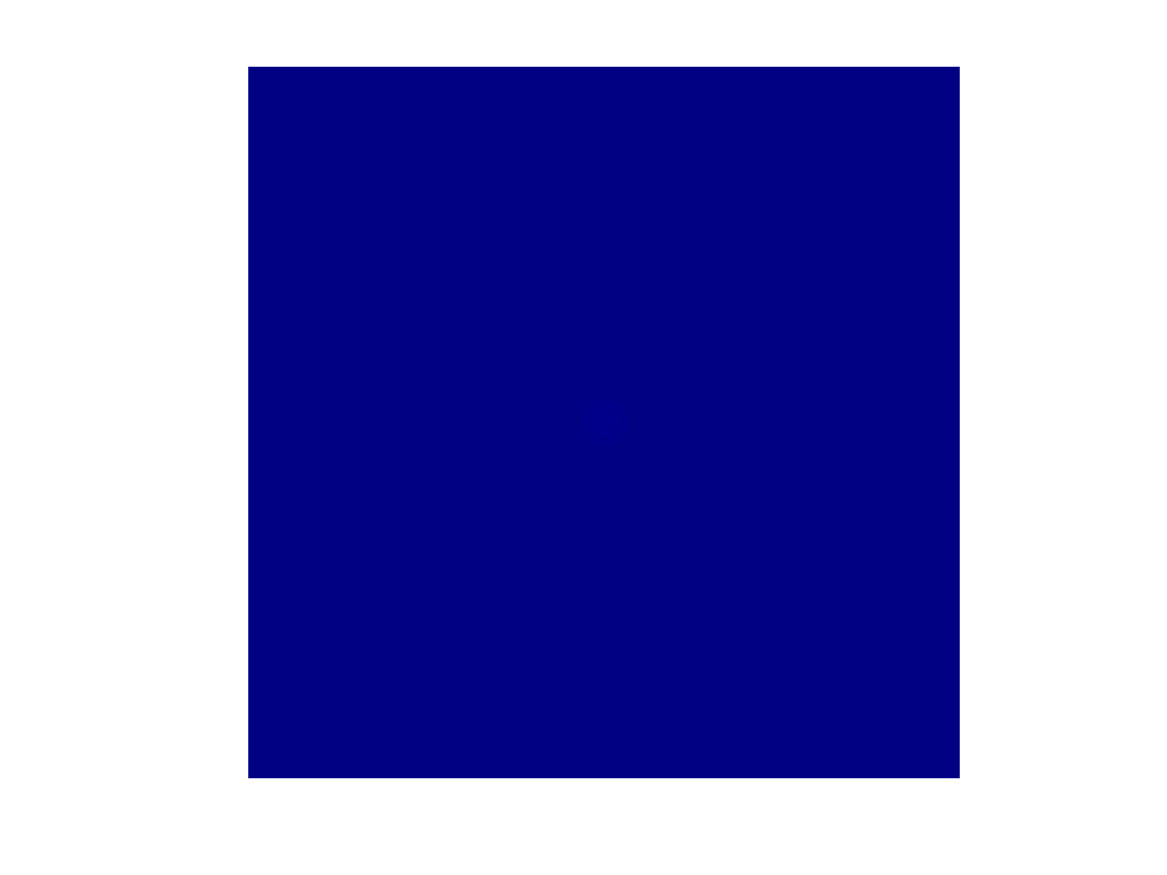}}
\caption{Snapshots of the simulated phase structures at the times $t = 0$, $20$, $80$, $120$, $180$,  and $200$ produced by the BDF2 scheme \eqref{BDF_2} for the shrinking bubble problem.}
\label{fig2}
\end{figure}

\begin{figure*}[htbp]
\begin{minipage}[t]{0.45\linewidth}
\centerline{\includegraphics[scale=0.35]{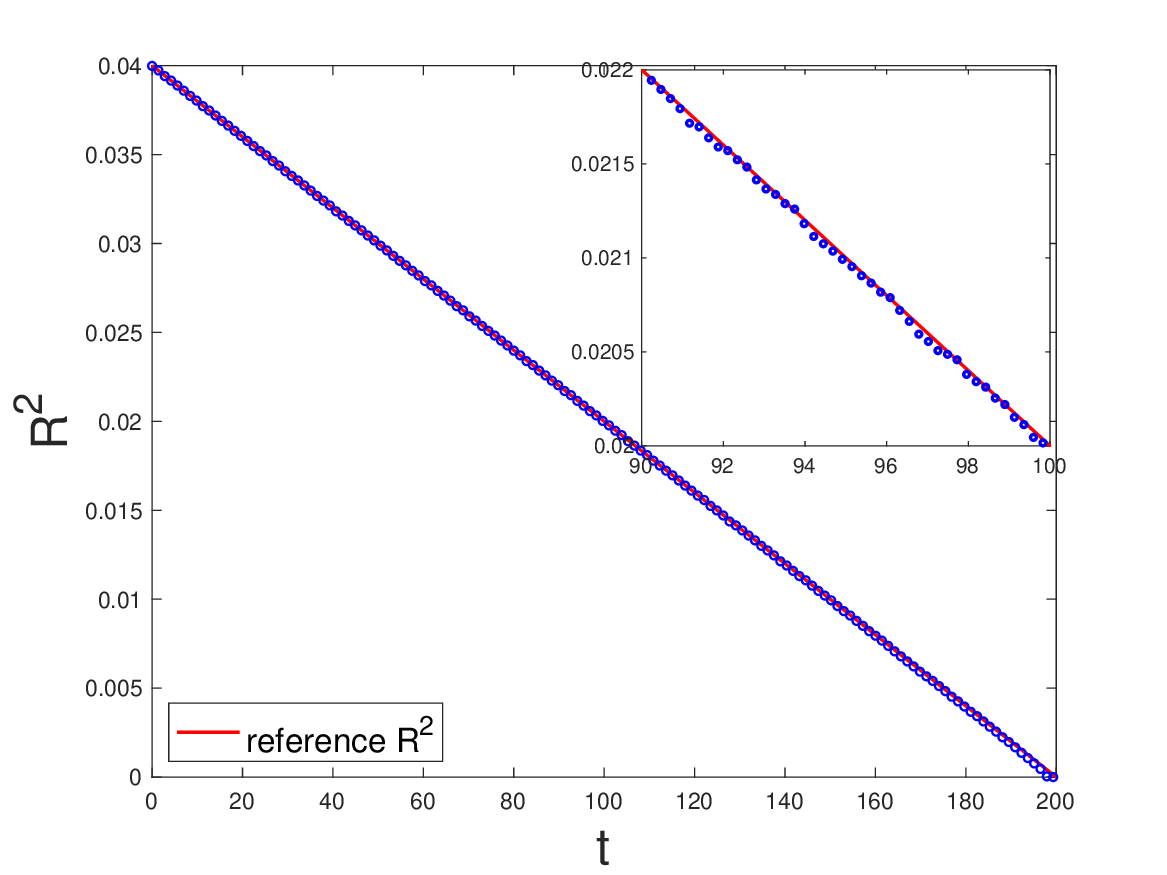}}
\centerline{(a) the radius  }
\end{minipage}
\begin{minipage}[t]{0.45\linewidth}
\centerline{\includegraphics[scale=0.35]{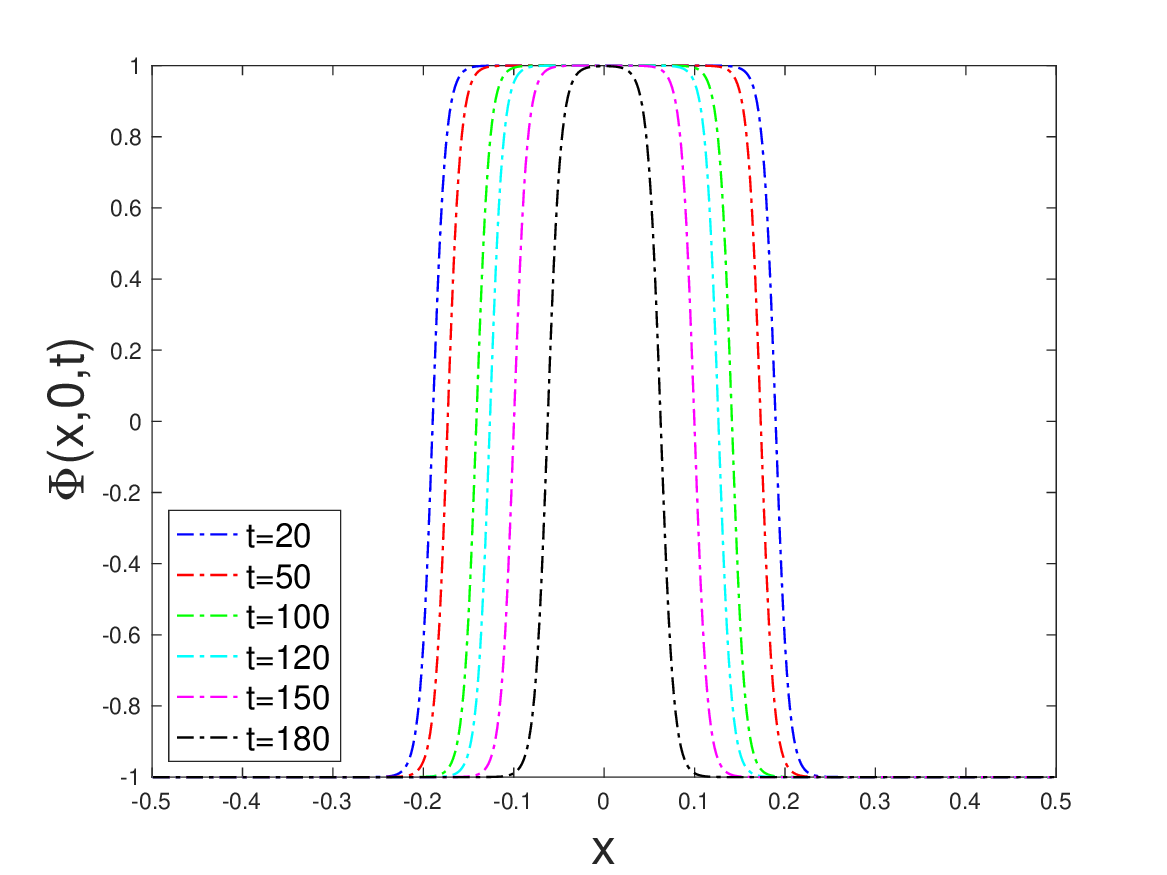}}
\centerline{(b) the cross-section view with $y=0$}
\end{minipage}
\vskip 3mm
\begin{minipage}[t]{0.45\linewidth}
\centerline{\includegraphics[scale=0.35]{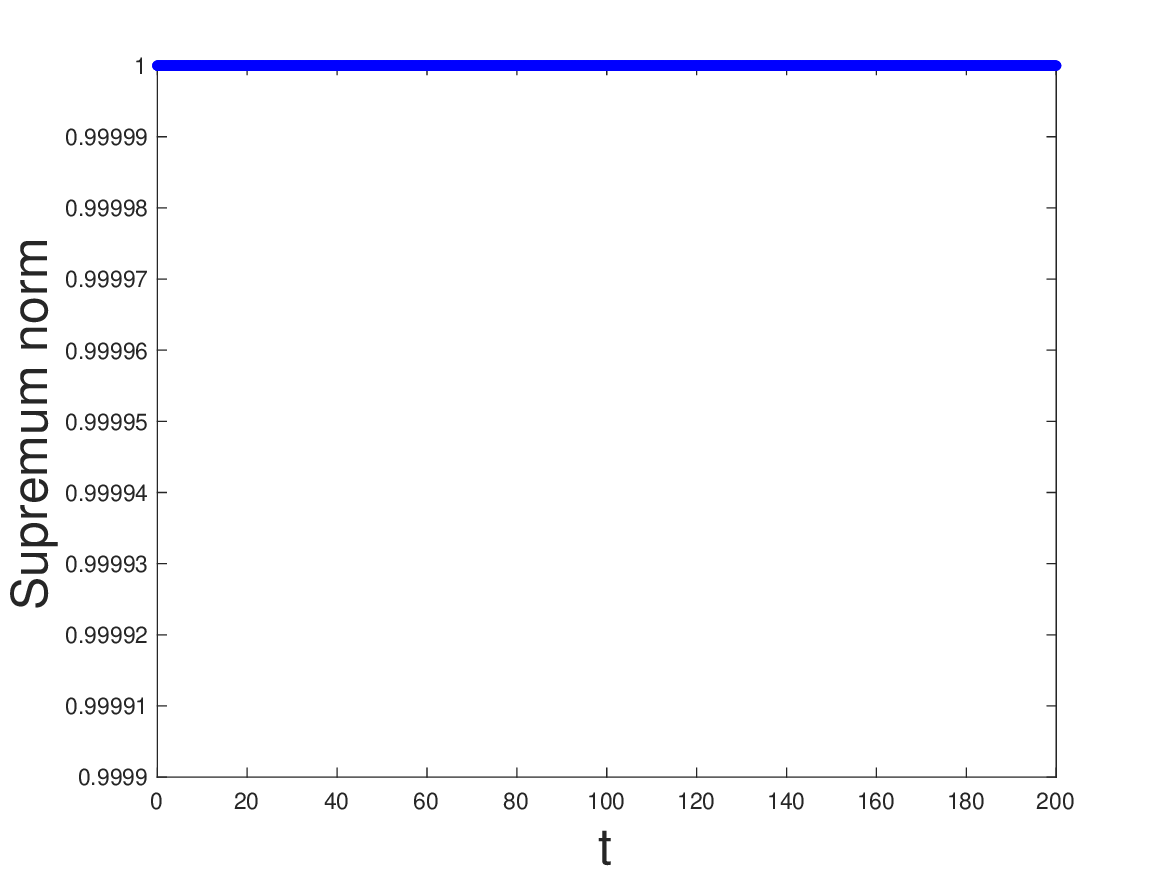}}
\centerline{(c) the supremum norm}
\end{minipage}
\begin{minipage}[t]{0.45\linewidth}
\centerline{\includegraphics[scale=0.35]{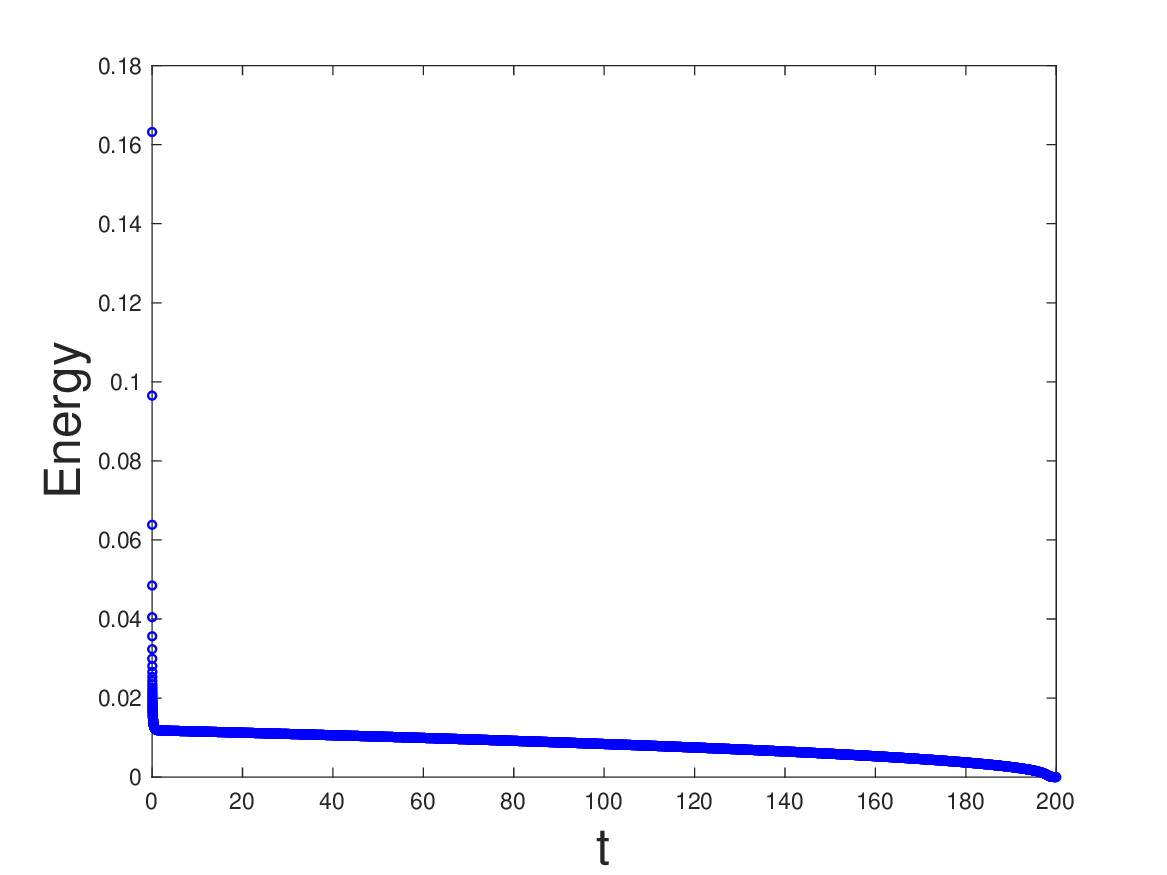}}
\centerline{(d) the energy}
\end{minipage}
\caption{The evolutions in time of the radius, the cross section of $y=0$, the supremum norm, and the energy of the simulated solution produced by the BDF2 scheme \eqref{BDF_2} with uniform time steps for the shrinking bubble problem.}\label{fig3}
\end{figure*}

The simulation is performed by the BDF2 scheme \eqref{BDF_2} with $h=1/512.$ The uniform time steps are used here with the time step size $\tau=\mathcal{G}(1)/[S+4\varepsilon^{2}/h^{2}]$, which is the maximum value satisfying the requirement \eqref{est_tau3}.
Snapshots of the simulated bubble at the times $t =0,20,80,120,180, 200$ are displayed in Fig. \ref{fig2}, which shows that the bubble  disappears at $t=200$ as expected.
Moreover, we plot the evolution in time of the radius of the simulated bubble in Fig. \ref{fig3}-(a), which matches the prediction \eqref{test_1} very well. Several cross-section views with $y=0$ for the simulated solution are presented  in Fig. \ref{fig3}-(b) and the evolution of its supremum norm  along with the time is displayed in Fig. \ref{fig3}-(c), which  demonstrate the MBP preservation of the proposed BDF2 scheme \eqref{BDF_2} during the whole simulation. Furthermore, it is also observed that the energy of the simulated solution is monotonically decreasing in time as shown in Fig. \ref{fig3}-(d).

\paragraph*{{\bf The grain coarsening dynamics with  a time adaptive strategy}}

Finally, we investigate the efficiency and the MBP preservation of the proposed BDF2 scheme  \eqref{BDF_2} with a time adaptive strategy for the simulation of the grain
coarsening. The coarsening dynamic process usually goes through several different stages within a long period:  changes quickly at the beginning and then rather slowly until it reaches a steady state.
In  particular, we consider the coarsening dynamics governed by the Allen-Cahn equation \eqref{prob} with the nonlinear degenerate mobility $M(\phi)=1-\phi^{2}$ and $\varepsilon=0.01$.
Particularly, it is of great importance to preserve the numerical solution $\phi\in[-1,1]$ in the numerical algorithm for such a nonlinear mobility function.
Otherwise, the numerical solutions may blow up during the time simulation.

\begin{figure*}[htbp]
\centerline{\includegraphics[scale=0.26]{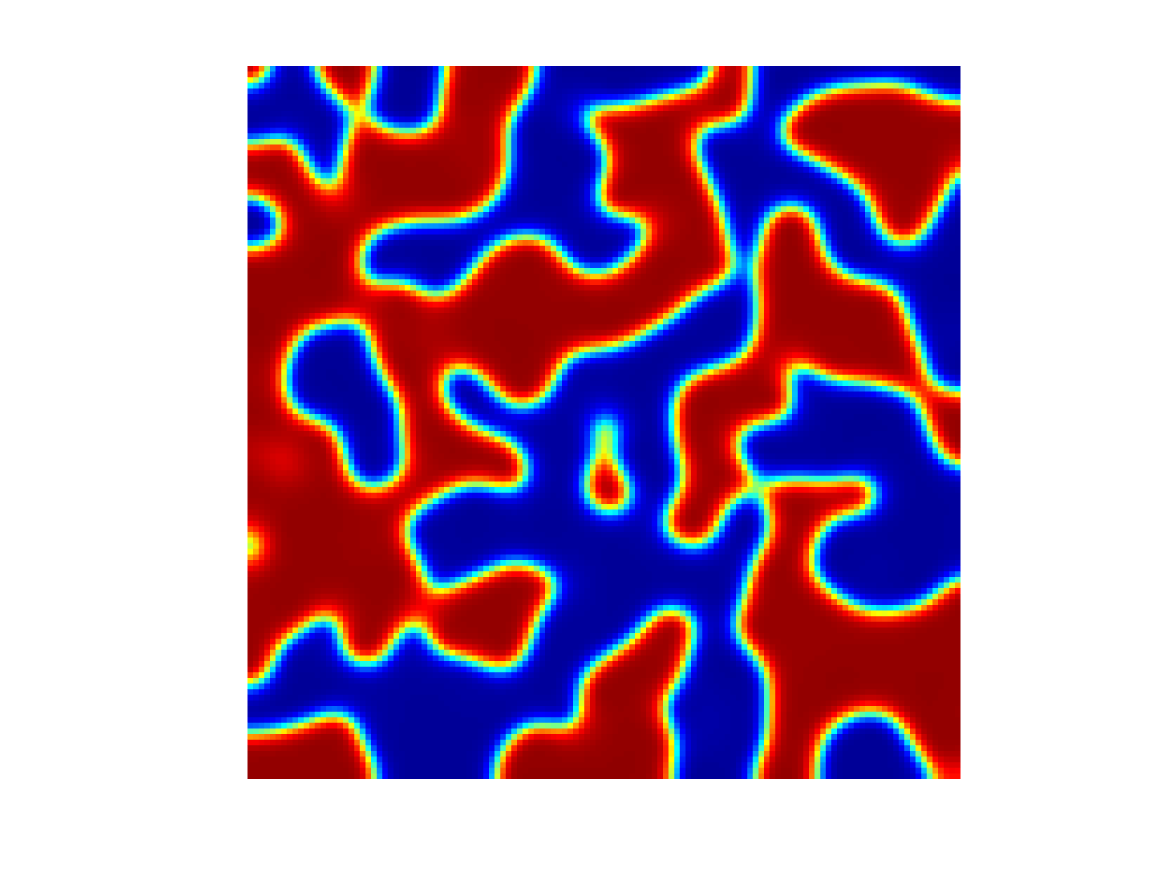}\hspace{-0.95cm}\includegraphics[scale=0.26]{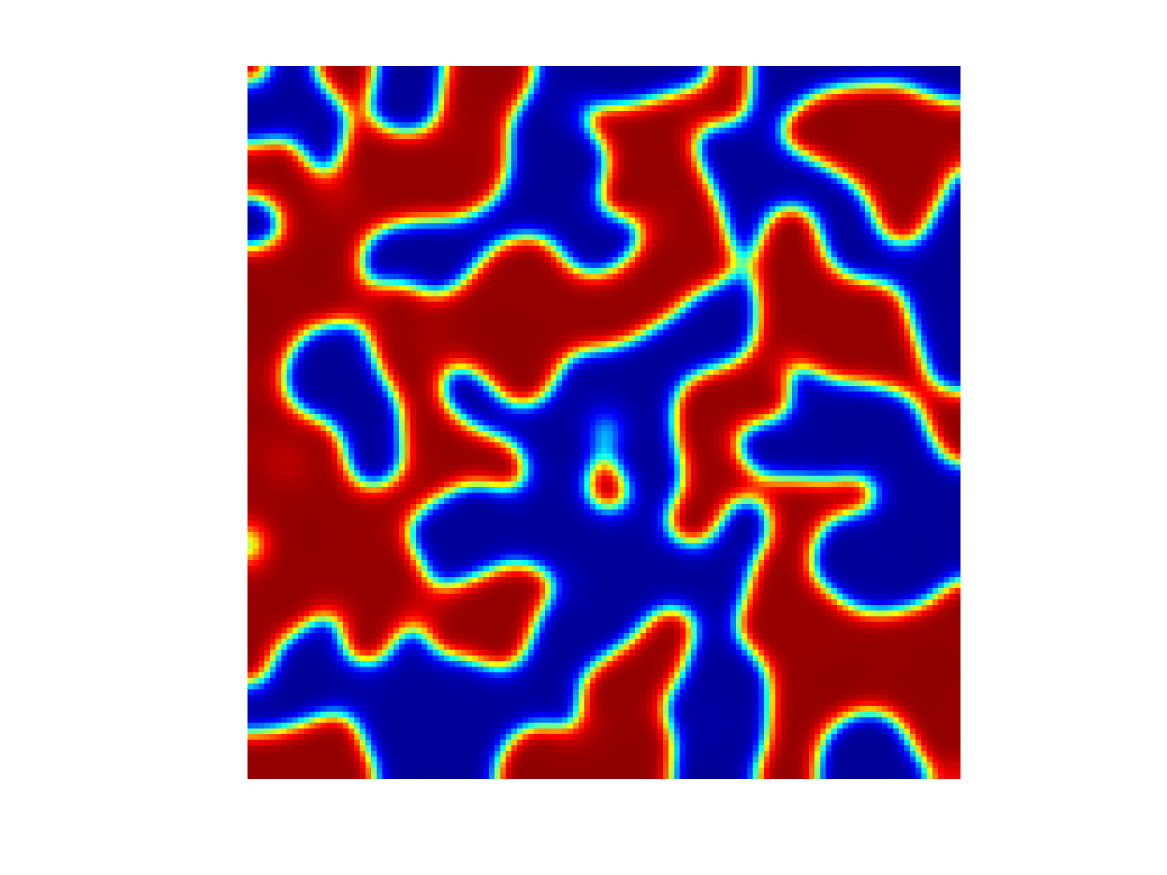}\hspace{-0.95cm}\includegraphics[scale=0.26]{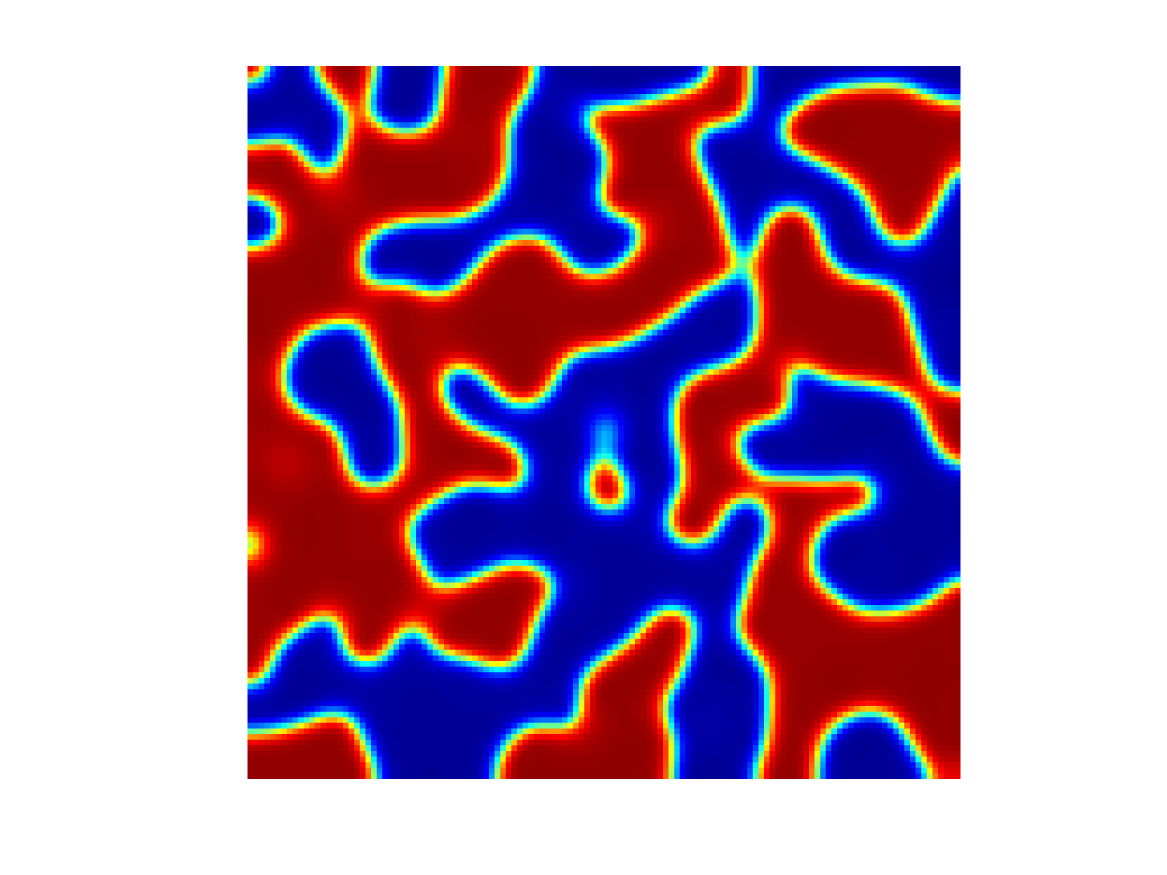}\hspace{-0.95cm}\includegraphics[scale=0.26]{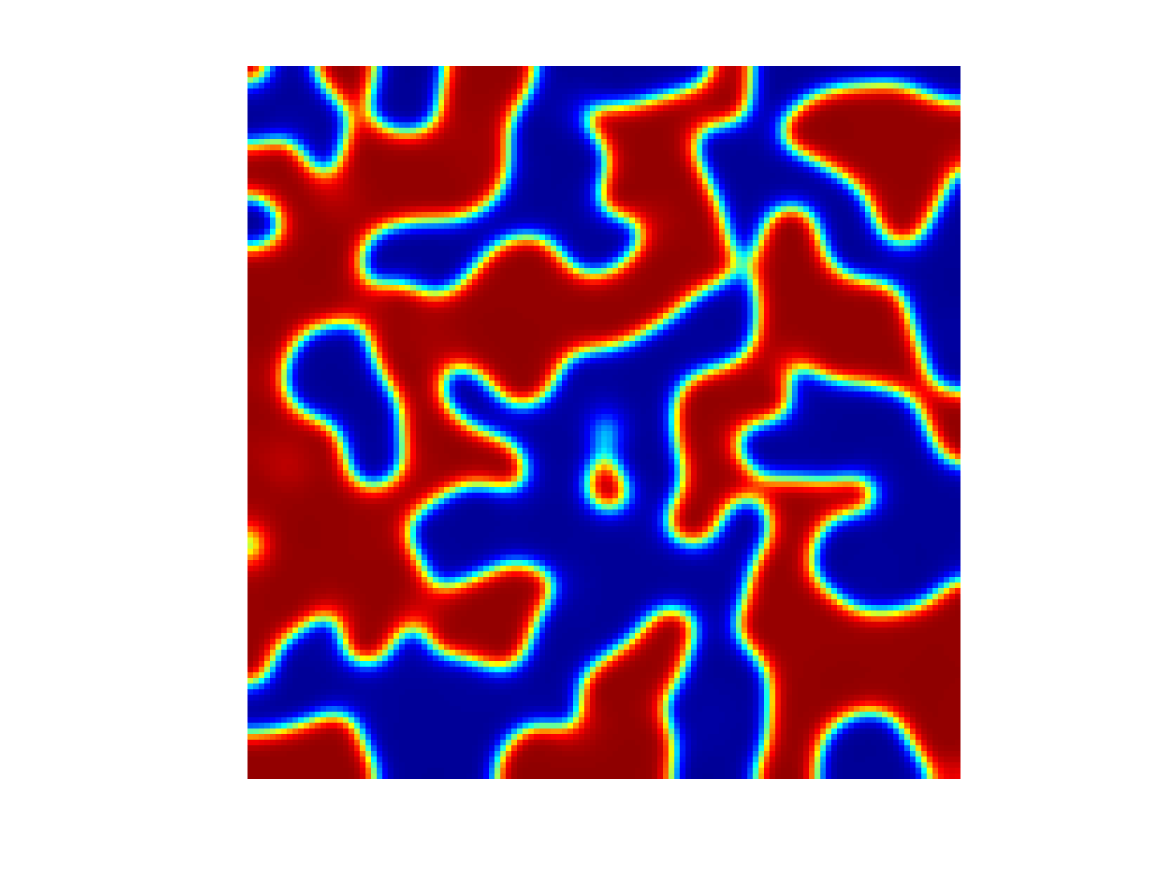}}
\vskip -3mm
\centerline{\includegraphics[scale=0.26]{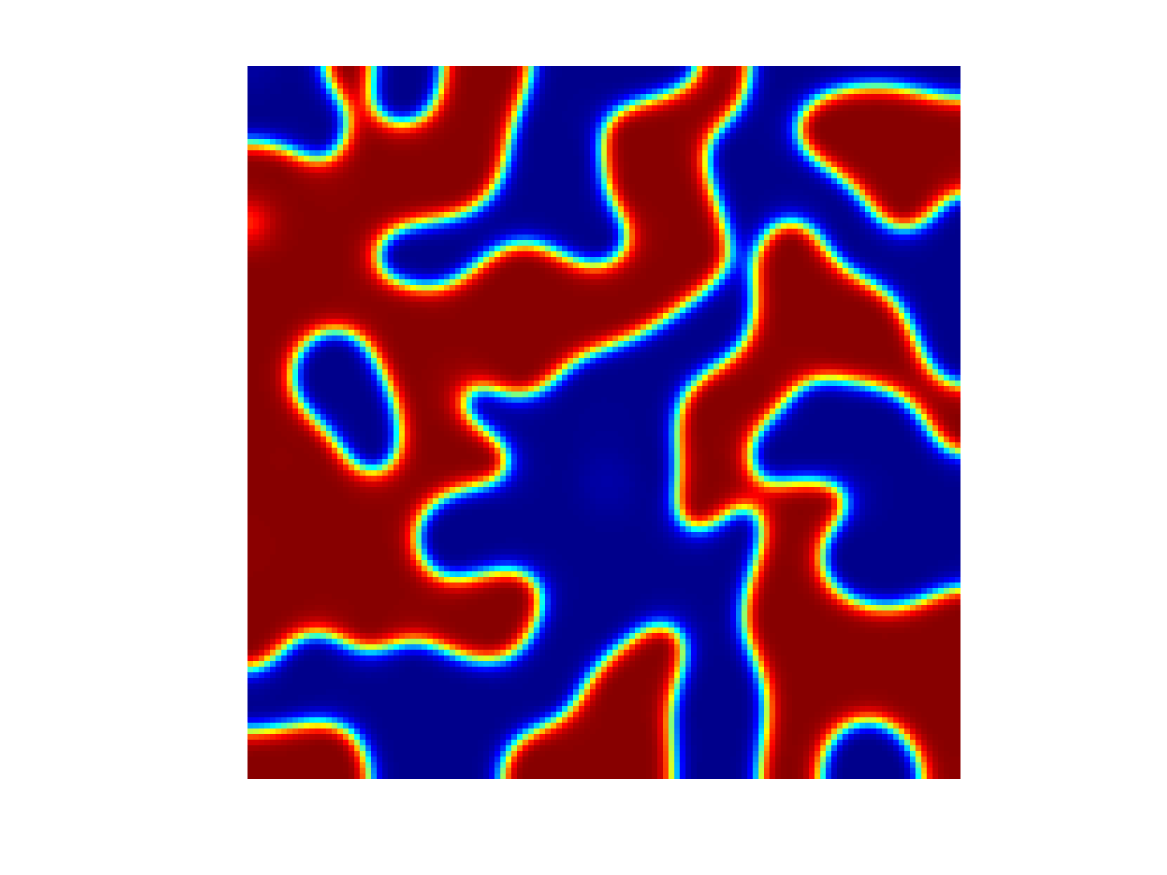}\hspace{-0.95cm}\includegraphics[scale=0.26]{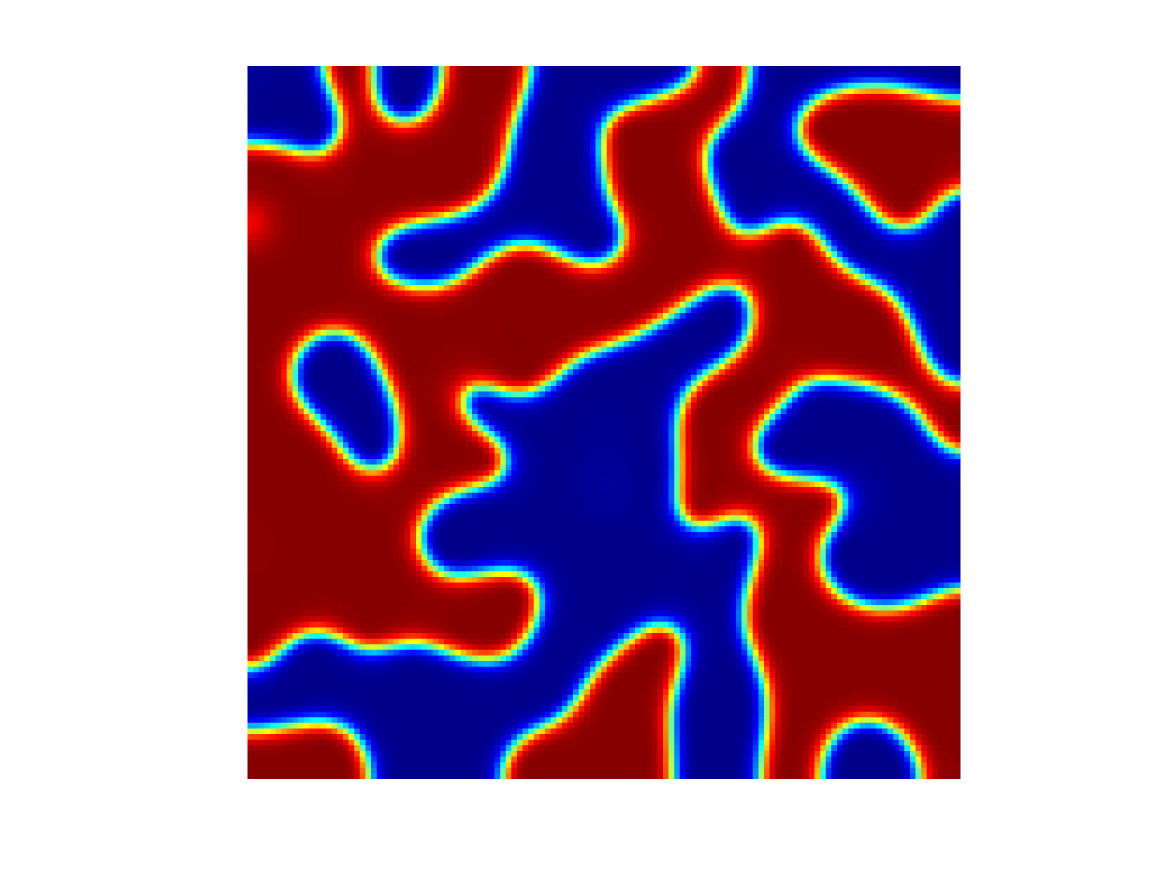}\hspace{-0.95cm}\includegraphics[scale=0.26]{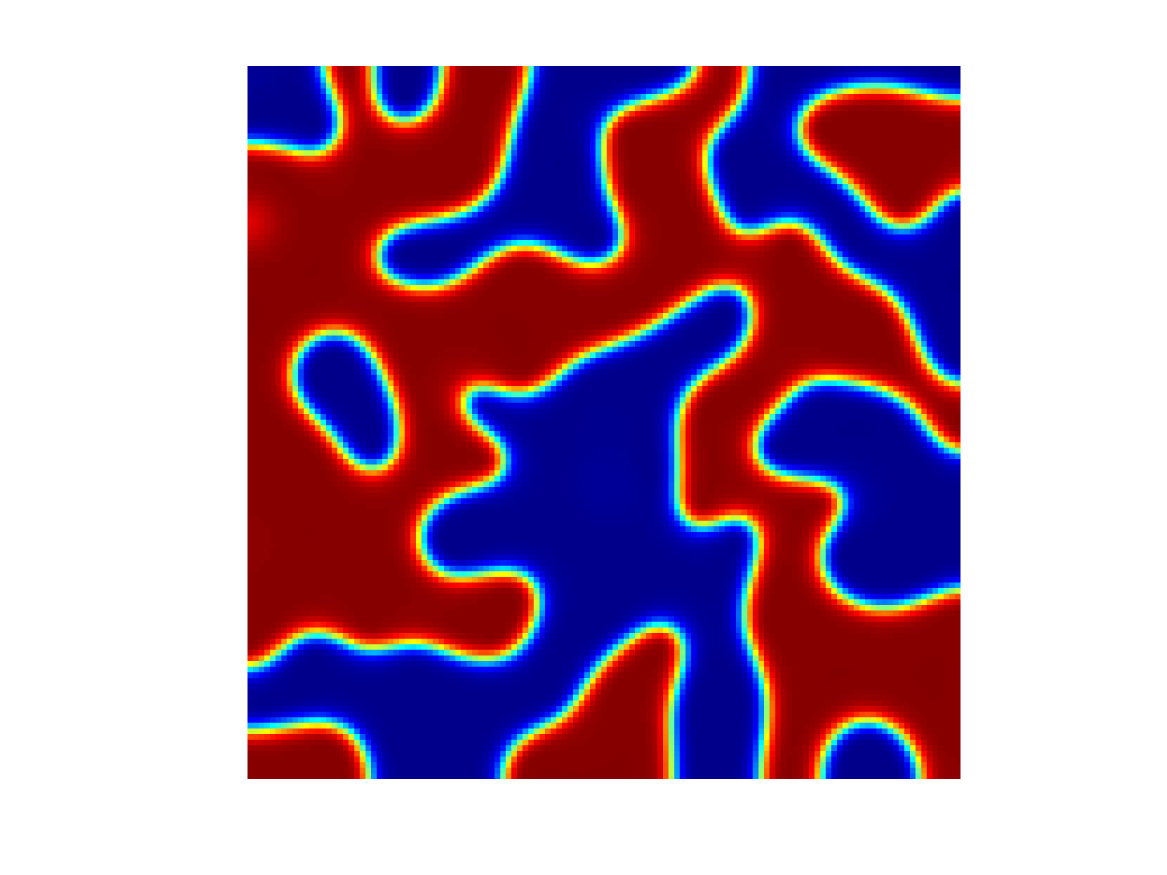}\hspace{-0.95cm}\includegraphics[scale=0.26]{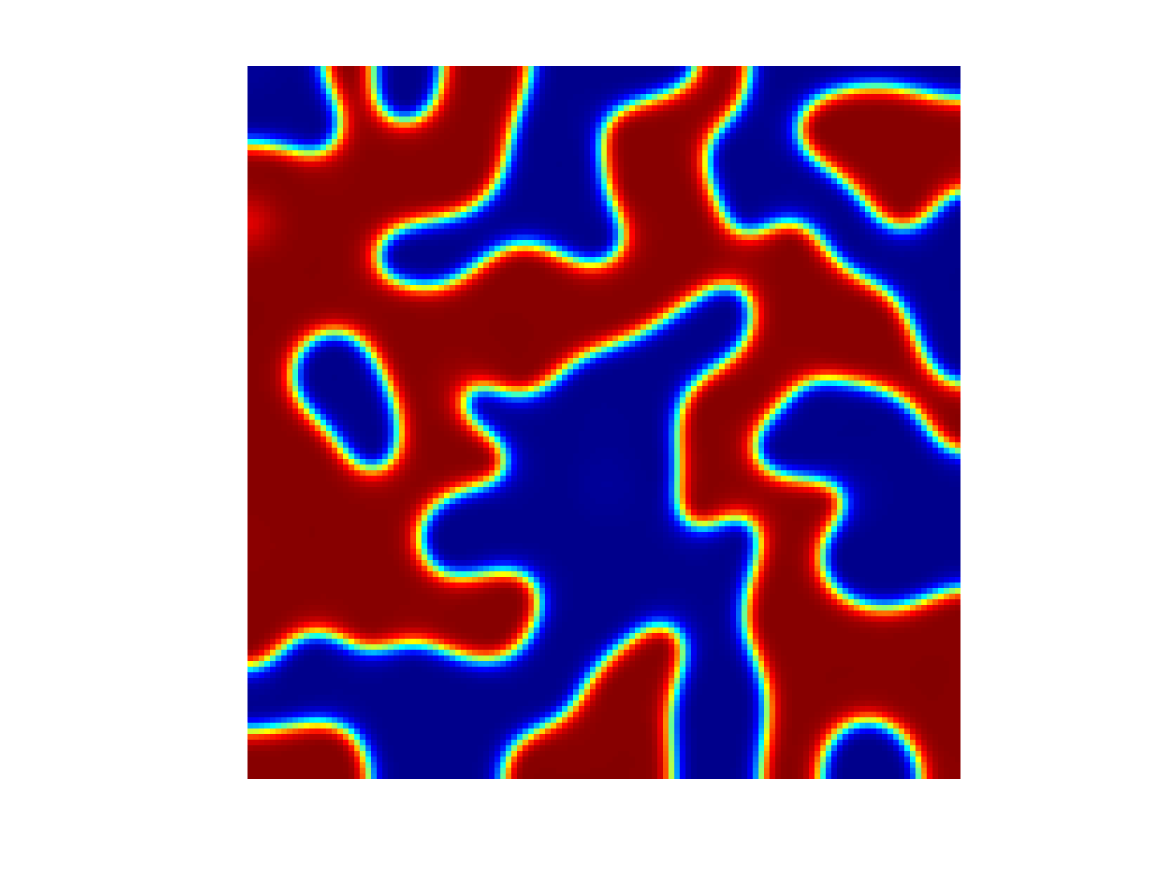}}
\vskip -3mm
\centerline{\includegraphics[scale=0.26]{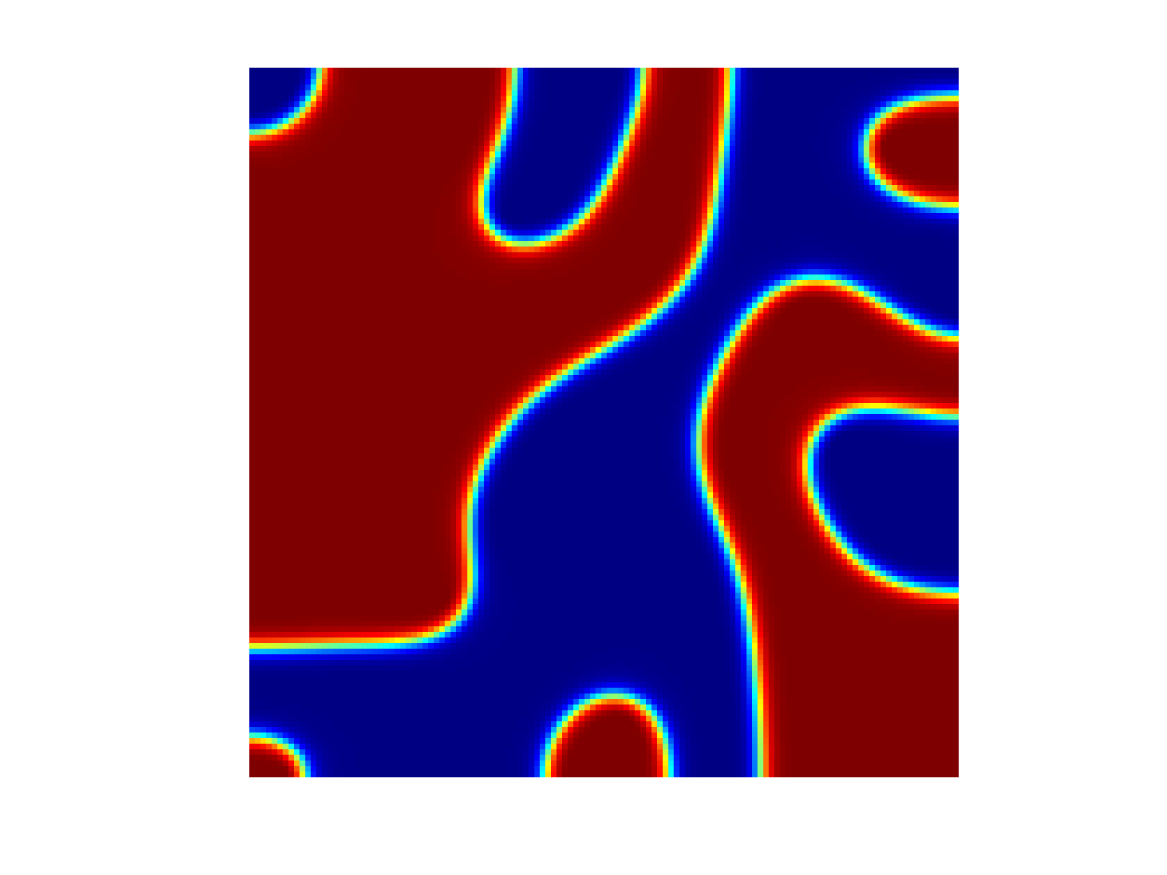}\hspace{-0.95cm}\includegraphics[scale=0.26]{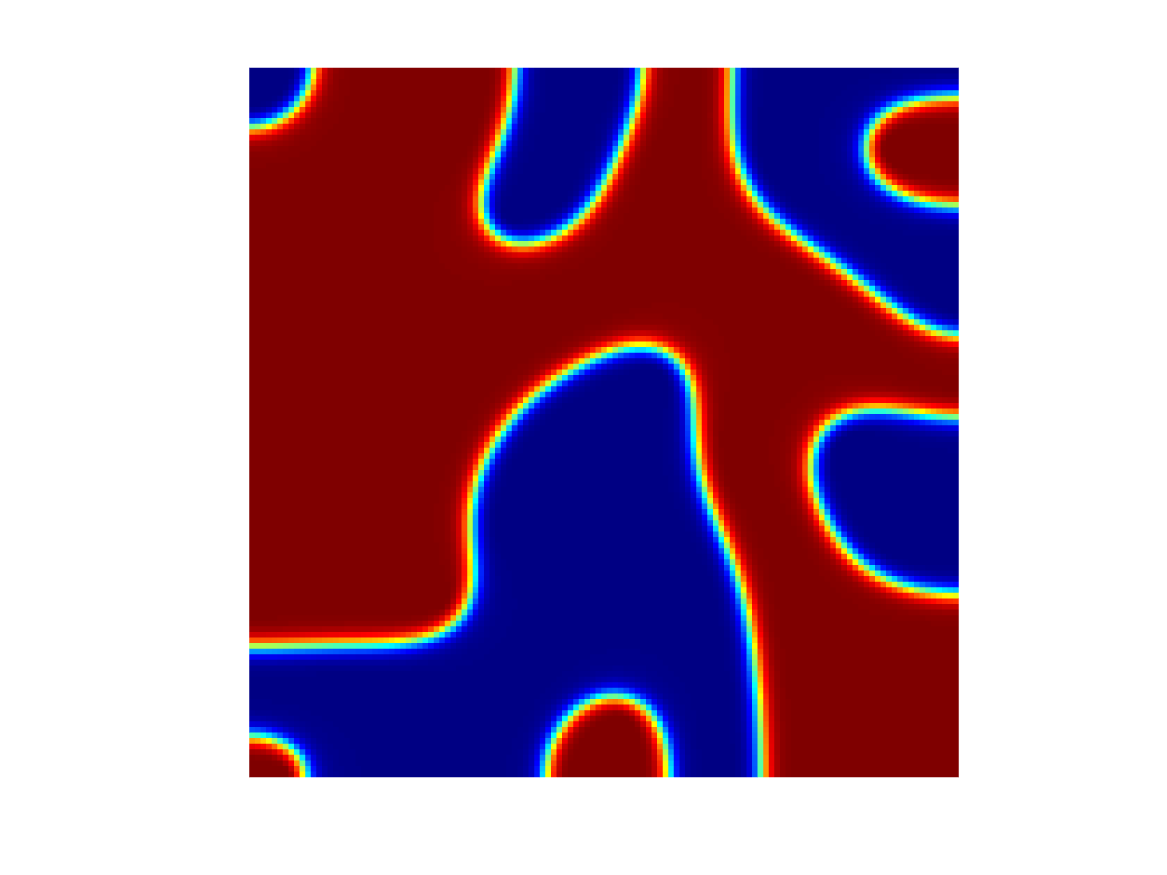}\hspace{-0.95cm}\includegraphics[scale=0.26]{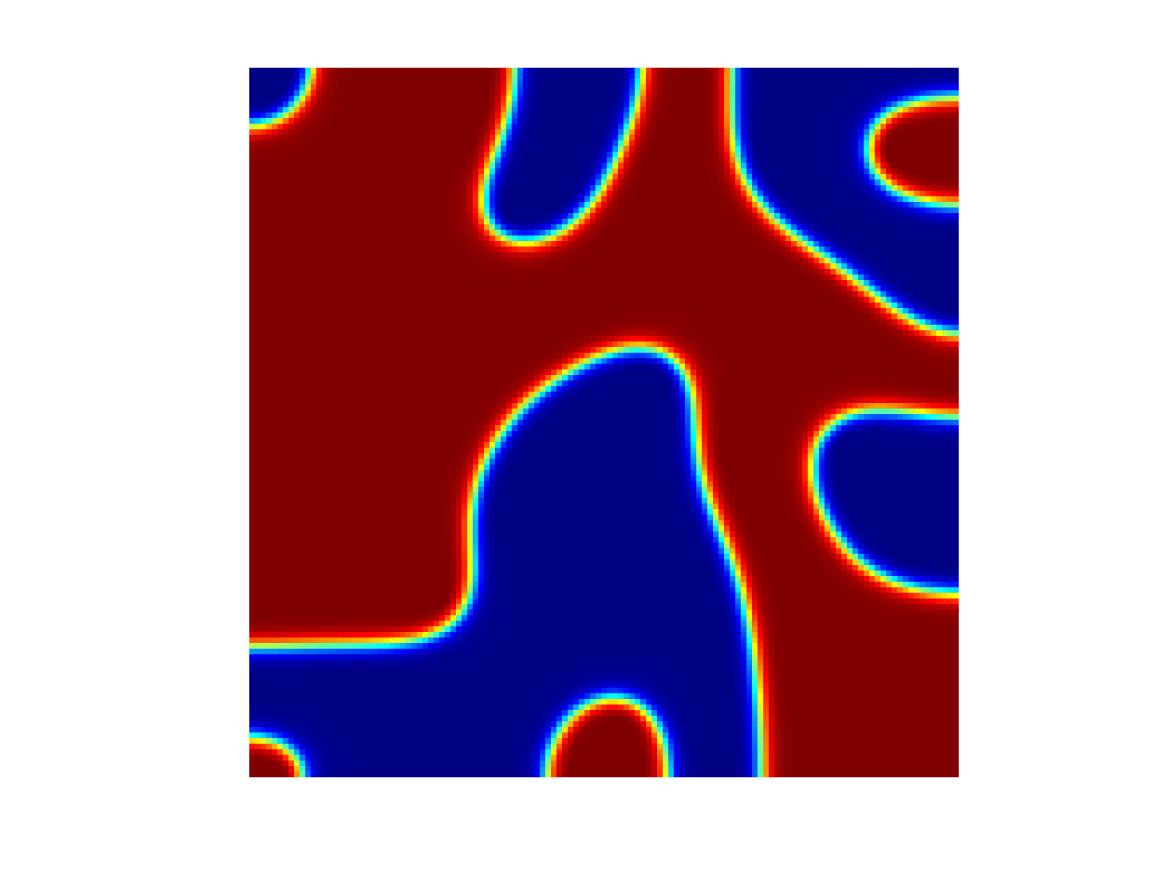}\hspace{-0.95cm}\includegraphics[scale=0.26]{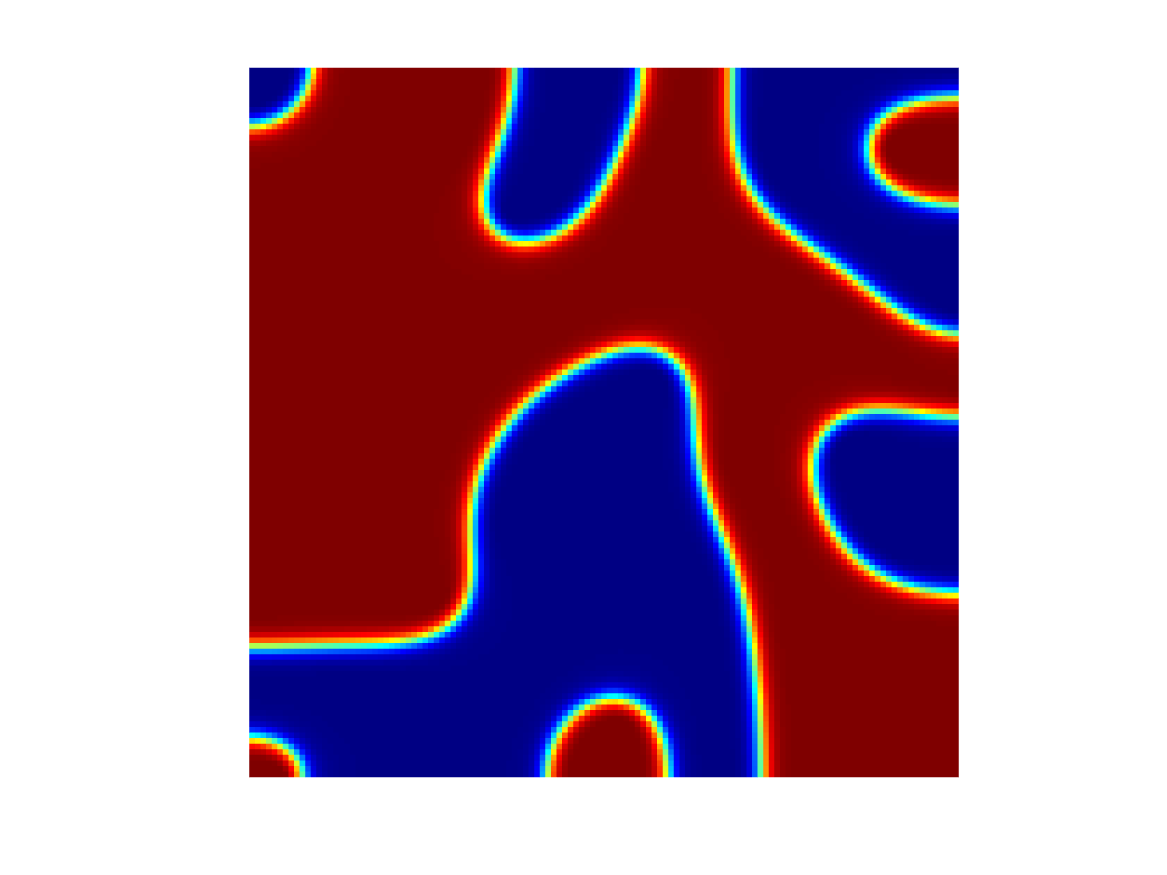}}
\vskip -3mm
\centerline{\includegraphics[scale=0.26]{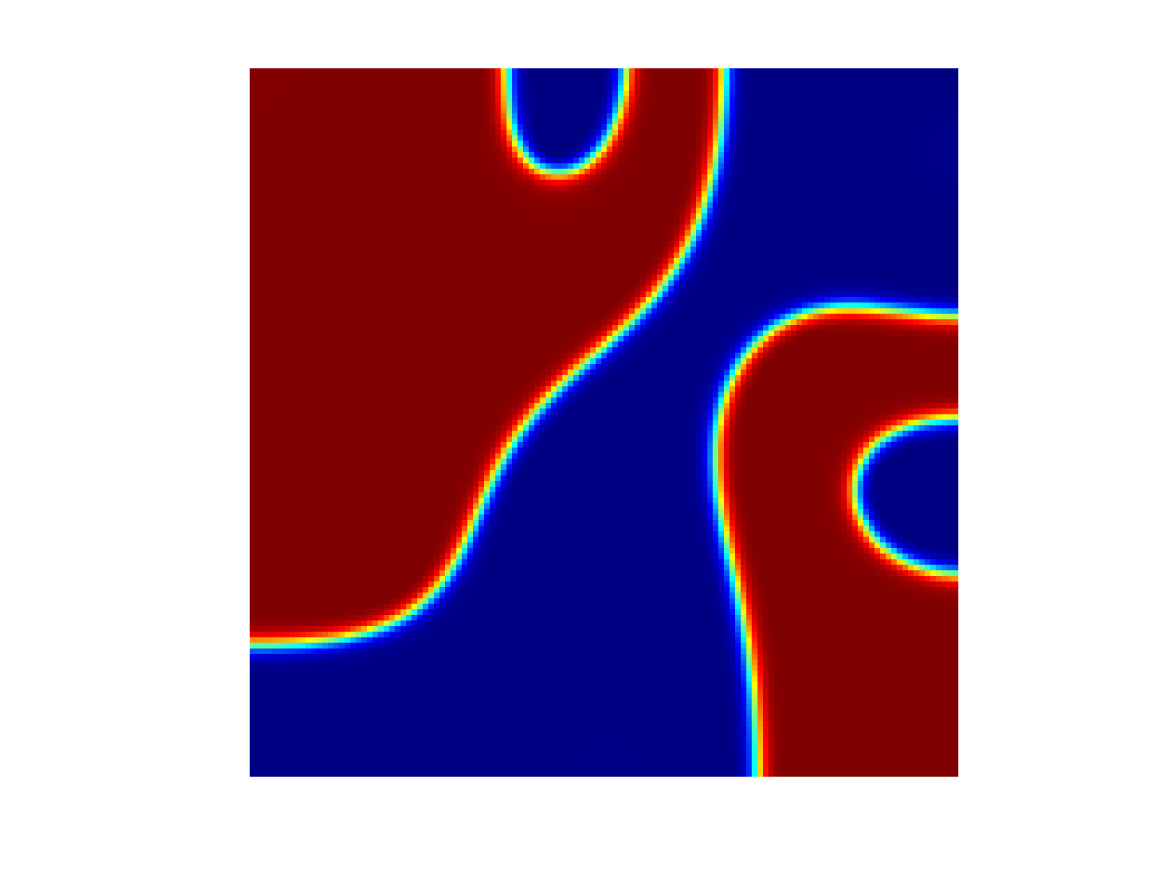}\hspace{-0.95cm}\includegraphics[scale=0.26]{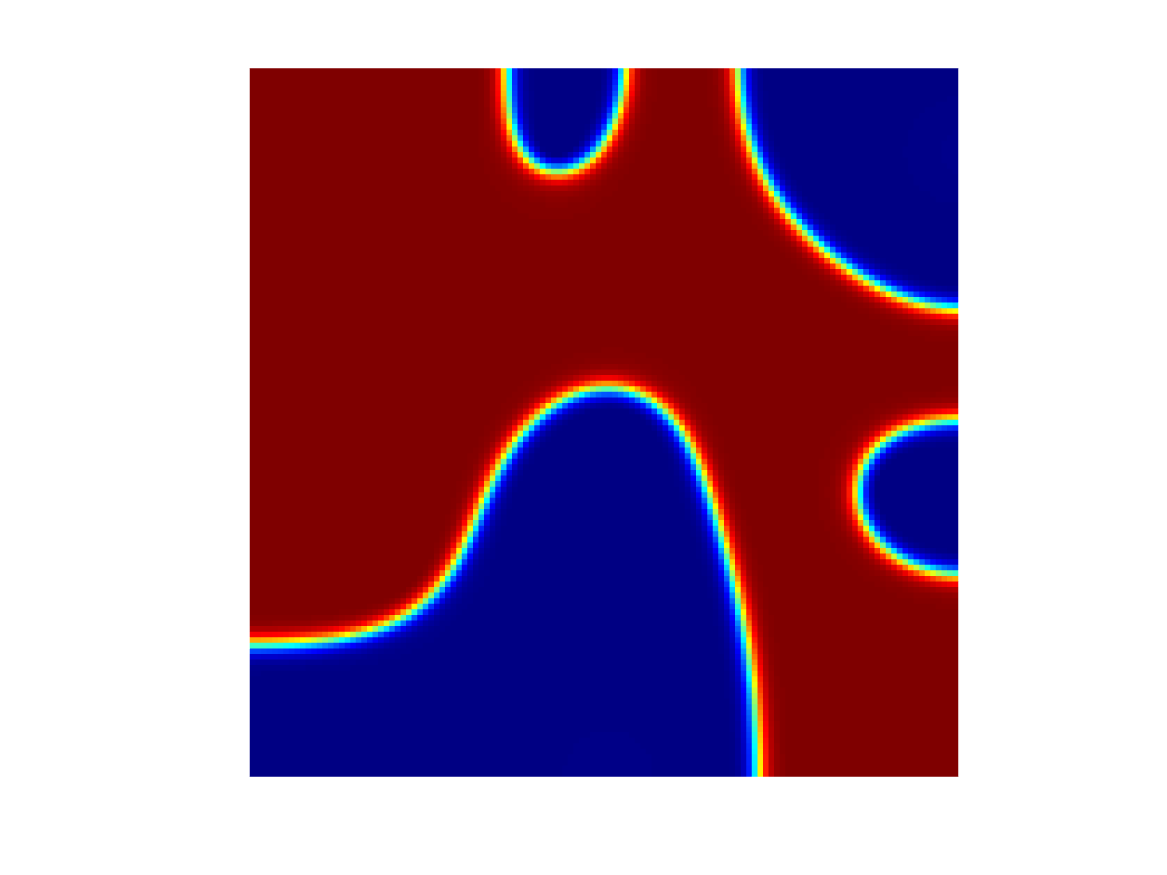}\hspace{-0.95cm}\includegraphics[scale=0.26]{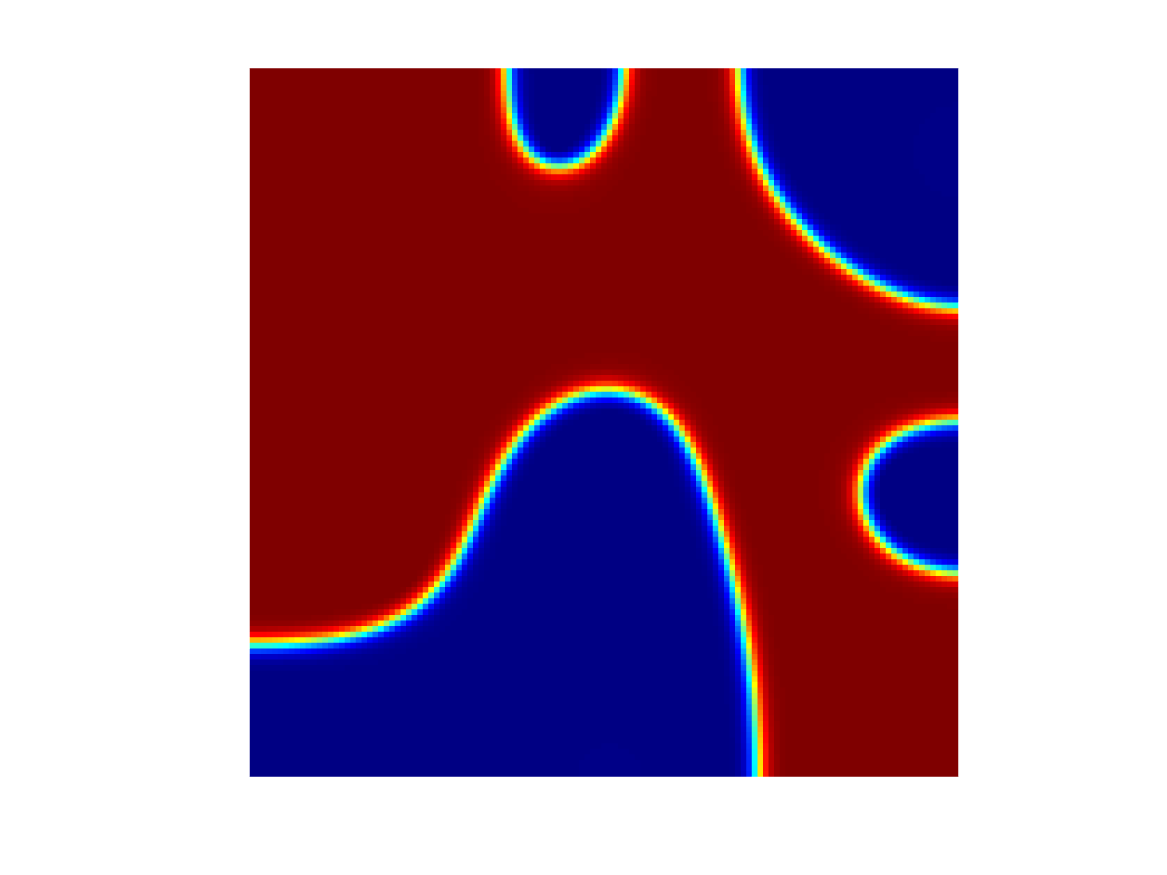}\hspace{-0.95cm}\includegraphics[scale=0.26]{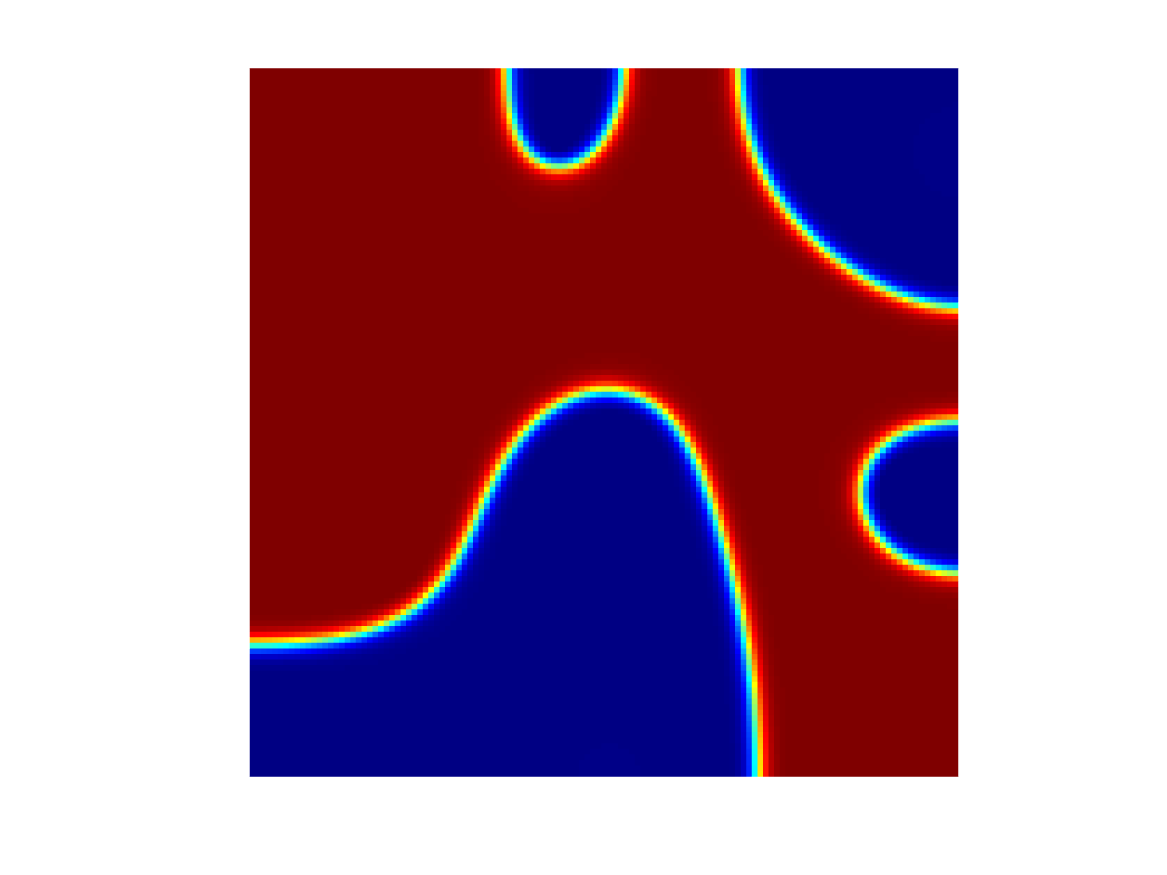}}
\vskip -3mm
\centerline{\includegraphics[scale=0.26]{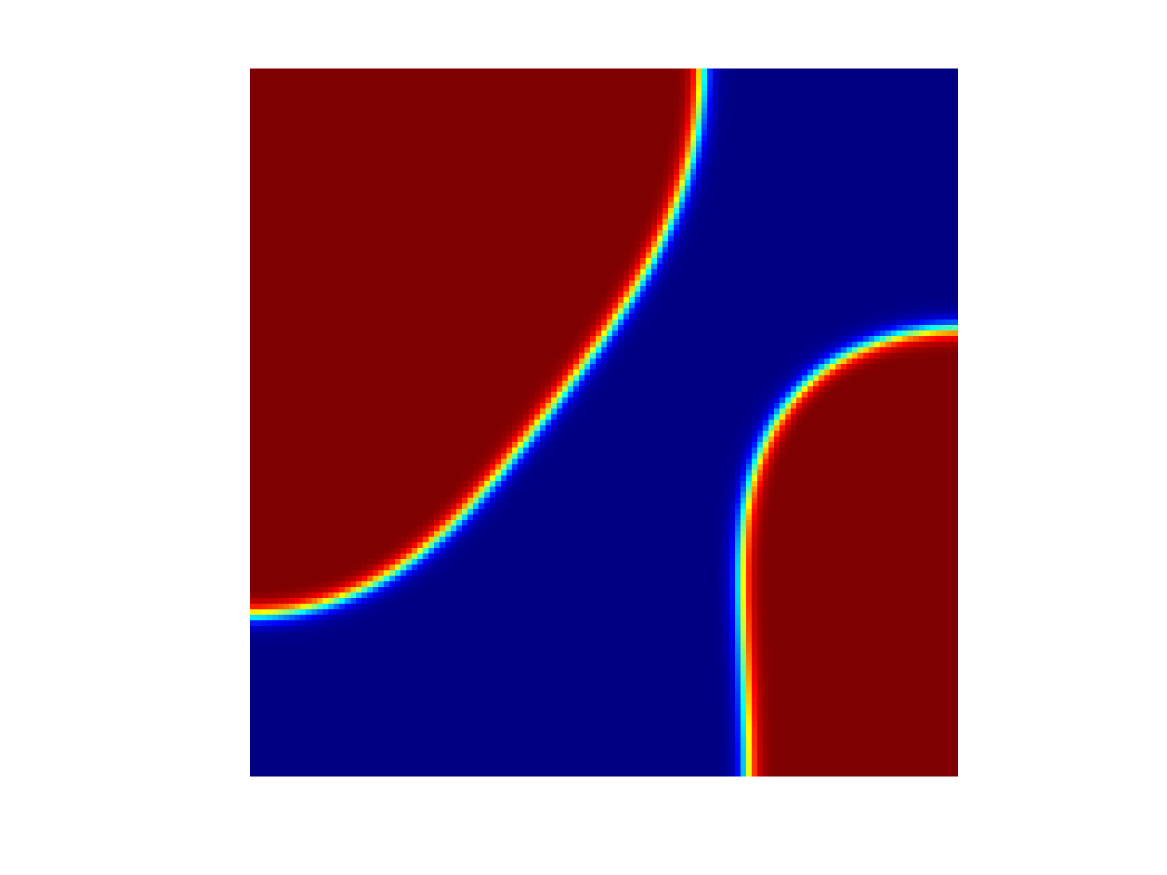}\hspace{-0.95cm}\includegraphics[scale=0.26]{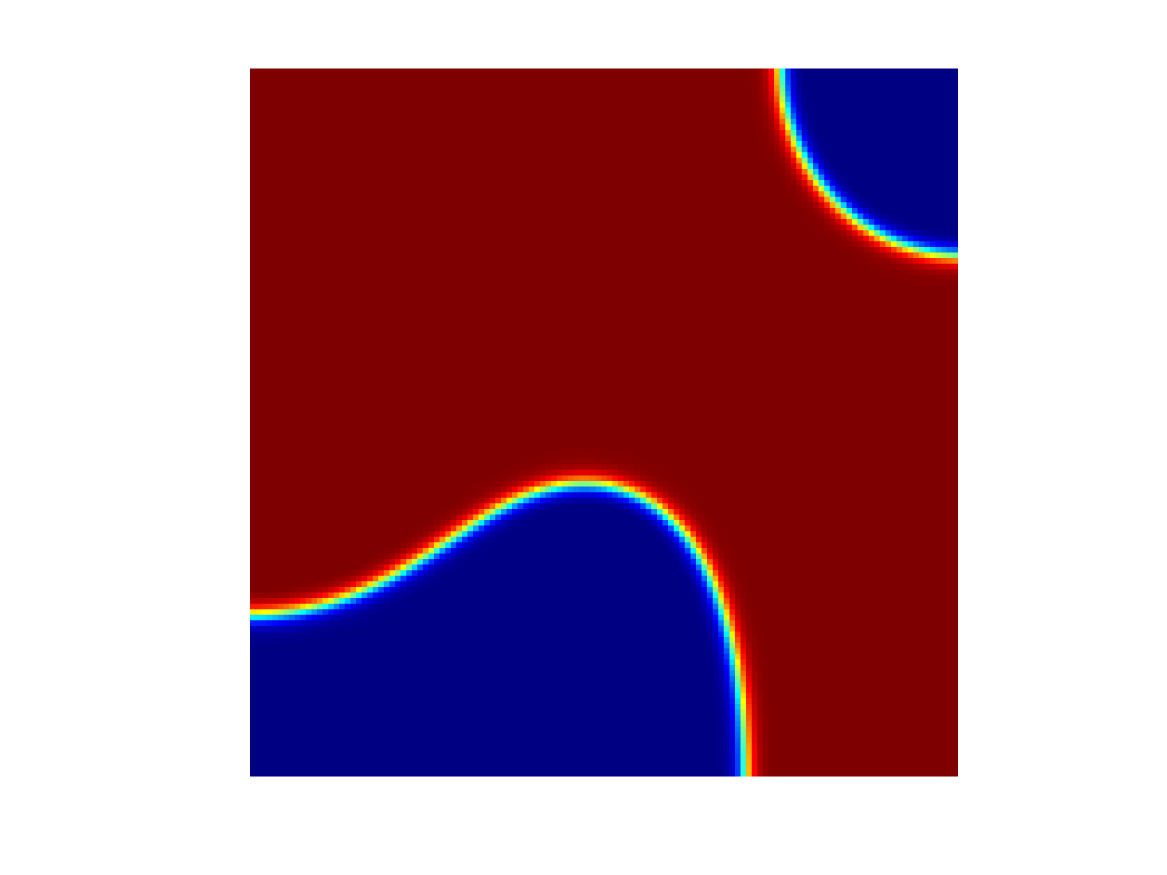}\hspace{-0.95cm}\includegraphics[scale=0.26]{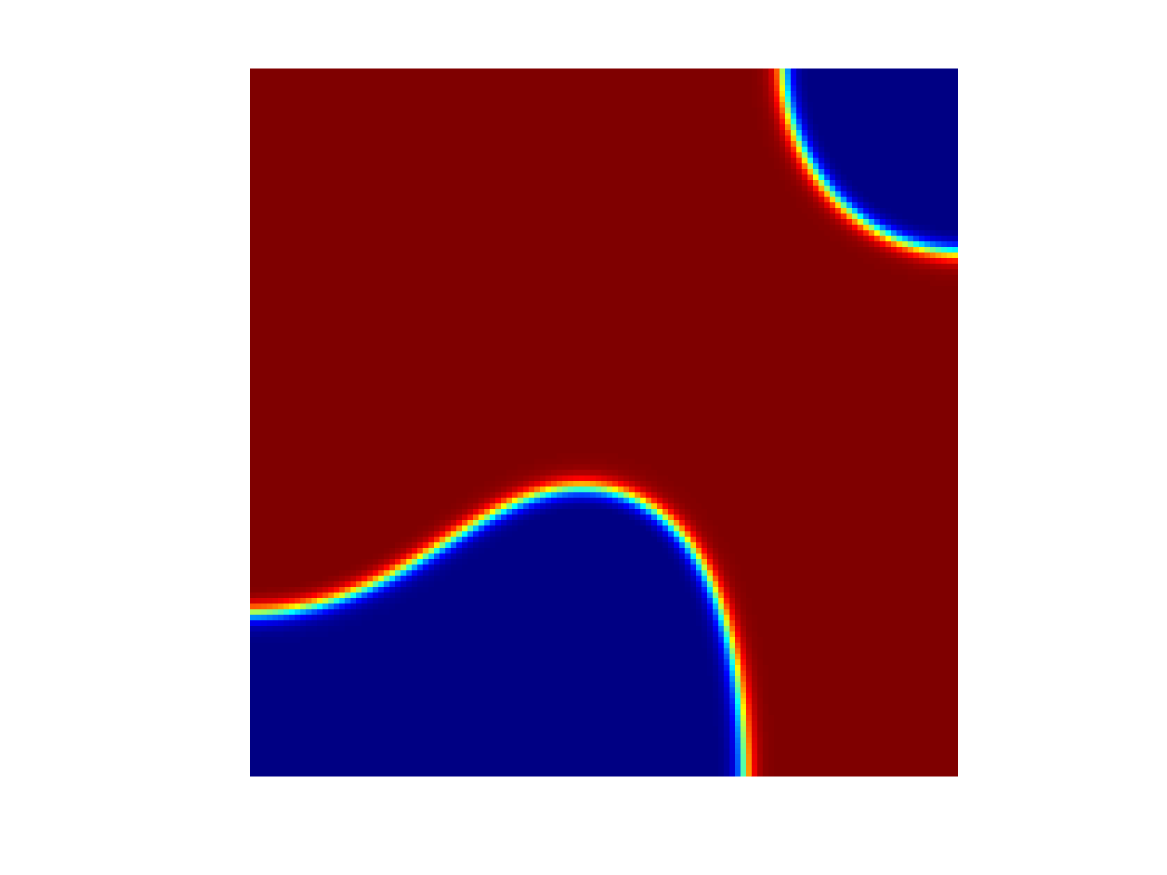}\hspace{-0.95cm}\includegraphics[scale=0.26]{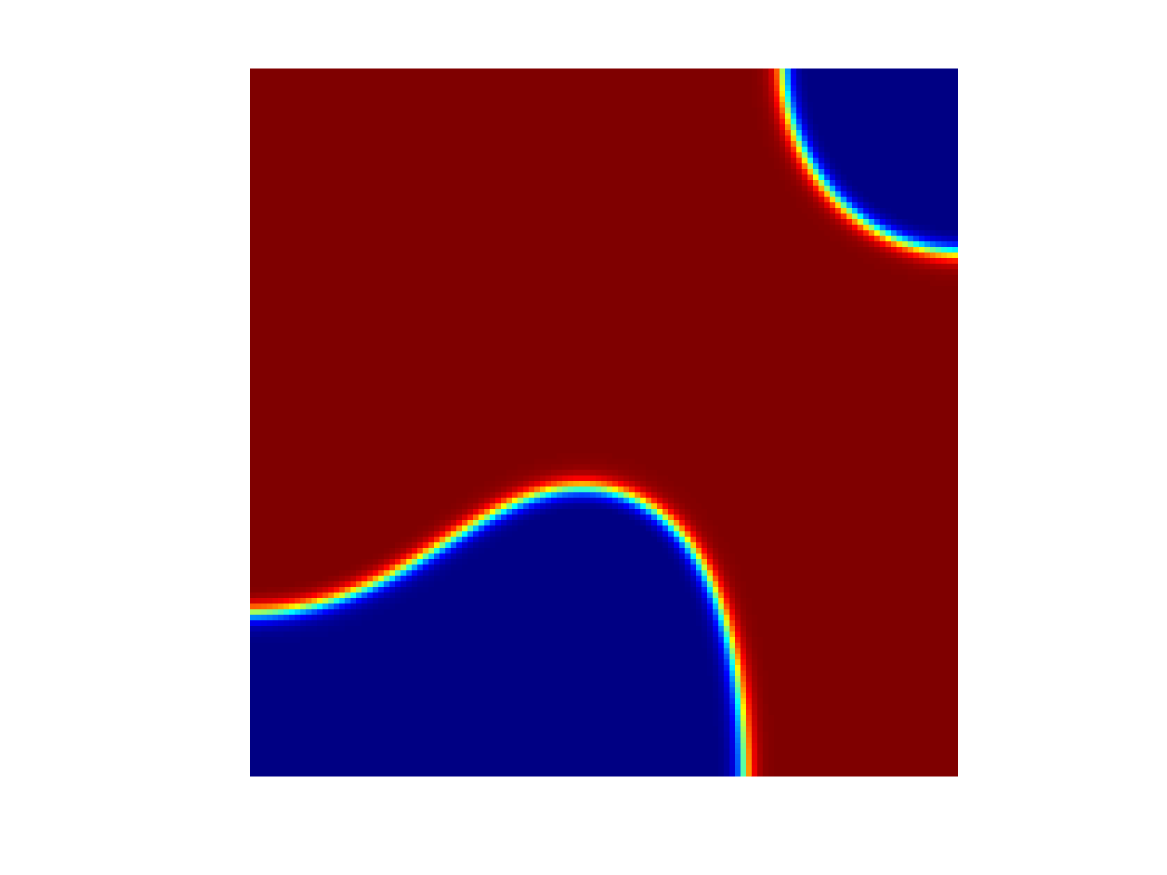}}
\caption{Snapshots of the simulated phase structures around the times $t=10$, $20$, $100$, $200$, and $500$ from top to bottom  produced by the BDF2 scheme \eqref{BDF_2} with  $h=1/128$ and four tested temporal meshes: the uniform time stepping  with fixed large time step size $\tau=0.1$ (the first column); the time adaptive strategy \eqref{adp} with $\tau_{max}=0.1$ (the second column);  the time adaptive strategy \eqref{adp} with $\tau_{max} =\mathcal{G}(1.5)/[S+4L\varepsilon^{2}/h^{2}]$ (the third column); the uniform time stepping with fixed small time step size $\tau=0.01$ (the last column).}
\label{fig2_1}
\end{figure*}

\begin{figure*}[htbp]
\begin{minipage}[t]{0.325\linewidth}
\centerline{\includegraphics[scale=0.28]{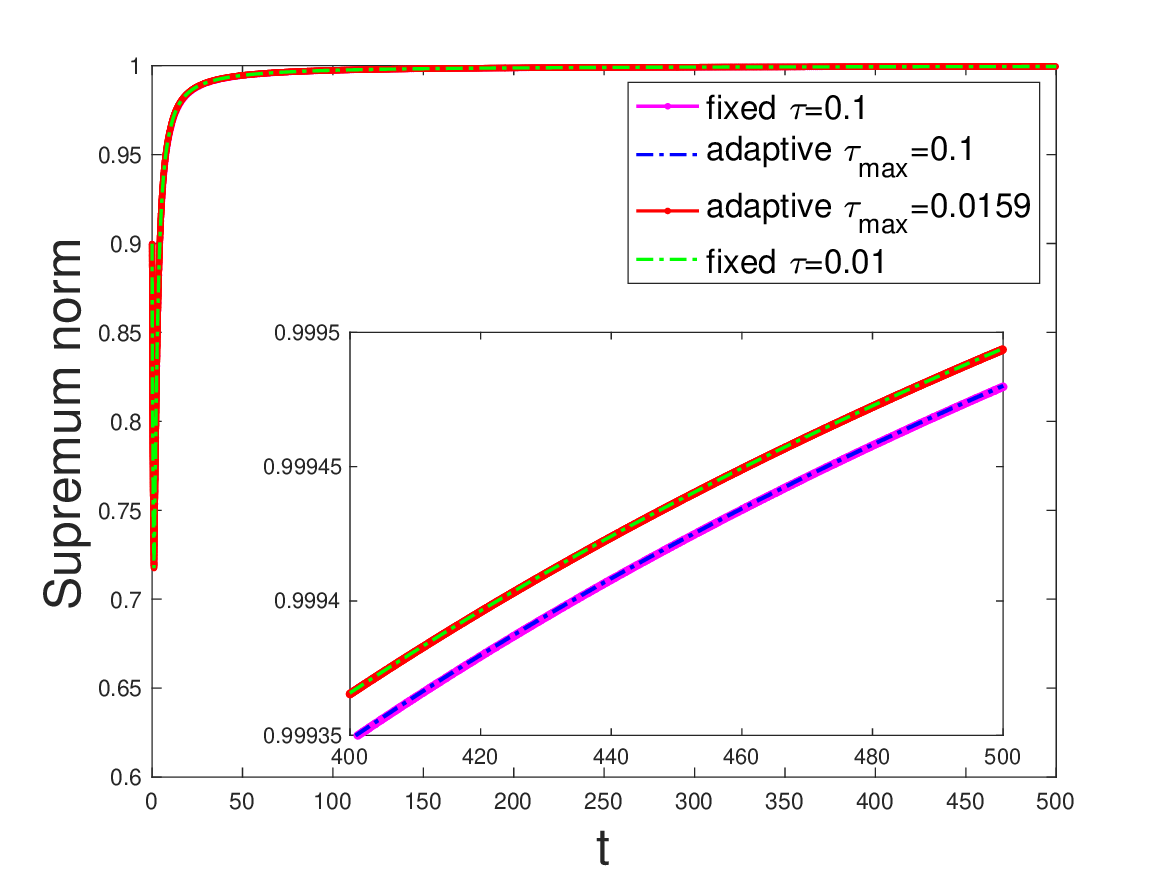}}
\centerline{(a) the supremum norm}
\end{minipage}
\begin{minipage}[t]{0.325\linewidth}
\centerline{\includegraphics[scale=0.28]{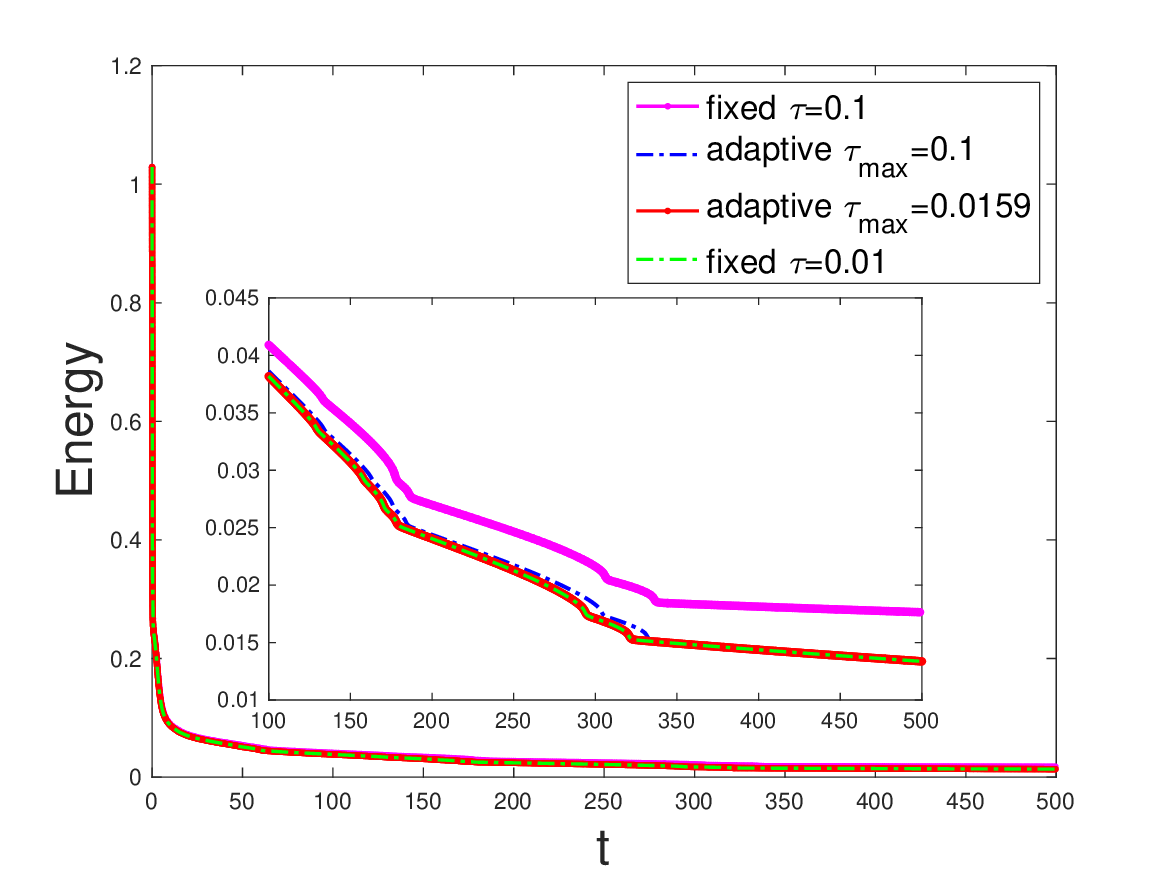}}
\centerline{(b) the energy }
\end{minipage}
\begin{minipage}[t]{0.325\linewidth}
\centerline{\includegraphics[scale=0.28]{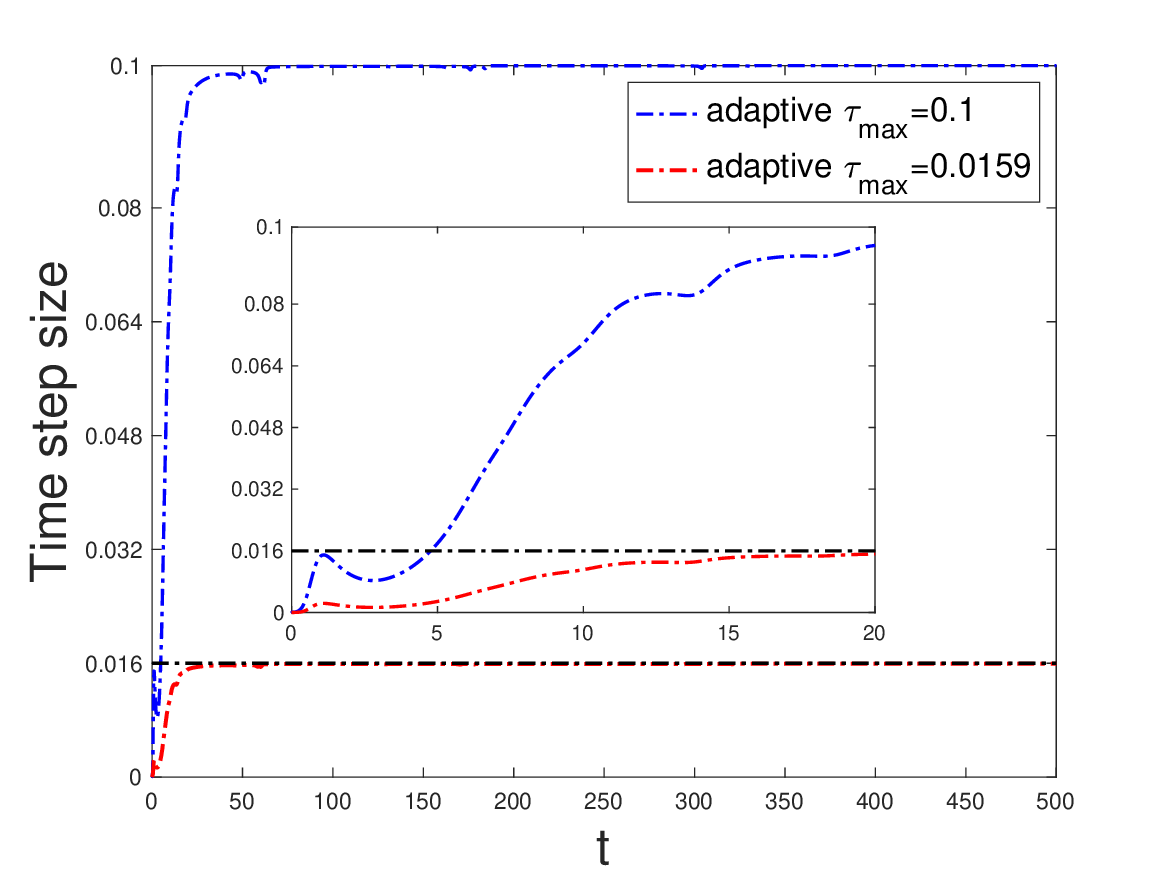}}
\centerline{(c) the time step size }
\end{minipage}
\caption{The evolutions in time of the supremum norm, the energy, and the time step sizes produced by the BDF2 scheme \eqref{BDF_2} with four types of temporal meshes for the grain coarsening  problem.}
\label{fig2_2}
\end{figure*}

The  domain is set to be $\Omega=(-0.5,0.5)^{2}$, and  the initial value configuration is given by a randomly sampled data ranging from $-0.9$ to $0.9$.
 There already exist several efficient time adaptive strategies \cite{LTZ20,Shen17_2,QZT11,STY16,GH11} available to be used together with numerical schemes with variable time steps.
 In this simulation, we will adopt the following robust time adaptive strategy based on the energy
variation proposed in \cite{QZT11}:
\bq\label{adp}
\dps\Dt_{n+1}=\min\Big(\max\big(\Dt_{min},\frac{\Dt_{max}}{\sqrt{1+\alpha |E^{'}(t)|^{2}}}\big),\gamma_{max}\Dt_{n}\Big),
\eq
where $\Dt_{min},\Dt_{max}$ denote the predetermined minimum and maximum time step sizes, $\gamma_{max}\in (0,1+\sqrt{2})$ is the predetermined maximum time step ratio, and $\alpha>0$ is a constant parameter. Such time adaptive strategy will automatically select large time steps when energy decays rapidly and small ones otherwise.
We numerically solve the coarsening dynamics problem using the  BDF2 scheme with four different types of temporal meshes, including the uniform
time stepping with a large step size $\tau=0.1$, two different ones from the time adaptive strategy \eqref{adp}, and the uniform time stepping with  a small step size $\tau=0.01$. For the time adaptive strategy \eqref{adp}, we always set $\gamma_{max}=1.5$, $\alpha=10^{5}$ and $\tau_{min}=10^{-5}$. Also
the predetermined maximum time step sizes are set to be $\tau_{max} =\mathcal{G}(1.5)/[S+4L\varepsilon^{2}/h^{2}]=0.0159$ satisfying the requirement \eqref{est_tau3} and a large one $\tau_{max}=0.1$ for the two tested adaptive temporal meshes, respectively.
 A visual comparison on the numerical solution evolution between these four types of temporal meshes is presented in Figs. \ref{fig2_1} and \ref{fig2_2}.
It is observed that
there is no obvious difference at about $t=10$ for the four tested temporal meshes as shown in the first line of Fig.  \ref{fig2_1},  in which the snapshot of the simulated phase structure with the uniform large time
 step size $\tau=0.1$ only differs from other three temporal grids in a few small details. As shown in Figs. \ref{fig2_1} and \ref{fig2_2}-(b), these minor phase-structure differences at $t=10$ gradually lead to inaccurate solution evolution and energy evolution for the case of the uniform large time
 step size $\tau=0.1$, while the tested two adaptive time strategies still produce correct coarsening pattern which is
 consistent with the numerical results computed by the small time step case $\tau=0.01$.
In Fig. \ref{fig2_2}-(a), we successfully verify the MBP-preserving property of the BDF2 scheme by displaying the evolution of the supremum norm of the numerical solution.
We also note that although both the uniform large time step case $\tau=0.1$ and the adaptive time strategy case with $\tau_{max}=0.1$ don't satisfy the condition \eqref{est_tau3}, they still maintain the MBP-preserving property. It suggests that the constraint \eqref{est_tau3} on the time step size may not be optimal for the proposed BDF2 scheme  \eqref{BDF_2} in term of  preserving the discrete MBP property. 
Furthermore, the evolutions of the energy and the adaptive time step sizes,  plotted in Fig. \ref{fig2_2}-(b)\&(c),
demonstrate  the monotonic energy dissipation and  the efficiency of the  BDF2 scheme with the time adaptive strategy.

\section {Concluding remarks}
In this paper we propose a second-order  BDF scheme with nonuniform time steps for the Allen-Cahn equation with a general mobility.
The MBP preservation of the proposed scheme is successfully established with mild restrictions on the time step sizes and the ratio of adjacent time step sizes.
Moreover, the discrete $H^{1}$ error estimate  and energy stability are rigorously derived for the constant mobility case  and so does the  $L^{\infty}$ error estimate for the general mobility case.
Finally, various  numerical experiments are carried out to validate the theoretical results and demonstrate  the performance of the proposed scheme adopted with a time adaptive strategy.
 It remains interest to further theoretically  explore  the discrete $H^1$ error analysis and energy stability for the general mobility case,
and study the Allen-Cahn equation  with  the logarithmic potential, instead of the double-well potential studied in this paper.
Moreover, we also would like extend the present work to the time-fractional Allen-Cahn equation, in which it is urgently desired to make use of variable-step structure-preserving high-order time stepping schemes to overcome the initial singularity from fractional derivatives.
 In addition, there are two non-constant coefficient Poisson-type equations to be solved at each time step for the model with non-constant mobility in the proposed linear BDF2 scheme. Consequently, it may not be  computationally cheaper and more accurate than a comparable second-order nonlinear scheme with the use of nonlinear multigrid method \cite{WSM10,CWWW19}.
 Thus, it is also an interesting future work to study nonlinear MBP-preserving numerical schemes for
the Allen-Cahn equation with variable mobility.

\bibliographystyle{plain}
\bibliography{ref}
\end{document}